\PassOptionsToPackage{prologue,dvipsnames}{xcolor}
\documentclass{article}

\usepackage{arxiv}

\usepackage[utf8]{inputenc} 
\usepackage[T1]{fontenc}    
\usepackage{hyperref}       
\usepackage{url}            
\usepackage{booktabs}       
\usepackage{amsfonts}       
\usepackage{nicefrac}       
\usepackage{microtype}      
\usepackage{lipsum}		
\usepackage{graphicx}
\usepackage[square,numbers]{natbib}
\usepackage{doi}

\usepackage{amssymb,amsmath,amsthm,graphicx,tikz,caption,subcaption,pgfplots,hyperref,float}
\usepackage{siunitx, comment, multirow}
\usepackage{bm,array,arydshln}
\pgfplotsset{compat=1.18} 
\usepackage[dvipsnames]{xcolor}

\usepackage{tikz-3dplot}
\usetikzlibrary{patterns}
\usetikzlibrary{graphs, arrows.meta}

\usepackage{algorithm}
\usepackage{algpseudocode}

\usepackage{booktabs} 

\usetikzlibrary{colorbrewer}

\DeclareMathOperator{\sfrac}{frac}

\renewcommand{\vec}[1]{\bm{#1}}

\newcommand{\secondranktensor}[1]{\bar{\bar{#1}}}

\theoremstyle{definition}
\newtheorem{mydefinition}{Definition} 

\newtheorem{example}{Example}

\algdef{SE}[DOWHILE]{Do}{doWhile}{\algorithmicdo}[1]{\algorithmicwhile\ #1}%

\algnewcommand\algorithmicforeach{\textbf{for each}}
\algdef{S}[FOR]{ForEach}[1]{\algorithmicforeach\ #1\ \algorithmicdo}

\definecolor{turquoiseblue}{rgb}{0.0, 1.0, 0.94}

\pgfplotscreateplotcyclelist{exotic}{
    teal,every mark/.append style={fill=teal!80!black},mark=*\\
    orange,every mark/.append style={fill=orange!80!black},mark=square*\\
    cyan!60!black,every mark/.append style={fill=cyan!80!black},mark=triangle*\\
    red!70!white,mark=star\\
    lime!80!black,every mark/.append style={fill=lime},mark=diamond*\\
    red,densely dashed,every mark/.append style={solid,fill=red!80!black},mark=*\\
    yellow!60!black,densely dashed,
        every mark/.append style={solid,fill=yellow!80!black},mark=square*\\
    black,every mark/.append style={solid,fill=gray},mark=otimes*\\
    blue,densely dashed,mark=star,every mark/.append style=solid\\
    red,densely dashed,every mark/.append style={solid,fill=red!80!black},mark=diamond*\\
}

\newcommand{\AxisRotator}[1][rotate=0]{%
    \tikz [x=0.25cm,y=0.60cm,line width=.2ex,-stealth,#1] \draw (0,0) arc (-150:150:1 and 1);%
}

\newcommand{\cercle}[4]{
\node[circle,inner sep=0,minimum size={2*#2}](a) at (#1) {};
\draw[black, thick,-stealth] (a.#3) arc (#3:#4:#2);
}

\title{Cell agglomeration strategy for cut cells in eXtended discontinuous Galerkin methods}

\author{Muhammed Toprak}

\author{Matthias Rieckmann}
\author{Florian Kummer}

\newif\ifuniqueAffiliation

\ifuniqueAffiliation 
\author{ \href{https://orcid.org/0000-0000-0000-0000}{\includegraphics[scale=0.06]{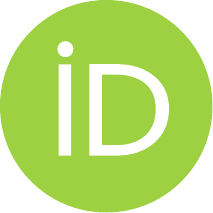}\hspace{1mm}David S.~Hippocampus}\thanks{Use footnote for providing further
		information about author (webpage, alternative
		address)---\emph{not} for acknowledging funding agencies.} \\
	Department of Computer Science\\
	Cranberry-Lemon University\\
	Pittsburgh, PA 15213 \\
	\texttt{hippo@cs.cranberry-lemon.edu} \\
	\And
	\href{https://orcid.org/0000-0000-0000-0000}{\includegraphics[scale=0.06]{orcid.pdf}\hspace{1mm}Elias D.~Striatum} \\
	Department of Electrical Engineering\\
	Mount-Sheikh University\\
	Santa Narimana, Levand \\
	\texttt{stariate@ee.mount-sheikh.edu} \\
}
\else
\usepackage{authblk}

\newbox{\orcid}\sbox{\orcid}{\includegraphics[scale=0.06]{orcid.pdf}} 
\author[1,2]{%
	\href{https://orcid.org/0009-0005-8817-5044}{\usebox{\orcid}\hspace{1mm}Muhammed Toprak\thanks{\texttt{toprak@fdy.tu-darmstadt.de}}}%
}
\author[1,2]{%
	\href{https://orcid.org/0000-0003-2199-9573}{\usebox{\orcid}\hspace{1mm}Matthias Rieckmann}%
}
\author[1,2]{%
	\href{https://orcid.org/0000-0002-2827-7576}{\usebox{\orcid}\hspace{1mm}Florian Kummer}%
}
\affil[1]{Chair of Fluid Dynamics, TU Darmstadt, Otto-Berndt-Str. 2, D-64287, Germany}
\affil[2]{Graduate School of Computational Engineering, TU Darmstadt, Dolivostr. 15, D-64293, Germany}

\fi

\renewcommand{\shorttitle}{Cell agglomeration strategy for cut cells in XDG methods}

\hypersetup{
pdftitle={A template for the arxiv style},
pdfsubject={q-bio.NC, q-bio.QM},
pdfauthor={David S.~Hippocampus, Elias D.~Striatum},
pdfkeywords={First keyword, Second keyword, More},
}

\begin{document}
\maketitle

\begin{abstract}
In this work, a cell agglomeration strategy for the cut cells arising in the extended discontinuous Galerkin (XDG) method is presented. Cut cells are a fundamental aspect of unfitted mesh approaches where complex geometries or interfaces separating sub-domains are embedded into Cartesian background grids to facilitate the mesh generation process. In such methods, arbitrary small cells occur due to the intersections of background cells with embedded geometries and lead to discretization difficulties due to their diminutive sizes. Furthermore, temporal evolutions of these geometries may lead to topological changes across different time steps.
Both of these issues, i.e., small-cut cells and topological changes, can be addressed with a cell agglomeration technique. 
In this work, a comprehensive strategy for the typical issues associated with cell agglomeration in three-dimensional and multiprocessor simulations is provided. The proposed strategy is implemented into the open-source software package BoSSS and tested with 2- and 3-dimensional simulations of immersed boundary flows.
\end{abstract}


\section{Introduction} \label{sec:Intro}




Numerical simulations play a critical role in a wide range of engineering applications, where the generation of a suitable mesh is of significant importance as it serves as a key component in these simulations. 
 However, ensuring the quality and precision of meshes involving complex geometries necessitates special considerations, particularly when irregular boundaries are involved. 
 While conforming the computational mesh to the boundaries may appear the easiest solution at first glance (i.e., body-fitted mesh), it poses significant difficulties, especially in the context of multi-phase or multi-physics systems~\cite{Ingram_Causon_Mingham_2003, Prenter_2023}. 
 The mesh creation process typically consumes a large amount of time, leading to repetitive efforts, and becomes even more challenging when dealing with multiple bodies or dynamic boundaries undergoing topological changes. 
 A further difficulty is introduced with curved meshes, which are present in more recent applications and high-order schemes.
 In such cases, embedding the complex geometry in a static fundamental grid, e.g. a Cartesian grid, provides a straightforward and relatively simple solution (i.e., unfitted mesh). 
 Moreover, the use of a Cartesian grid results in regular cell boundaries for most parts of the mesh, which can significantly reduce the computational requirements for calculating quadratures~\cite{Berger_2017}. In addition, it offers greater accuracy because high-order terms nullify each other due to the symmetry of the mesh structure. 

The notion of embedding bodies into a non-conforming fundamental mesh can be traced back to 1972, with Peskin~\cite{Peskin_1972} proposing the Immersed Boundary Method (IBM) to simulate blood in cardiac flows.
Subsequently, a wide variety of modifications and alternative approaches have emerged under various names such as immersed, embedded, cut cell, unfitted, or extended, each sharing the common objective of efficiently handling numerical simulations without the constraints of a conforming mesh. 
The first applications of such a methodology into Finite Element Methods (FEMs) can be found in the works of Melenk and Babu\u{s}ka~\cite{Melenk_1996, Babuska_1997} or \text{Mo{\"e}s} et al.~\cite{Moes_1999} who introduced extended FEM, so-called XFEM, to simulate crack growth in composite materials. 
Later, the finite cell method (FCM)~\cite{Parvizian_2007, Duster_2008}, CutFEM~\cite{Burman_2015}, and AgFEM~\cite{Badia_2018} can be counted as notable examples of the successor methods that utilize FEMs to solve partial differential equations on a given background grid without mesh alignment.

Within the realm of geometrically unfitted approaches, discontinuous Galerkin (DG) methods have recently garnered special attention due to the flexibility offered by the weak enforcement between the cells, whose initial development can be dated back to 1973 with the work of Reed and Hill~\cite{ReedHill1973}, as well as the works of Babu\u{s}ka~\cite{Babuska_1970,Babuska_Zlamal_1973}, Nitsche~\cite{Nitsche1971}, and Arnold~\cite{Arnold_1980}. 
 Compared to the standard FEMs, in which solutions are defined globally, this facilities local adaptivity by allowing a better treatment for cell refinement and agglomeration through simple modifications of the local shape function, without additional measures for stability~\cite{Prenter_2023}.
 Moreover, DG methods are inherently conservative schemes and require coupling only between the immediate neighbor cells, thus streamlining implementation efforts for data structures and parallelization. 
 Nevertheless, DG methods necessitate a higher number of degrees of freedom (DOF) for the same problem, juxtaposed with their advantageous properties.
 The first extended DG (XDG) method, the so-called unfitted DG method, was introduced by Bastian and Engwer~\cite{Bastian_2009} to simulate fluid dynamics within porous media to account for complex geometries. 
 Since then, its variants have been developed and applied to different fields, including two-phase flows~\cite{Heimann_2013,kummer_extended_2017, Henneaux_2023}, fluid-structure interaction~\cite{Muller_2017, SAYE_2017, Zonca_2018}, acoustics~\cite{Schoeder_2020}, and shock capturing~\cite{Geisenhofer_2019}.

 Another crucial aspect of the numerical simulations is the representation of physical boundaries, a feature relevant to numerous engineering applications.
In, for example, CFD community, the approaches to address multiphase flows (e.g. fluid-structure interaction~\cite{BORAZJANI_2008, MOKBEL_2018}, fluid-fluid~\cite{Hirt_1981, Sussman_1994, Jacqmin_1999, Tryggvason_2001, Olsson_2005} or particle-laden flows~\cite{Uhlmann_2005}) can be broadly categorized into diffuse and sharp interface methods~\cite{Sotiropoulos_2014, Mirjalili2017InterfacecapturingMF}, depending on how they treat the discontinuity between phases. 
Diffuse interface methods entail a mathematical framework to manage abrupt changes (e.g. jumps and kinks) and geometric irregularities occurring at physical boundaries with the help of relaxation operators, creating a smooth interface of finite thickness distributed over multiple grid nodes. 
Nonetheless, numerically resolving the interface thickness in diffuse interface methods proves arduous due to its considerably smaller magnitude~\cite{Ding_2007}.
Conversely, sharp interface methods take a more explicit approach to treating the interface with a theoretical zero thickness, where the behavior of each phase is determined individually. 
Hence, they are often preferable for the calculation of physical phenomena relying on interfacial interactions like phase change, surface tension, and multiphase heat transfer.
Yet, they introduce higher complexity.  
This work represents a sharp interface concept accompanied by the extended DG method that can automatically handle arbitrary complex shapes on a Cartesian background grid.

When an embedded (a.k.a. immersed) geometry or interface intersects the given background grid, it creates so-called \textit{cut cells} which are split into two separate disjoint domains.
 Typically the behavior of those cut cells is defined implicitly, for example, by a level-set function, and can lead to almost arbitrary sizes and shapes since the background mesh is not fitted to the geometry.
 As a consequence, various problems can arise in numerical calculations, which are often referred to as the \textit{small-cut problem}~\cite{Prenter_2023,Burman_2021}. 
 For instance, cut cells may have sizes of several orders of magnitude smaller than the typical elements, which challenges not only the conditioning of discretizations but also the time step restrictions for explicit schemes. 
 Moreover, due to their irregular shapes, these cut cells could introduce additional inaccuracies and complications in integration. Therefore, addressing the small-cut problem is imperative to ensure the stability and accuracy of the numerical solutions.
 
To this date, different solutions have been proposed to overcome the small-cut problem in the literature, such as h-box method~\cite{Berger_2012}, pre-conditioning~\cite{Lehrenfeld_2017, Prenter_2019}, ghost penalty formulation~\cite{Burman_2010, Burman_2010b}, flux distribution~\cite{Chern_1987, Colella_2006}, DoD stabilization~\cite{May_2022}, and cell agglomeration (a.k.a. merging or aggregation)~\cite{Bayyuk_1993, Coirier_1995, Qin_2013}. 
Among these solutions, cell agglomeration stands out as the most straightforward approach~\cite{Berger_2012}, in which small-cut cells are merged with suitable neighbors to form larger cells, providing a convenient solution to the associated problems. In addition, it can also be employed to regulate topological changes stemming from evolving interfaces, as proposed by \cite{Kummer_2018_time}. These changes occur due to the alterations in the domain of phases and pose a conceptual problem in matching discretizations across time steps. Hence, the proposed agglomeration aims to create topologically consistent grid structures by merging the non-matched cells with appropriate cells, thus offering a pragmatic solution to the conceptual and computational difficulties associated with evolving interfaces. In this study, we present a cell agglomeration technique for resolving both the small-cut and the topological inconsistency problems in XDG methods, which can also be applied to a wider array of methodologies utilizing cell agglomeration.






In DG methods, cell agglomeration is performed by simply extending the support of basis functions of a neighbor cell to cover undesired cells, replacing the original polynomial space. 
However, creating appropriate agglomeration mappings for these cells becomes notably challenging in 3-dimensional (3D) space since there is a high degree of neighborship between cells.
Moreover, agglomeration can cause cumbersome problems in sizeable computational simulations such as the formation of agglomeration chains and ineffective information exchange, due to which it is often considered to have drawbacks in implementation~\cite{ May_2022,Engwer_2019}. 
Another criticism of the cell agglomeration approach is the insufficient research available~\cite{Gurkan_2020}.
In the meantime, an efficient parallelizable algorithm for cell agglomeration is crucial for large-scale simulations, which are often performed on computer clusters with multiprocessors to meet the demanding computational requirements.

For these reasons, in this study, we elaborate on the cell agglomeration strategy presented in the previous works~\cite{Muller_2017, Kummer_2021_BoSSS} to provide a comprehensive solution for highly dynamic and parallelizable simulations, which can also deal with 3D meshes. Specifically, our focus is on developing a cell agglomeration strategy that mitigates issues related to cut cells such as agglomeration chains, as well as implementational efforts like inter-processor agglomerations. Therefore, we present a general recipe for cell agglomeration by providing complementary algorithms. 
The proposed strategy is implemented using Message Passing Interface (MPI) into the open-source XDG solver BoSSS and tested with multiprocessor simulations of immersed boundary flows in both 2D and 3D spaces.

\section{E(x)tended discontinuous Galerkin method} \label{sec:XDG}
In this section, we introduce the XDG method for a trivial problem by employing the scalar transport equation to display its discretization framework and omit the details of the variational formulation for the sake of simplicity.
 This equation can be interpreted as the continuity equation in the context of fluid mechanics.
 Its expansion to more general settings, such as the Navier-Stokes equation, can be found in the work of Kummer~\cite{kummer_extended_2017}. 

In the absence of production terms, the transport equation for a scalar parameter, which is denoted by $c$, reads as:
\begin{equation} \label{eq:scalarConservation}
	\frac{\partial c}{\partial t} + \nabla \vec{f}(c) = 0,
\end{equation}
where $\vec{f}$ denotes the flux and $c$ is a function of space and time, i.e. $c=c(\vec{x},t)$.
 The space vector $\vec{x}$ is defined in spatial domain $\Omega \subset \mathbb{R}^D$ ($D \in \{ 2, 3 \}$), whereas the time $t$ is a non-negative real number, i.e. $t \in \mathbb{R}^{+}_0$.
 The domain $\Omega$ is discretized into the background grid $\Omega_h$, which is characterized by the length scale of the coarsest background cell, denoted as $h$.
 $\Omega_h$ is formed by a collection of non-overlapping (i.e. $\int_{K_i \cap K_j} 1 dV = \emptyset \text{ for } i \neq j $) and regular-shaped cells as:
	\begin{equation}
		\Omega_h =  \scalebox{1.5}{ $\mathbin{\mathaccent\cdot\cup}$}_{i} \, K_i.
	\end{equation}
 Subsequently, the computational mesh $\mathfrak{K}_h$ is defined as the set of all cells, $\mathfrak{K}_h = \{K_1,...,K_N\}$.
 The set of all the corresponding edges is defined as $\Gamma := \scalebox{1.4}{ $\mathbin{\mathaccent\cdot\cup}$}_{i} \, \partial K_i$ and it is divided into three subsets: $\Gamma = \Gamma_\mathrm{int}(t) \cup \Gamma_\mathrm{D} \cup \Gamma_\mathrm{N}$. The internal edges are defined as $\Gamma_\mathrm{int} := \Gamma \setminus \partial \Omega$, whilst the Dirichlet and Neumann boundary conditions are denoted by $\Gamma_\mathrm{D}$ and $\Gamma_\mathrm{N}$, respectively. 

The broken polynomial space (a.k.a. DG space) with a total degree $p$ is defined as:
\begin{equation} \label{eq:DGspace}
	\mathbb{P}_p(\mathfrak{K}_h) := \{ \phi \in L^2(\Omega); 
	\forall K \in \mathfrak{K}_h : \phi|_K \text{ is a polynomial} \nonumber \\ \text{ and deg}(\phi|_K) \leq p \}.
\end{equation}
The so-called weak formulation for the single cell $K_i$ is then obtained from Eq.~\ref{eq:scalarConservation} via multiplication with the test function $\phi_{i,m}$ and performing integration by parts:
\begin{equation} \label{eq:weakFormulation}
	\int_{K_i} \frac{\partial c} {\partial t} \phi_{i,m} dV + \int_{\partial K_i} (\vec{f}(c) \cdot \vec{n}_{\partial K_i}) \phi_{i,m} dA - \int_{K_i} \vec{f}(c) \cdot \nabla \phi_{i,m} dV = 0,
\end{equation}
where $\vec{n}_{\partial K_i}$ denotes the outward unit normal vector.
Using the Galerkin ansatz, the test and the trial functions are selected to be identical, indicated by subscripts $m$ and $n$, respectively. As a result, a scalar field \(c \in \mathbb{P}_p(\mathfrak{K}_h)\) is represented as a linear combination of the same polynomial basis as:
\begin{equation}
	c_i(\vec{x},t) \approx \sum_{n} \phi_{i,n}(\vec{x}) \tilde{c}_{i,n}(t) = \phi_{i,-} (\vec{x}) \cdot \tilde{c}_{i,-}(t)= \underline{\phi}_{i} (\vec{x}) \cdot \underline{\tilde{c}_{i}}(t),	
\end{equation}
where \(\tilde{c}_{i,n}\) is a DG coefficient representing a degree of freedom (DOF), and \(\phi_{i,n}\) is a basis function of \(\mathbb{P}_p(\mathfrak{K}_h)\) with \(\text{supp}(\phi_{i,n}) = K_i\). 
Thus, a cell-local basis can be expressed as a row vector of the basis functions with $\underline{\phi}_i = \left( \phi_{i,1}, ..., \phi_{i,M} \right)$. Furthermore, Equation~\ref{eq:weakFormulation} can be simplified into temporal and spatial components in the following semi-discrete form:
\begin{equation}
	\label{eq:simplifiedWeak}
	\mathcal{M}_i \int_{K_i} \frac{\partial \tilde{c}} {\partial t} dV + \vec{F}_i(c)  = 0,
\end{equation}
where $\mathcal{M}_i$ is the cell-local mass matrix, which is symmetric with dimensions of $M \times M$, and can be obtained by:
\begin{equation}
	\label{eq:MassMatrix}
	\mathcal{M}_i = \int_{K_i} \underline{\phi}_{i}^T \cdot \underline{\phi}_{i} d V,
\end{equation}
meanwhile, $\vec{F}_i(c)$ encompasses the contributions from the volume and surface integrals originating from the flux term. 
In DG discretizations, the inner and outer values of a scalar parameter $c$ are described as:
\begin{align}
		c^{-}(\vec{x}) &:= \lim_{\xi \to 0} \; \; c(\vec{x} - \xi \vec{n}_\Gamma)  \text{ for} \; \, \vec{x} \in \Gamma , \\
		c^{+}(\vec{x}) &:= \lim_{\xi \to 0} \; \; c(\vec{x} + \xi \vec{n}_\Gamma)  \text{ for} \; \, \vec{x} \in \Gamma_\mathrm{int},
\end{align}
where $\vec{n}_\Gamma$ is the outward unit normal vector of the edge.
Jump and averaging operators are respectively introduced to describe the variation of the parameters along cell edges and material interfaces as:
\begin{align}
	[[c]]:=&  
	\begin{cases}
		c^{-} &\text{ on } \; \, \; \, \; \, \partial \Omega \\
		c^{-} - c^{+} &\text{ on } \; \, \; \, \; \, \Gamma_\mathrm{int}
	\end{cases}
	, \\
	\{\{c\}\}:=&
	\begin{cases}
		c^{-} &\text{ on } \; \, \partial \Omega \\
		\frac{1}{2} (c^{-} + c^{+}) &\text{ on } \; \, \Gamma_\mathrm{int}
	\end{cases}
	.
\end{align} 

Furthermore, the entire domain is divided into two time-dependent disjoint subdomains $\mathfrak{A}$ and $\mathfrak{B}$, along with their interface $\mathfrak{I}(t)$:
\begin{equation}
	\Omega = \mathfrak{A}(t) \mathbin{\mathaccent\cdot\cup} \mathfrak{I}(t) \mathbin{\mathaccent\cdot\cup} \mathfrak{B}(t).
\end{equation}
The subdomains $\mathfrak{A}$ and $\mathfrak{B}$ can be interpreted as the species or phases (e.g. water-oil or solid-fluid) in the context of multiphase flows.
 The boundary of the total domain $\Omega$ is denoted by $\partial \Omega$, while both the boundary $\partial \Omega$ and the interface $\mathfrak{I}$ are $D-1$ dimensional. 

 In this work, the behavior of the subdomains and the interface are controlled by a sufficiently smooth level-set function $\psi$ as:
 \begin{align}
	\mathfrak{I}(t) &:= \{\vec{x} \in \Omega \mid \psi(\vec{x},t) = 0\}, \\
	\mathfrak{A}(t) &:= \{\vec{x} \in \Omega \mid \psi(\vec{x},t) < 0\}, \\
	\mathfrak{B}(t) &:= \{\vec{x} \in \Omega \mid \psi(\vec{x},t) > 0\},
\end{align}
while other forms of interface representation (e.g. volume of fluid method or CAD data) are also applicable to the presented agglomeration strategy.

\begin{figure}[t!]
	\centering
	\begin{tikzpicture}

    \pgfmathsetmacro{\xSize}{3}
    \pgfmathsetmacro{\ySize}{3}


  \draw[domain=-1.5:1.5, samples=100, cyan, thick] plot ({3-((\x)^3)*0.8}, {(\x+1.5)});
  \fill[cyan!10] (9, 3) -- (6, 3) -- plot[domain=1.5:-1.5, samples=100] ({3-((\x)^3)*0.8}, {(\x+1.5)}) -- (6,0) -- (9,0) -- cycle;
  
  \draw[thick] (-\xSize,0) rectangle (\xSize,\ySize);  
  \draw[thick] (0,0) rectangle (\xSize,\ySize);
  \draw[thick] (\xSize,0) rectangle (2*\xSize,\ySize);
  \draw[thick] (2*\xSize,0) rectangle (3*\xSize,\ySize);

  \node[text=cyan] at (\xSize*11/12,\ySize*6.5/12) {$\mathfrak{I}$};
  
  \node at (\xSize*-11/12,\ySize*0.1) {$\mathfrak{A}$};
  \node[text=cyan!90] at (\xSize*35/12,\ySize*0.1) {$\mathfrak{B}$};

  \node at (-\xSize/2,\ySize/2) {$K_{i-1,\mathfrak{A}}$};    
  \node at (\xSize/2,\ySize/2) {$K_{i,\mathfrak{A}}$};
  \draw (\xSize*5/12,\ySize*11.5/12) -- (\xSize*3/12,\ySize*10/12) node[left,below=0pt] {$K_{i,\mathfrak{B}}$};
  \node at (5*\xSize/2,\ySize/2) {$K_{i+2,\mathfrak{B}}$};


   \draw (\xSize*19/12,\ySize*0.5/12) -- (\xSize*21/12,\ySize*2/12)node[right,above=-2pt] {$K_{i+1,\mathfrak{A}}$};
  
  \node at (3/2*\xSize,\ySize/2) {$K_{i+1,\mathfrak{B}}$};
  
\draw [decorate,decoration={brace,amplitude=5pt,mirror}] (0*\xSize,-0.1) -- (2*\xSize,-0.1) node[midway,below=5pt] {Cut cells};

\draw [dashed, ->] (\xSize*11/12,0.95*\ySize) --  (\xSize*16/12,0.95*\ySize)  node[midway,below] {\scriptsize $Agg_\mathfrak{B}$};
  \fill (\xSize*11/12,0.95*\ySize) circle (1pt);

\draw [dashed, ->]  (\xSize*13/12,0.1*\ySize) -- (\xSize*8/12,0.1*\ySize)node[midway,above] {\scriptsize $Agg_\mathfrak{A}$};
  \fill  (\xSize*13/12,0.1*\ySize) circle (1pt);

\end{tikzpicture}
	\caption{Cut cell agglomeration on the XDG space for an arbitrary interface $\mathfrak{I}$. The cut cells are indicated by $K_{i,\mathfrak{A}}$, $K_{i,\mathfrak{B}}$, $K_{i+1,\mathfrak{A}}$, $K_{i,1\mathfrak{B}}$. The small-cut cells are agglomerated to neighbor elements as $K_{i+1,\mathfrak{A}} \to K_{i,\mathfrak{A}}$ and $K_{i,\mathfrak{B}} \to K_{i+1,\mathfrak{B}}$. \label{fig:CutCell}}
\end{figure}
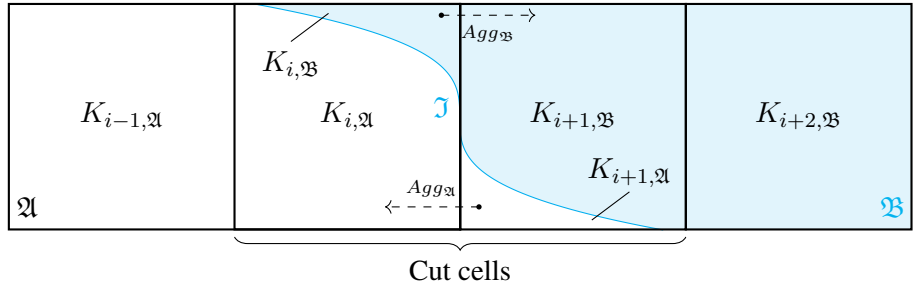

Consequently, the original background cells are subdivided based on their intersection with species $\mathfrak{s} \in \{ \mathfrak{A}, \mathfrak{B}\} $ as:
\begin{equation}
	K_{i, \mathfrak{s}} := K_i \cap \mathfrak{s}(t). 
\end{equation}
In this context, we denote $K_{i, \mathfrak{s}}$ as cells exclusively associated with the species or phases. Hence, they are collectively referred to as \textit{phase cells} to describe them within their respective domain. 
 When a background cell is cut by the interface $\mathfrak{I}$, i.e. $\oint_{K_i \cap \mathfrak{I}(t)} 1 dS > 0$, it produces smaller disjointed phase cells by dividing the background cell into two, as illustrated in Fig.~\ref{fig:CutCell}. These cells, which constitute a special subset of the phase cells, are referred to as \textit{cut cells} and occupy only a portion of the background grid. The time-dependent extended mesh $\mathfrak{K}_h^X(t)$ that collects all the phase cells is defined as:
\begin{equation}
	\mathfrak{K}_h^X(t) := \{ K_{1,\mathfrak{A}},K_{1,\mathfrak{B}}, ... , K_{N,\mathfrak{A}},K_{N,\mathfrak{B}} \}. 
\end{equation}
Accordingly, the DG space in Eq.~\eqref{eq:DGspace} is modified to form the XDG space as:
\begin{align}
	\mathbb{P}_p^X(\mathfrak{K}_h,t) := \mathbb{P}_p(\mathfrak{K}_h^X(t)) = \{ \phi \in L^2(\Omega); 
	\forall K \in \mathfrak{K}_h^X(t) : \phi |_{K \cap \mathfrak{A}}, \, \phi |_{K \cap \mathfrak{B}} \text{ are polynomial}, \nonumber \\ \text{deg}(\phi |_{K \cap \mathfrak{A}}), \, \text{deg}(\phi|_{K \cap \mathfrak{B}}) \leq p \}.
\end{align}
For the integration in cut cells, the method proposed by Saye~\cite{Saye_2015} is employed.
The parameters of the equations are piece-wisely defined in subdomains $\mathfrak{A}$ and $\mathfrak{B}$, which gives for the scalar parameter $c$:
\begin{equation}
	c(\vec{x},t) =\begin{cases} 
		c_\mathfrak{A}, \text{ for} \; \, \vec{x} \in \mathfrak{A}(t) \\
		c_\mathfrak{B}, \text{ for} \; \, \vec{x} \in \mathfrak{B}(t).
	\end{cases}
\end{equation}
Hence, the jump between the bulk phases (i.e. $\Omega \setminus \mathfrak{I}$) at interface $\mathfrak{I}$ reads as	$[[c]] = (c_\mathfrak{B} - c_\mathfrak{A}) \cdot \vec{n_\mathfrak{I}}$ for the interface normal vector $\vec{n_\mathfrak{I}}$, which points from $\mathfrak{A}$ to $\mathfrak{B}$.


\section{Cell agglomeration} \label{sec:agg}
The cell agglomeration strategy is explained in several sections. Initially, the theoretical background along with relevant definitions for the agglomerated XDG space is presented. Then, the algebra of agglomeration and its practical applications are deliberated upon. Lastly, details regarding the implementation and the utilized algorithms are provided.

\subsection{Graph description of cell agglomeration}
In this section, we formally introduce essential definitions for describing an agglomerated space as well as an agglomeration mapping on a numerical mesh with respect to graph theory.
\begin{figure*}[ht!]
	\centering
	\resizebox{0.8\textwidth}{!}{
	\begin{tikzpicture}
    \draw[step=1,thick,dashed] (0,2) grid (13,6);    
    \begin{scope}[every node/.style={circle, thick, draw, minimum size=1.5em}]
        \node[fill=gray!0] (A) at (0.5,3.5) {};
        \node[fill=BrickRed!30] (B) at (1.5,3.5) {};
        \node[fill=NavyBlue!30] (C) at (3.5,4.5) {};
        \node[fill=NavyBlue!30] (D) at (4.5,5.5) {};
        \node[fill=NavyBlue!30] (E) at (4.5,4.5) {};
        \node[fill=NavyBlue!30] (F) at (3.5,3.5) {};
        \node[fill=gray!0] (G) at (4.5,3.5) {};
        \node[fill=NavyBlue!30] (H) at (5.5,3.5) {};
        \node[fill=gray!0] (I) at (5.5,2.5) {};
        \node[fill=BrickRed!30] (J) at (4.5,2.5) {};
        \node[fill=gray!0] (A2) at (7.5,3.5) {};
        \node[fill=BrickRed!30] (A3) at (8.5,3.5) {};
        \node[fill=gray!0] (A4) at (9.5,3.5) {};
        \node[fill=gray!0] (A5) at (8.5,4.5) {};
        \node[fill=gray!0] (A6) at (8.5,2.5) {};
        \node[fill=gray!0] (B2) at (11.5,3.5) {};
        \node[fill=NavyBlue!30] (B3) at (12.5,3.5) {};
        \node[fill=NavyBlue!30] (B4) at (12.5,2.5) {};
        \node[fill=NavyBlue!30] (B5) at (11.5,2.5) {};
        \node[fill=BrickRed!30] (B6) at (11.5,4.5) {};
    \end{scope}
    
    \begin{scope}[>={Stealth[NavyBlue]},
                  every edge/.style={draw=black, thick}]
        \path [->] (C) edge[NavyBlue] (E);
        \path [->] (D) edge[NavyBlue] (E);
        \path [->] (E) edge[NavyBlue] (G);
        \path [->] (F) edge[NavyBlue] (G);
        \path [->] (G) edge[NavyBlue] (J);
        \path [->] (H) edge[NavyBlue] (I);
        \path [->] (B3) edge[NavyBlue] (B2);
        \path [->] (B4) edge[NavyBlue] (B3);
        \path [->] (B5) edge[NavyBlue, bend left] (B6);
    \end{scope}
    
    \begin{scope}[>={Stealth[black]},
                  every edge/.style={draw=black, thick}]
        \path [->] (A) edge (B);
        \path [->] (I) edge (J);
        \path [->] (G) edge (J);
        \path [->] (A2) edge (A3);
        \path [->] (A5) edge (A3);
        \path [->] (A6) edge (A3);
        \path [->] (A5) edge (A3);
        \path [->] (A4) edge (A3);
        \path [->] (B2) edge (B6);
    \end{scope}
    \end{tikzpicture}
	}
	\caption{A sample illustration of agglomeration mapping on a regular mesh. Each tree represents an agglomeration group, with the first tree also representing an agglomeration pair. Final targets are displayed in red, while chain agglomeration edges and their source cells are displayed in blue. The black edges indicate the direct agglomeration pairs and the white vertices indicate their source cells.    \label{fig:aggGraph}}

\end{figure*}
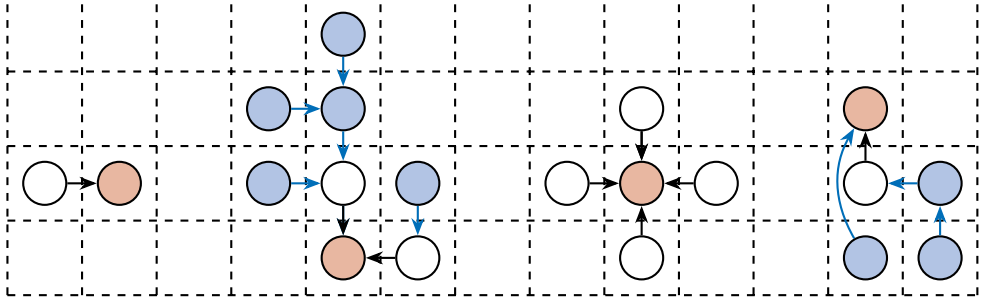
\mydefinition [Agglomeration group and mapping]{With respect to some mesh $\mathfrak{K}$, some set $\mathfrak{A}_\mathrm{grp} \in \mathfrak{K} \times \mathfrak{K}$ is an agglomeration group, if the graph $(\mathfrak{K}_\mathrm{grp},\mathfrak{A}_\mathrm{grp})$, with vertices $\mathfrak{K}_\mathrm{grp}$ and edges $\mathfrak{A}_\mathrm{grp}$, is a directed tree, where each edge directs towards the final vertex (i.e., a directed rooted tree). The forest of agglomeration groups in a mesh is also called agglomeration mapping and is denoted by $\mathfrak{A}_\mathrm{map}$.
}
%

\mydefinition[Agglomeration source and target]{For some pair $(K_\mathrm{src},K_\mathrm{tar}) \in \mathfrak{A}_\mathrm{map}$, the first entry $K_\mathrm{src}$ is called agglomeration source, while the second entry $K_\mathrm{tar}$ is called agglomeration target. }

\mydefinition[Direct and chain agglomeration] {Some pair of cells $(K_\mathrm{src},K_\mathrm{tar}) \in \mathfrak{A}_\mathrm{grp}$ is called a direct agglomeration if $K_\mathrm{tar}$ is the final target (i.e., the root vertex) in $\mathfrak{A}_\mathrm{grp}$ and they share a common boundary (i.e., $\small \oint_{\partial K_\mathrm{src} \cap \partial K_\mathrm{tar}} 1 dS > 0$  ). Otherwise, $(K_\mathrm{src},K_\mathrm{tar})$ is called a chain agglomeration.
}
\mydefinition[Agglomerated cell] {Let $\mathfrak{A}_\mathrm{grp} \subset \mathfrak{A}_\mathrm{map}$ be an agglomeration group, i.e., there is no other edge $(K_\mathrm{src}, K_\mathrm{tar}) \in \mathfrak{A}_\mathrm{grp}$, which has a connection to the edges in $\mathfrak{A}_\mathrm{map}$. Let cell $K_\mathrm{agg}$ be the final target of the agglomeration group. 
Then, the agglomerated cell for $\mathfrak{A}_\mathrm{grp}$ is the union of all cells in it, i.e.,
	\begin{equation}
		K_{\mathrm{agg}} := \bigcup_{ (K_\mathrm{src},K_\mathrm{tar}) \in \mathfrak{A}_\mathrm{grp}  } (K_\mathrm{src} \cup K_\mathrm{tar}).
	\end{equation}
}

\mydefinition[The agglomerated mesh and space]{The agglomerated mesh $\mathfrak{K}_h^{\mathrm{agg}}$ is defined as the set of all the ordinary cells, which are not part of an agglomerated cell, and the agglomerated cells. Hence, the agglomerated DG space is defined as a subspace of the original space:
\begin{equation} \label{eq:subSpace}
	\mathbb{P}_p^{\mathrm{agg}}(\mathfrak{K}_h) := \mathbb{P}_p(\mathfrak{K}_h^{\mathrm{agg}}).
\end{equation}
}
Note: 
\begin{itemize}
	\item[1)] By these definitions, an agglomeration mapping may not contain cycles because it solely consists of trees.
	\item[2)] Chain agglomerations are typically the result of pairs concatenating with direct agglomerations and may be replaced by equivalent pairs. (See §\ref{sec:aggLevels})
	\item[3)] Due to the directing of the graph, the root (a.k.a. final target) of all agglomeration groups is unique.
\end{itemize}


\subsection{Agglomeration algebra} \label{sec:aggAlgebra}
\begin{figure}[ht!]
	\centering
	\resizebox{0.45\textwidth}{!}{ 
		\begin{subfigure}[b]{0.6\textwidth}
			\centering
			\begin{tikzpicture}
    \begin{axis}[
      xlabel={$x$},
      axis lines=none,
      xmin=-0.1, xmax=3.5,
      ymin=-0.5, ymax=2,
      xtick={},
      ytick={},
      ticks=none,
      grid=none,
    ]
    
    \addplot[BrickRed, ultra thick, densely dotted, domain=0:1] {1};
    \addplot[NavyBlue, ultra thick, densely dotted, domain=0:1] {x/2 + 0.25};

    \addplot[BrickRed, ultra thick, domain=1:2] {1}  node[pos=-0.5,above] {$\phi_{agg,1}$};
    \addplot[NavyBlue, ultra thick,domain=1:2] {x/2 + 0.25}  node[pos=-0.5,above=2pt] {$\phi_{agg,2}$};

    \addplot[BrickRed, ultra thick, densely dotted, domain=2:3] {1};
    \addplot[NavyBlue, ultra thick, densely dotted, domain=2:3] {x/2 + 0.25};

    \draw[loosely dashed, black] (axis cs: 0, 0) -- (axis cs: 0, 2);
    \draw[loosely dashed, black] (axis cs: 1, 0) -- (axis cs: 1, 2);
    \draw[loosely dashed, black] (axis cs: 2, 0) -- (axis cs: 2, 2);
    \draw[loosely dashed, black] (axis cs: 3, 0) -- (axis cs: 3, 2);


    \draw[line width=2pt] (0,0) -- (3,0) node[pos=0.5, below, very thick, font=\large] {$K_\mathrm{agg}$}; 
    \draw[line width=2pt] (0,-1/10) -- (0,1/10); 
    \draw[line width=2pt] (3,-1/10) -- (3,1/10); 

    \end{axis}
  \end{tikzpicture}
		\end{subfigure}
	}
		\hfill
		\resizebox{0.45\textwidth}{!}{
		\begin{subfigure}[b]{0.6\textwidth}
			\centering
			\begin{tikzpicture}
    \begin{axis}[
      xlabel={$x$},
      ylabel={$\phi(x)$},
      axis lines=none,
      xmin=-0.1, xmax=3.1,
      ymin=-0.5, ymax=2,
      xtick={},
      ytick={},
      grid=none,
    ]
    \draw[line width=2pt] (0,0) -- (1,0) node[pos=0.5, below, thick, font=\large] {$K_1$}; 
    \draw[line width=2pt] (0,-1/10) -- (0,1/10); 
    \draw[line width=2pt] (1,-1/10) -- (1,1/10); 

    \draw[line width=2pt] (1,0) -- (2,0) node[pos=0.5, below, thick, font=\large] {$K_2$}; 
    \draw[line width=2pt] (2,-1/10) -- (2,1/10); 

    \draw[line width=2pt] (2,0) -- (3,0) node[pos=0.5, below, thick, font=\large] {$K_3$}; 
    \draw[line width=2pt] (3,-1/10) -- (3,1/10); 

    \addplot[BrickRed, ultra thick, domain=0:1] {1.25} node[pos=0.5,above] {$\phi_{1,1}$};
    \addplot[NavyBlue, ultra thick, domain=0:1] {x/2 + 0.25 +0.5} node[pos=0.2, below] {$\phi_{1,2}$};


    \addplot[BrickRed, ultra thick, domain=1:2] {1} node[pos=0.2,above=-2.5pt] {$\phi_{2,1}$};
    \addplot[NavyBlue, ultra thick, domain=1:2] {x/2 - 0.25 +0.5} node[pos=0.2,below] {$\phi_{2,2}$};

    \addplot[BrickRed, ultra thick, domain=2:3] {0.75} node[pos=0.8,above] {$\phi_{3,1}$};
    \addplot[NavyBlue, ultra thick, domain=2:3] {x/2 - 0.75 +0.5} node[pos=0.2,above=1pt] {$\phi_{3,2}$};



    \draw[loosely dashed, black] (axis cs: 0, 0) -- (axis cs: 0, 2);
    \draw[loosely dashed, black] (axis cs: 1, 0) -- (axis cs: 1, 2);
    \draw[loosely dashed, black] (axis cs: 2, 0) -- (axis cs: 2, 2);
    \draw[loosely dashed, black] (axis cs: 3, 0) -- (axis cs: 3, 2);



    
    \end{axis}
  \end{tikzpicture}
		\end{subfigure}
	}
	\caption{The formation of the agglomerated basis $\phi_{\mathrm{agg},m}$ (left) with respect to the original bases $\phi_{i,m}$ (right) with $\text{supp}(\phi_{i,m}) = K_i$. \label{fig:DgExt}}
\end{figure}
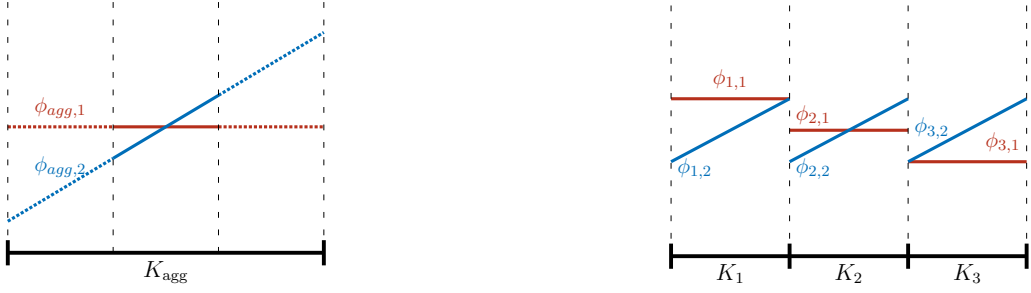
Agglomeration is performed by extending the solution basis of target cells to encompass their source cells so that the target cells merge with the source cells in their agglomeration group and represent the unified entity. Once the agglomerated basis is obtained, the problem is carried on the agglomerated space, replacing the original space. With respect to Eq.\eqref{eq:subSpace}, the agglomerated XDG space, denoted by $\mathbb{P}_p^{\mathrm{X},\mathrm{agg}}$, is hence defined as:
\begin{equation}
	\mathbb{P}_p^{\mathrm{X},\mathrm{agg}}(\mathfrak{K}_h^{\mathrm{X}}) := \mathbb{P}_p^{\mathrm{X}}(\mathfrak{K}_h^{\mathrm{X},\mathrm{agg}}).
\end{equation}
Since the agglomerated XDG space is a subspace of the original XDG space, its basis functions can be derived through a linear combination of the original space. This is achieved by utilizing a global coupling matrix, alternatively referred to as an injection operator (from the agglomerated space to the original space), denoted by $\bm{Q}$:
\begin{equation}
	\underline{\phi}^{\mathrm{X},\mathrm{agg}} = \underline{\phi}^{\mathrm{X}} \bm{Q},
\end{equation} 
where the corresponding global basis functions are defined as the row concatenation of cell-local bases $\underline{\phi}_i^{\mathrm{X}} = \left( \phi_{i,1}, ..., \phi_{i,M} \right)$, as $\underline{\phi}^{\mathrm{X},\mathrm{agg}}=[\underline{\phi}_1 ... \, \underline{\phi}_{ N - N_{\mathrm{src} }}]$ and $\underline{\phi}^{\mathrm{X}}=[\underline{\phi}_1 ... \, \underline{\phi}_{ N }]$, respectively. The coupling matrix $\bm{Q}$ is a real-valued matrix of the size $ M N \times M ( N - N_{\mathrm{src}} )$ and holding the coefficients of the basis transformation for the global basis functions. Consequently, the cell agglomeration can be performed by using the operator $\bm{Q}$.

When needed, it is also possible to re-establish the solution on the original space by projecting agglomerated space back. In general, agglomeration can be applied to both cut and uncut cells regardless of their designation as a source cell. The following section presents a basic illustration of agglomeration algebra using a simple example, without loss of generality (w.l.o.g.). 

\begin{example} \label{ex:aggAlgebra}
Consider an uncut 1-D mesh with $\mathfrak{K}=\{ K_1, K_2, K_3 \}$ as depicted in Figure~\ref{fig:DgExt}, where the agglomeration edges $E=\{( K_1, K_2), ( K_3, K_2) \}$ with $K_2$ as the chosen target cell of the agglomeration group.
 For the sake of simplicity, let us consider a first-order polynomial space (i.e. $p=1$) where the global basis vector $\underline{\phi}=[\underline{\phi}_{1} \; \underline{\phi}_{2} \; \underline{\phi}_{3}]$ is defined as the row concatenation of the local bases denoted by $\underline{\phi}_{i}=[\vec{\phi}_{i,1} \; \vec{\phi}_{i,2}]$. Hence, the agglomerated basis can be obtained by transforming the local bases as:
\begin{equation}
	\label{eq:aggBasis}
	\underline{\phi}^\mathrm{agg} = \underline{\phi}\bm{Q} = \underline{\phi}_1  \bm{Q}_{1,2} + \underline{\phi}_2 \bm{Q}_{2,2} + \underline{\phi}_3 \bm{Q}_{3,2},
\end{equation}
where $\bm{Q}_{i,j}$ is the local coupling matrix from $K_i$ to $K_j$ with dimensions $M \times M$. In essence, the coupling matrix $\bm{Q}_{i,j}$ functions as an extension/extrapolation operator, mapping from a source to a target. It is calculated based on the properties of the local bases and can be computed via $\bm{Q}_{i,j,m,n}=\int_{K_i} \phi_{i,m} \phi_{j,n} dV $ for orthonormal bases. 

\subsection{Forming chains} \label{sec:FormingChains}
\begin{figure}[ht!]
	\centering
	\begin{subfigure}[b]{0.45\textwidth}
		\centering
		\begin{tikzpicture}
    \pgfmathsetmacro{\xSize}{0.75}
    \pgfmathsetmacro{\ySize}{\xSize}
    \pgfmathsetmacro{\magicRatio}{1.01}

    \coordinate (p1) at (\xSize/2*0,\ySize*4*\magicRatio);
    \coordinate (p12) at (\xSize,\ySize*3);

    \coordinate (p2) at (\xSize*3/2,\ySize*5/2*\magicRatio);
    \coordinate (p3) at (\xSize*5/2,\ySize*5/2*\magicRatio);
    \coordinate (p4) at (\xSize*7/2,\ySize*5/2*\magicRatio);
    \coordinate (p5) at (\xSize*9/2,\ySize*5/2*\magicRatio);
    \coordinate (p6) at (\xSize*11/2,\ySize*5/2*\magicRatio);
    \coordinate (p7) at (\xSize*13/2,\ySize*5/2*\magicRatio);
    \coordinate (p78) at (\xSize*14/2,\ySize*6/2*\magicRatio);
    \coordinate (p8) at (\xSize*16/2,\ySize*4*\magicRatio);
    \draw[cyan,thick] plot [smooth] coordinates {(p1) (p12) (p2) (p3) (p4) (p5) (p6) (p7)  (p78) (p8)};
    \fill[cyan!10] (p1) -- plot [smooth] coordinates {(p1) (p12) (p2) (p3) (p4) (p5) (p6) (p7)  (p78) (p8)} -- (p8) -- (p8) -- cycle;

    \coordinate (p1) at (\xSize/2,\ySize*3/2*0);
    \coordinate (p12) at (\xSize,\ySize*1);
    \coordinate (p2) at (\xSize*3/2,\ySize*3/2);
    \coordinate (p3) at (\xSize*5/2,\ySize*3/2);
    \coordinate (p4) at (\xSize*7/2,\ySize*3/2);
    \coordinate (p5) at (\xSize*9/2,\ySize*3/2);
    \coordinate (p6) at (\xSize*11/2,\ySize*3/2);
    \coordinate (p7) at (\xSize*13/2,\ySize*3/2);
    \coordinate (p78) at (\xSize*14/2,\ySize*2/2);
    \coordinate (p8) at (\xSize*15/2,\ySize*3/2*0);
    
    \draw[cyan,thick] plot [smooth] coordinates {(p1) (p12) (p2) (p3) (p4) (p5) (p6) (p7) (p78) (p8)};
    \fill[cyan!10] (p1) -- plot [smooth] coordinates {(p1) (p12) (p2) (p3) (p4) (p5) (p6) (p7) (p78) (p8)} -- (p8) -- (7,0) -- cycle;
     
    
      \pgfmathsetmacro{\yRatio}{1.75}
    \draw [->]  (\xSize*13/6,\yRatio*\ySize) -- ++(-\xSize/2,0);
    \fill (\xSize*13/6,\yRatio*\ySize) circle (1pt);
    \draw [->]  (\xSize*19/6,\yRatio*\ySize) -- ++(-\xSize/2,0);
    \fill (\xSize*19/6,\yRatio*\ySize) circle (1pt);
    \draw [->]  (\xSize*25/6,\yRatio*\ySize) -- ++(-\xSize/2,0);
    \fill (\xSize*25/6,\yRatio*\ySize) circle (1pt);
    \draw [->]  (\xSize*31/6,\yRatio*\ySize) -- ++(-\xSize/2,0);
    \fill (\xSize*31/6,\yRatio*\ySize) circle (1pt);
  
      \pgfmathsetmacro{\yRatio}{2.25}
      \draw [->]  (\xSize*13/6,\yRatio*\ySize) -- ++(-\xSize/2,0);
      \fill (\xSize*13/6,\yRatio*\ySize) circle (1pt);
      \draw [->]  (\xSize*19/6,\yRatio*\ySize) -- ++(-\xSize/2,0);
      \fill (\xSize*19/6,\yRatio*\ySize) circle (1pt);
      \draw [->]  (\xSize*25/6,\yRatio*\ySize) -- ++(-\xSize/2,0);
      \fill (\xSize*25/6,\yRatio*\ySize) circle (1pt);
      \draw [->]  (\xSize*31/6,\yRatio*\ySize) -- ++(-\xSize/2,0);
      \fill (\xSize*31/6,\yRatio*\ySize) circle (1pt);
      
    \draw[thick,step=\xSize] (0,0) grid (8*\xSize,4*\xSize);
    
    \node at (0.4*\xSize,1.3*\ySize) {$\mathfrak{A}$};
    \node[text=cyan!90] at (7/2*\xSize,7/2*\ySize) {$\mathfrak{B}$};
  
  \end{tikzpicture}

	\end{subfigure}
	\hfill
	\begin{subfigure}[b]{0.45\textwidth}
		\centering
		\begin{tikzpicture}
  \pgfmathsetmacro{\xSize}{0.75}
  \pgfmathsetmacro{\ySize}{\xSize}
  
  \coordinate (p1) at (\xSize/2*0,\ySize*4);
  \coordinate (p12) at (\xSize,\ySize*3);

  \coordinate (p2) at (\xSize*3/2,\ySize*5/2);
  \coordinate (p3) at (\xSize*5/2,\ySize*5/2);
  \coordinate (p4) at (\xSize*7/2,\ySize*5/2);
  \coordinate (p5) at (\xSize*9/2,\ySize*5/2);
  \coordinate (p6) at (\xSize*11/2,\ySize*5/2);
  \coordinate (p7) at (\xSize*13/2,\ySize*5/2);
  \coordinate (p78) at (\xSize*14/2,\ySize*6/2);
  \coordinate (p8) at (\xSize*16/2,\ySize*4);
  \draw[cyan,thick] plot [smooth] coordinates {(p1) (p12) (p2) (p3) (p4) (p5) (p6) (p7)  (p78) (p8)};
  \fill[cyan!10] (p1) -- plot [smooth] coordinates {(p1) (p12) (p2) (p3) (p4) (p5) (p6) (p7)  (p78) (p8)} -- (p8) -- (p8) -- cycle;

  \coordinate (p1) at (\xSize/2,\ySize*3/2*0);
  \coordinate (p12) at (\xSize,\ySize*1);
  \coordinate (p2) at (\xSize*3/2,\ySize*3/2);
  \coordinate (p3) at (\xSize*5/2,\ySize*3/2);
  \coordinate (p4) at (\xSize*7/2,\ySize*3/2);
  \coordinate (p5) at (\xSize*9/2,\ySize*3/2);
  \coordinate (p6) at (\xSize*11/2,\ySize*3/2);
  \coordinate (p7) at (\xSize*13/2,\ySize*3/2);
  \coordinate (p78) at (\xSize*14/2,\ySize*2/2);
  \coordinate (p8) at (\xSize*15/2,\ySize*3/2*0);
  
  \draw[cyan,thick] plot [smooth] coordinates {(p1) (p12) (p2) (p3) (p4) (p5) (p6) (p7) (p78) (p8)};
  \fill[cyan!10] (p1) -- plot [smooth] coordinates {(p1) (p12) (p2) (p3) (p4) (p5) (p6) (p7) (p78) (p8)} -- (p8) -- (7,0) -- cycle;
   
  
    \pgfmathsetmacro{\yRatio}{1.75}
    \draw [->]  (\xSize*13/6,\yRatio*\ySize) -- ++(-\xSize/2,0);
    \fill (\xSize*13/6,\yRatio*\ySize) circle (1pt);
    \draw [->]  (\xSize*19/6,\yRatio*\ySize) -- ++(-\xSize/2,0);
    \fill (\xSize*19/6,\yRatio*\ySize) circle (1pt);
    \draw [->]  (\xSize*29/6,\yRatio*\ySize) -- ++(\xSize/2,0);
    \fill (\xSize*29/6,\yRatio*\ySize) circle (1pt);
    \draw [->]  (\xSize*35/6,\yRatio*\ySize) -- ++(\xSize/2,0);
    \fill (\xSize*35/6,\yRatio*\ySize) circle (1pt);

    \pgfmathsetmacro{\yRatio}{2.25}
    \draw [->]  (\xSize*13/6,\yRatio*\ySize) -- ++(-\xSize/2,0);
    \fill (\xSize*13/6,\yRatio*\ySize) circle (1pt);
    \draw [->]  (\xSize*19/6,\yRatio*\ySize) -- ++(-\xSize/2,0);
    \fill (\xSize*19/6,\yRatio*\ySize) circle (1pt);
    \draw [->]  (\xSize*29/6,\yRatio*\ySize) -- ++(\xSize/2,0);
    \fill (\xSize*29/6,\yRatio*\ySize) circle (1pt);
    \draw [->]  (\xSize*35/6,\yRatio*\ySize) -- ++(\xSize/2,0);
    \fill (\xSize*35/6,\yRatio*\ySize) circle (1pt);
    
  \draw[thick,step=\xSize] (0,0) grid (8*\xSize,4*\xSize);
  
  \node at (0.4*\xSize,1.3*\ySize) {$\mathfrak{A}$};
  \node[text=cyan!90] at (7/2*\xSize,7/2*\ySize) {$\mathfrak{B}$};

\end{tikzpicture}

	\end{subfigure}
	\caption{Unweighted (left) and weighted (right) chain agglomerations in a cut cell mesh. The arrows indicate the preferred way of agglomeration for the small-cut cells smaller than the half of background cell (shown with dots) in $\mathfrak{A}$. The subdomains $\mathfrak{A}$ and $\mathfrak{B}$ are shown in white and cyan, respectively.  \label{fig:ChainCutCell}}

 \end{figure}
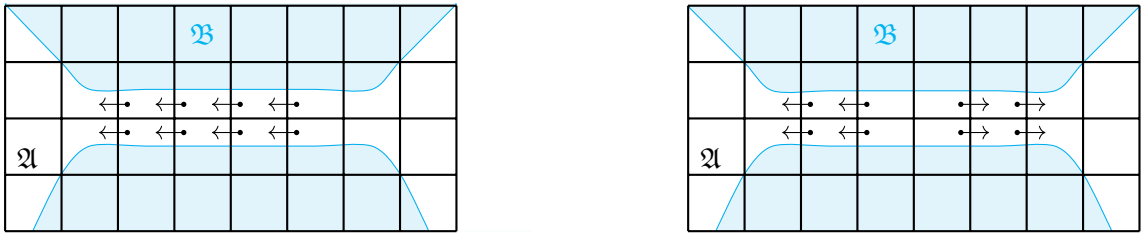

In complicated geometries, multiple adjacent source cells in the intersection regions may require to be agglomerated into a single target cell. These source cells may share a direct logical edge with the target cell, a process referred to as direct agglomeration or they can link together to form chains by connecting one another, a process referred to as chain agglomeration. Essentially, chains represent a particular subtype of agglomeration groups in which neighboring source cells are agglomerated to a single suitable target cell without sharing an edge. It should be noted that agglomeration groups can also be formed by multiple cells that share an edge with the target cell like in Example~\ref{ex:aggAlgebra}. Figure~\ref{fig:ChainCutCell} provides an illustrative example of the chain-forming process, which is frequently observed in phenomena like coalescence or the breakup of subregions.
In the case of direct agglomeration, selecting the target cell is relatively straightforward and it is often based on the cell size when multiple adjacent target cells are available to a single source cell. However, this simplicity does not apply to the chain agglomeration, where the absence of a shared logical edge and the inter-cell distance become important considerations.

As depicted in Figure~\ref{fig:ChainCutCell}, various selection criteria may yield several potential alternatives for forming different agglomeration chains. Furthermore, the selection process for chains unfolds over multiple steps due to the unpredictability of linking between source and target cells. In this context, we observed that the best agglomeration graphs are achieved by the selection criteria weighted on the distance and size. Therefore, our agglomeration strategy aims to target the closest cell to the source among the candidates available. In cases where multiple candidates are equidistant to the source, the target cell is selected based on its size, with preference given to the one with the largest fraction. If multiple candidates still exist, the selection is then determined by the lowest cell number. Because the weights of candidate cells are unique, at least the cell number, the resulting agglomeration mapping remains unique.

\subsection{Agglomeration levels} \label{sec:aggLevels}
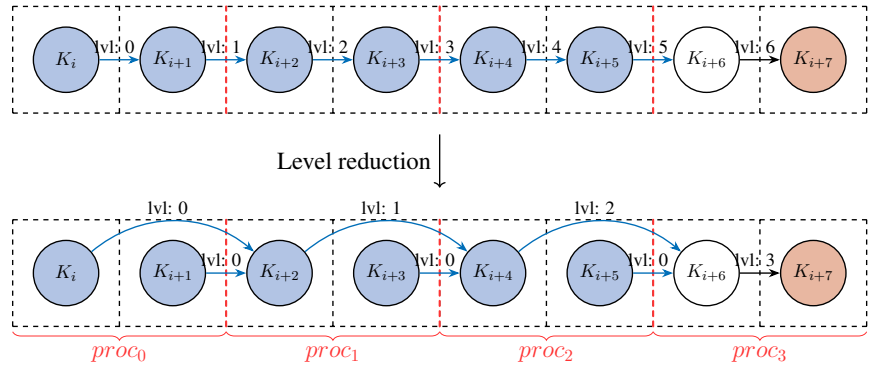
\begin{figure}[ht!]
	\centering
	\resizebox{0.7\textwidth}{!}{
	\begin{tikzpicture}
    \pgfmathsetmacro{\xSize}{2.2}
    \pgfmathsetmacro{\yMid}{2.5}

\draw[step=\xSize,thick,dashed] (0,\xSize*2) grid (\xSize*8,\xSize*3);    

\draw[red!80 , very thick, dashed, line cap=round,line join=round] (2*\xSize,\xSize*2) -- (2*\xSize,\xSize*3);
  \draw[red!80 , very thick, dashed, line cap=round,line join=round] (4*\xSize,\xSize*2) -- (4*\xSize,\xSize*3);
  \draw[red!80 , very thick, dashed, line cap=round,line join=round] (6*\xSize,\xSize*2) -- (6*\xSize,\xSize*3);
  
\begin{scope}[every node/.style={circle, thick, draw, minimum size=3.5em}]
    \node[fill=NavyBlue!30] (A) at  (\xSize*1/2,\xSize*\yMid) {$K_{i}$};
    \node[fill=NavyBlue!30] (B) at  (\xSize*3/2,\xSize*\yMid) {$K_{i+1}$};
    \node[fill=NavyBlue!30] (C) at  (\xSize*5/2,\xSize*\yMid) {$K_{i+2}$};
    \node[fill=NavyBlue!30] (D) at  (\xSize*7/2,\xSize*\yMid) {$K_{i+3}$};
    \node[fill=NavyBlue!30] (E) at  (\xSize*9/2,\xSize*\yMid) {$K_{i+4}$};
    \node[fill=NavyBlue!30] (F) at  (\xSize*11/2,\xSize*\yMid){$K_{i+5}$};
    \node[fill=gray!0] (G) at   (\xSize*13/2,\xSize*\yMid) {$K_{i+6}$};
    \node[fill=BrickRed!30] (H) at   (\xSize*15/2,\xSize*\yMid) {$K_{i+7}$};
\end{scope}

    
\begin{scope}[>={Stealth[NavyBlue]},
              every edge/.style={draw=NavyBlue, thick}]
    \path [->] (A) edge node[midway, above] { lvl: 0 \;} (B);
    \path [->] (B) edge node[midway, above] { lvl: 1 \;} (C); 

    \path [->] (C) edge node[midway, above] { lvl: 2 \;} (D);
    \path [->] (D) edge node[midway, above] { lvl: 3 \;} (E);
    \path [->] (E) edge node[midway, above] { lvl: 4 \;} (F);
    \path [->] (F) edge node[midway, above] { lvl: 5 \;} (G);

\end{scope}

\begin{scope}[>={Stealth[black]},
              every edge/.style={draw=black, thick}]

    \path [->] (G) edge node[midway, above] { lvl: 6 \;} (H);

\end{scope}

  \draw [->, thick]  (4*\xSize,1.8*\xSize) -- ++(0,-0.5*\xSize) node[midway,left=1pt] {\Large Level reduction};
    \pgfmathsetmacro{\xSize}{2.2}
    \pgfmathsetmacro{\yMid}{0.5}

\draw[step=\xSize,thick,dashed] (0,\xSize*0) grid (\xSize*8,\xSize*1);    
  \draw[red!80 , very thick, dashed, line cap=round,line join=round] (2*\xSize,\xSize*0) -- (2*\xSize,\xSize*1);
  \draw[red!80 , very thick, dashed, line cap=round,line join=round] (4*\xSize,\xSize*0) -- (4*\xSize,\xSize*1);
  \draw[red!80 , very thick, dashed, line cap=round,line join=round] (6*\xSize,\xSize*0) -- (6*\xSize,\xSize*1);
  
\begin{scope}[every node/.style={circle, thick, draw, minimum size=3.5em}]
    \node[fill=NavyBlue!30] (A2) at  (\xSize*1/2,\xSize*\yMid) {$K_{i}$};
    \node[fill=NavyBlue!30] (B2) at  (\xSize*3/2,\xSize*\yMid) {$K_{i+1}$};
    \node[fill=NavyBlue!30] (C2) at  (\xSize*5/2,\xSize*\yMid) {$K_{i+2}$};
    \node[fill=NavyBlue!30] (D2) at  (\xSize*7/2,\xSize*\yMid) {$K_{i+3}$};
    \node[fill=NavyBlue!30] (E2) at  (\xSize*9/2,\xSize*\yMid) {$K_{i+4}$};
    \node[fill=NavyBlue!30] (F2) at  (\xSize*11/2,\xSize*\yMid){$K_{i+5}$};
    \node[fill=gray!0] (G2) at   (\xSize*13/2,\xSize*\yMid) {$K_{i+6}$};
    \node[fill=BrickRed!30] (H2) at   (\xSize*15/2,\xSize*\yMid) {$K_{i+7}$};
\end{scope}

\begin{scope}[>={Stealth[NavyBlue]},
              every edge/.style={draw=NavyBlue, thick}]
    \path [->] (A2) edge[bend left=40] node[midway, above] { lvl: 0 \;} (C2);
    \path [->] (B2) edge node[midway, above] { lvl: 0 \;} (C2); 

    \path [->] (C2)  edge[bend left=40] node[midway, above] { lvl: 1 \;} (E2);
    \path [->] (D2) edge node[midway, above] { lvl: 0 \;} (E2);
    \path [->] (E2)  edge[bend left=40] node[midway, above] { lvl: 2 \;} (G2);
    \path [->] (F2) edge node[midway, above] { lvl: 0 \;} (G2);

\end{scope}

\begin{scope}[>={Stealth[black]},
              every edge/.style={draw=black, thick}]
    \path [->] (G2) edge node[midway, above] { lvl: 3 \;} (H2);
\end{scope}

  \draw [red!80, decorate,decoration={brace,amplitude=5pt,mirror}] (0*\xSize,-0.05*\xSize) -- (2*\xSize,-0.05*\xSize) node[text=red!80, midway,below=5pt] {\Large $proc_0$};

  \draw [red!80, decorate,decoration={brace,amplitude=5pt,mirror}] (2*\xSize,-0.05*\xSize) -- (4*\xSize,-0.05*\xSize) node[text=red!80, midway,below=5pt] {\Large $proc_1$};

  \draw [red!80, decorate,decoration={brace,amplitude=5pt,mirror}] (4*\xSize,-0.05*\xSize) -- (6*\xSize,-0.05*\xSize) node[text=red!80, midway,below=5pt] {\Large $proc_2$};

  \draw [red!80, decorate,decoration={brace,amplitude=5pt,mirror}] (6*\xSize,-0.05*\xSize) -- (8*\xSize,-0.05*\xSize) node[text=red!80, midway,below=5pt] {\Large $proc_3$};

  \end{tikzpicture}
	}
	\caption{Illustration of an inter-processor agglomeration chain on a 2D domain decomposed into four regions. The top chain displays the sequential pairs, while the bottom chain features the level-reduced equivalent of the top chain. The red lines indicate the processor boundaries, whereas the arrows in the cells indicate the agglomeration pairs from source to target. Final targets are indicated by red vertices, while the source cells in direct and chain agglomerations are displayed with white and blue vertices, respectively.
	\label{fig:InterProcCutCell}}
\end{figure}
 Agglomeration sources can be categorized into multiple levels based on their order in the respective mathematical operations, particularly when forming agglomeration chains.
 In cases of agglomeration groups consisting solely of direct agglomeration pairs, where sources and targets are adjacent, agglomeration occurs without any dependence. Hence, the source cells can be merged with the target cells in a single step without the need for a specific order using an extension operator as outlined in Eq.~\ref{eq:aggBasis}. 
 However, in cases of chain agglomeration, the source cells need to be agglomerated through the cells in their path to final targets while maintaining a sequence of operations, as a source in one pair may also function as a target cell in another pair. 
 
 For instance, let us consider Example~\ref{ex:aggAlgebra} with a revised agglomeration using $K_1 \rightarrow K_2$ and $K_2 \rightarrow K_3$ with $K_3$ as the target cell. In this case, the formal expression for the basis extrapolation would read as:
 \begin{align}
	\underline{\phi}^\mathrm{agg}_2 &= \underline{\phi}_1 \bm{Q}_{1,2} +  \underline{\phi}_2 \bm{Q}_{2,2},  \\
	\underline{\phi}^\mathrm{agg}_3 &=  \underline{\phi}^\mathrm{agg}_2 \bm{Q}_{2,3} + \underline{\phi}_3 \bm{Q}_{3,3},
 \end{align}
 where $K_1$ is first agglomerated to $K_2$ with the help of coupling matrix $\bm{Q}_{1,2}$ and subsequently the union of $K_1$ and $K_2$ is agglomerated to $K_3$ with $\bm{Q}_{2,3}$. In contrast to Eq.~\ref{eq:aggBasis}, the above formulation is computationally more expensive because it introduces additional operations at each sub-agglomeration, and the cost increases with longer chains. Moreover, for a correct calculation, it is necessary to perform the agglomeration from $K_1$ to $K_2$ and its respective operations before moving on to the agglomeration from $K_2$ to $K_3$. This sequencing is crucial to prevent the loss of information regarding the coordinate vector of $K_1$. Otherwise, $K_2$ would not contain the necessary information about $K_1$ 
 during its agglomeration process.
 
Therefore, to ensure a clear and logical sequential order, we establish agglomeration levels commencing at zero and increasing with each subsequent agglomeration pair. This level indicates one step higher than the highest level of the preceding agglomerations. Consequently, operations are systematically started from the lowest level and proceed to the highest consecutively. 
At first, the zero-th-level operations are conducted directly between source and target cells, without any dependencies on previous agglomerations. The remaining agglomeration pairs are only performed after the pairs with lower levels are agglomerated. It is also important to note that, a source cell can have several preceding agglomerations with different levels. 

 Furthermore, by converting chain agglomerations into equivalent pairs while ensuring connectivity, it is possible to decrease agglomeration levels and thus reduce computational effort for agglomeration chains. For instance, the aforementioned agglomerated cell can be obtained by substituting the agglomeration mapping with $K_1 \rightarrow K_3$ and $K_2 \rightarrow K_3$. This also enables a direct computation of coupling matrix $\bm{Q}_{1,3}$ between bases of $K_1$ and $K_3$, rather than establishing the relationship through the local coupling matrices of the intermediate element by $\bm{Q}_{1,2} \bm{Q}_{2,3}$.
 Since the substituted pairs of agglomerations do not involve lower levels, they are both at the zero-th level and the corresponding operations in computations can be executed without dependence on preceding agglomerations as:
  \begin{equation}
	\label{eq:aggBasisChainLvl}
	\underline{\phi}^\mathrm{agg} = \underline{\phi}_1 \bm{Q}_{1,3}  + \underline{\phi}_2  \bm{Q}_{2,3}+  \underline{\phi}_3  \bm{Q}_{3,3}.
\end{equation}
Hence, it is preferred to reduce levels of chain agglomerations if possible, but, it is not always feasible to decrease agglomeration levels.
 
In parallel executions, the computational domain is decomposed into multiple regions where the respective operations are performed separately on different processors before being coupled. Typically, the operations are coupled by utilizing a thin layer of cells at the processor boundaries (i.e., ghost cells), with no explicit information exchange occurring between the other cells. Therefore, it is not possible to replace inter-processor agglomeration pairs with their lower-level equivalents without compromising the parallelization. Figure~\ref{fig:InterProcCutCell} illustrates an example graph for the inter-processor agglomeration chains with levels. As the downstream agglomeration pairs require information about the upstream pairs that are in other processors, it is required to establish agglomeration levels to maintain the sequential order of operations. In such situations, the pairs are agglomerated to the processor boundary cells of the upper-level pairs, with further levels of agglomeration carried out by the respective processor. Thus, the maximum level of agglomeration is controlled by the number of processors utilized. Assuming a chain passes through a processor only once, the maximum level of agglomeration achievable is equal to the number of processors decreased by one. For these inter-processor agglomeration chains, we employ MPI-capable sparse matrix operations. 
\end{example}
\subsection{Small-cut agglomeration}
The cell agglomeration approach mainly aims to eliminate small-cut cells arising from the interface intersection by agglomerating them into adjacent cells, as shown in Figure~\ref{fig:CutCell}. This is done to ensure that the condition numbers are scaled with the size of regular cells, rather than being influenced by arbitrarily shaped small cells introduced during the numerical simulations. As such, agglomeration can be viewed as a preconditioner to the operator and mass matrices. Furthermore, small-cut agglomeration promotes the stability of discretizations by removing the troublesome elements. 

In essence, the small-cut agglomeration maps cut cells with a fraction below the user-defined threshold $0\leq \alpha \leq 1$ to the appropriate targets. The fraction of a cut cell is described by the ratio of its $D$-dimensional volume to the background cell:
 \begin{equation}
	 \sfrac(K_{i,\mathfrak{s}}) := |  K_{i,\mathfrak{s}} | / | K_{i} |.
 \end{equation}
Accordingly, the cells with $\sfrac(K_{i,\mathfrak{s}}) < \alpha$ are designated to be agglomeration sources and the set of these agglomeration sources is denoted by $\mathfrak{K}_{\mathrm{sml},\mathfrak{s}}$. 
 Hereby, the choice of $\alpha$ usually falls within the range of 0.1 and 0.3, as noted by the prior studies~\cite{kummer_extended_2017,Muller_2017} and our numerical tests confirm this finding (see §\ref{sec:Results}).

 \subsection{Temporal changes} 
Similar to the small-cut problem, the notion of agglomeration can be employed to accommodate temporal changes due to interface movements in the dynamic cut-cell mesh structure~\cite{Kummer_2018_time}. Within the time-dependent cut cell mesh structure, the subregions $\mathfrak{A}$ and $\mathfrak{B}$ are controlled by the implicitly defined level-set function $\psi(t)$, which can modify their domains over time. These modifications in the subregions can lead to alterations in their grid structure and respective discretizations, potentially causing issues such as mismatches in the dimensions or degrees of freedom across different time steps. 

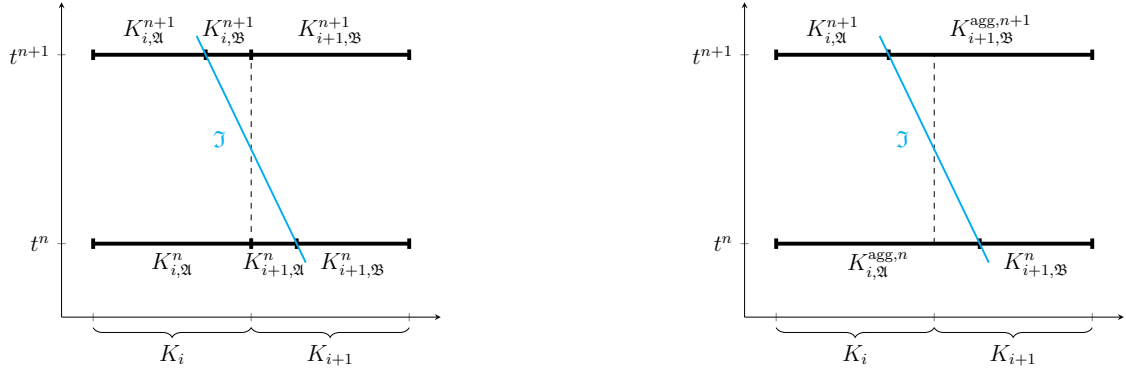
\begin{figure}[t!]
	\centering
	\resizebox{0.45\textwidth}{!}{ 
		\begin{subfigure}[b]{0.6\textwidth}
			\centering
			\begin{tikzpicture}
  \pgfmathsetmacro{\tDown}{0.22}
  \pgfmathsetmacro{\tUp}{0.4}
  \pgfmathsetmacro{\tH}{\tUp - \tDown}
  \pgfmathsetmacro{\slope}{0.8}
\begin{axis}[
  xlabel={},
  ylabel={},
  legend style={at={(0.5,-0.2)},anchor=north},
  ytick={\tDown,\tUp}, 
  yticklabels={$t^n$,$t^{n+1}$}, 
  xtick={0,0.5,1},
  xticklabels={}, 
  xmin=-0.1, xmax=1.1,
  ymin=0.15, ymax=0.45,
  clip=false,
  font=\large,
  axis lines=left 
]

\addplot[black, thick] coordinates {(0,\tDown) (1,\tDown)};
\addplot[black, thick] coordinates {(0,\tUp) (1,\tUp)};

\addplot[black, dashed] coordinates {(0.5,\tDown) (0.5,\tUp)};

\draw[line width=2pt] (axis cs: 0,\tDown-0.005) -- (axis cs: 0,\tDown+0.005); 
\draw[line width=2pt] (axis cs: 0,\tDown) -- (axis cs: 0.5,\tDown) node[pos=0.5, below] {$K_{i,\mathfrak{A}}^n$};
\draw[line width=2pt] (axis cs: 0.5,\tDown-0.005) -- (axis cs: 0.5,\tDown+0.005); 


\draw[line width=2pt] (axis cs: 0.5,\tDown) -- (axis cs: 0.5+\slope*\tH,\tDown) node[pos=0.5, below] {$K_{i+1,\mathfrak{A}}^n$}; 
\draw[line width=2pt] (axis cs:  0.5+\slope*\tH,\tDown-0.005) -- (axis cs:  0.5+\slope*\tH,\tDown+0.005); 

\draw[line width=2pt] (axis cs: 0.5+\slope*\tH,\tDown) -- (axis cs: 1,\tDown) node[pos=0.5, below] {$K_{i+1,\mathfrak{B}}^n$}; 
\draw[line width=2pt] (axis cs:  1,\tDown-0.005) -- (axis cs: 1,\tDown+0.005); 

\draw[line width=2pt] (axis cs: 0,\tUp-0.005) -- (axis cs: 0,\tUp+0.005); 
\draw[line width=2pt] (axis cs: 0,\tUp) -- (axis cs: 0.5-\slope*\tH,\tUp)  node[pos=0.5, above] {$K_{i,\mathfrak{A}}^{n+1}$};
\draw[line width=2pt] (axis cs:  0.5-\slope*\tH,\tUp-0.005) -- (axis cs:  0.5-\slope*\tH,\tUp+0.005); 

\draw[line width=2pt] (axis cs: 0.5-\slope*\tH,\tUp) -- (axis cs: 0.5,\tUp)  node[pos=0.5, above] {$K_{i,\mathfrak{B}}^{n+1}$};
\draw[line width=2pt] (axis cs:  0.5,\tUp-0.005) -- (axis cs:  0.5,\tUp+0.005); 

\draw[line width=2pt] (axis cs: 0.5,\tUp) -- (axis cs: 1,\tUp) node[pos=0.5, above] {$K_{i+1,\mathfrak{B}}^{n+1}$};
\draw[line width=2pt] (axis cs: 1,\tUp-0.005) -- (axis cs: 1,\tUp+0.005); 

\draw [decorate,decoration={brace,amplitude=5pt,mirror}] (axis cs: 0,0.14) -- (axis cs: 0.5,0.14) node[midway,below=5pt] {$K_{i}$};
\draw [decorate,decoration={brace,amplitude=5pt,mirror}] (axis cs: 0.5,0.14) -- (axis cs: 1,0.14) node[midway,below=5pt] {$K_{i+1}$};

\addplot[cyan, thick, line width=1pt] coordinates {(0.5+\slope*\tH*1.2),\tDown-\tH*0.1) (0.5-\slope*\tH*1.2,\tUp+\tH*0.1)};
\node[text=cyan, line width=1pt] at (axis cs: 0.4, 0.32) {$\mathfrak{I}$};


\end{axis}
\end{tikzpicture}
		\end{subfigure}
	}
		\hfill
		\resizebox{0.45\textwidth}{!}{
		\begin{subfigure}[b]{0.6\textwidth}
			\centering
			\begin{tikzpicture}
  \pgfmathsetmacro{\tDown}{0.22}
  \pgfmathsetmacro{\tUp}{0.4}
  \pgfmathsetmacro{\tH}{\tUp - \tDown}
  \pgfmathsetmacro{\slope}{0.8}
\begin{axis}[
  xlabel={},
  ylabel={},
  legend style={at={(0.5,-0.2)},anchor=north},
  ytick={\tDown,\tUp}, 
  yticklabels={$t^n$,$t^{n+1}$}, 
  xtick={0,0.5,1},
  xticklabels={}, 
  xmin=-0.1, xmax=1.1,
  ymin=0.15, ymax=0.45,
  clip=false,
  font=\large,
  axis lines=left 
]

\addplot[black, thick] coordinates {(0,\tDown) (1,\tDown)};
\addplot[black, thick] coordinates {(0,\tUp) (1,\tUp)};

\addplot[black, dashed] coordinates {(0.5,\tDown) (0.5,\tUp)};

\draw[line width=2pt] (axis cs: 0,\tDown-0.005) -- (axis cs: 0,\tDown+0.005); 
\draw[line width=2pt] (axis cs: 0,\tDown) -- (axis cs: 0.5+\slope*\tH,\tDown) node[pos=0.5, below] {$K_{i,\mathfrak{A}}^{\mathrm{agg},n}$};
\draw[line width=2pt] (axis cs:  0.5+\slope*\tH,\tDown-0.005) -- (axis cs:  0.5+\slope*\tH,\tDown+0.005); 

\draw[line width=2pt] (axis cs: 0.5+\slope*\tH,\tDown) -- (axis cs: 1,\tDown) node[pos=0.5, below] {$K_{i+1,\mathfrak{B}}^n$}; 
\draw[line width=2pt] (axis cs:  1,\tDown-0.005) -- (axis cs: 1,\tDown+0.005); 

\draw[line width=2pt] (axis cs: 0,\tUp-0.005) -- (axis cs: 0,\tUp+0.005); 
\draw[line width=2pt] (axis cs: 0,\tUp) -- (axis cs: 0.5-\slope*\tH,\tUp)  node[pos=0.5, above] {$K_{i,\mathfrak{A}}^{n+1}$};
\draw[line width=2pt] (axis cs:  0.5-\slope*\tH,\tUp-0.005) -- (axis cs:  0.5-\slope*\tH,\tUp+0.005); 


\draw[line width=2pt] (axis cs: 0.5-\slope*\tH,\tUp) -- (axis cs: 1,\tUp) node[pos=0.5, above] {$K_{i+1,\mathfrak{B}}^{\mathrm{agg},n+1}$};
\draw[line width=2pt] (axis cs: 1,\tUp-0.005) -- (axis cs: 1,\tUp+0.005); 

\draw [decorate,decoration={brace,amplitude=5pt,mirror}] (axis cs: 0,0.14) -- (axis cs: 0.5,0.14) node[midway,below=5pt] {$K_{i}$};
\draw [decorate,decoration={brace,amplitude=5pt,mirror}] (axis cs: 0.5,0.14) -- (axis cs: 1,0.14) node[midway,below=5pt] {$K_{i+1}$};

\addplot[cyan, thick, line width=1pt] coordinates {(0.5+\slope*\tH*1.2),\tDown-\tH*0.1) (0.5-\slope*\tH*1.2,\tUp+\tH*0.1)};
\node[text=cyan, line width=1pt] at (axis cs: 0.4, 0.32) {$\mathfrak{I}$};


\end{axis}
\end{tikzpicture}
		\end{subfigure}
	}
	\caption{A visual representation of the temporal evolution of the moving interface $\mathfrak{I}$ on a 1D cut-cell mesh. As time progresses from $t = t^n$ to $t = t^{n+1}$, the cut cell $K_{i+1,\mathfrak{A}}^n$ vanishes, and concurrently, the new cut cell $K_{i,\mathfrak{B}}^{n+1}$ emerges. The left figure displays the original cells, while the right figure displays the agglomeration for topological consistency. \label{fig:VanishingNewbornCutCell}}
\end{figure}

In general, the motion of the interface $\mathfrak{I}$ can be categorized into two distinct types of topological change between the time steps $t=t^n$ and $t=t^{n+1}$, as shown in Figure~\ref{fig:VanishingNewbornCutCell}. These are:
\begin{itemize}
  \item vanishing cells if $\sfrac(K_{i,\mathfrak{s}}^{n}) \geq 0$ and $\sfrac(K_{i,\mathfrak{s}}^{n+1}) = 0$,
  \item newborn cells if $\sfrac(K_{i,\mathfrak{s}}^{n}) = 0$ and $\sfrac(K_{i,\mathfrak{s}}^{n+1}) \geq 0$. 
\end{itemize}
In simpler terms, phase cells must disappear in the following time step if the interface movement reduces the domain of interest, and conversely, new phase cells must emerge if the domain expands. As a result, topological structure changes with the appearance and disappearance of cells, ultimately causing a conceptual problem with the temporal calculations. 

The BoSSS software package offers two approaches to discretize time for the cases involving dynamic interfaces, namely \textit{splitting} and \textit{moving interface} approaches, implemented based on \cite{Kummer_2018_time}. In the splitting approach, firstly the interface is moved and then the values at time step $t^n$ are extrapolated to the new mesh at time step $t^{n+1}$. Subsequently, temporal integration is performed as if the interface was stationary between time steps. This approach can be demonstrated using the explicit Euler time scheme with regards to Equation~\ref{eq:simplifiedWeak} as:
\begin{equation}
	\text{Splitting approach:} \quad \mathcal{M}^{n+1}(\underline{\tilde{c}}^{n+1}-\underline{\tilde{c}}^{n}) / \Delta t + \mathcal{F}_\mathrm{sp} = 0,
\end{equation}
where $\mathcal{F}_\mathrm{sp}$ represents the terms related to the spatial discretization. In this approach, newborn cells introduce singular blocks into the mass matrix. To preserve consistent topology across the time steps, these newborn cells are hence agglomerated with suitable targets. 

On the other hand, the moving interface approach considers the interface as a dynamic entity and recovers its behavior within a time step by introducing a space-time notation. Accordingly, its temporal discretization with an explicit Euler time scheme leads to the following form with respect to Eq.~\ref{eq:simplifiedWeak}:
  \begin{align}
	\text{Moving interface approach:} \quad & (\mathcal{M}^{n+1}\underline{\tilde{c}}^{n+1}-\mathcal{M}^{n}\underline{\tilde{c}}^{n}) / \Delta t + \mathcal{F}_\mathrm{mov} = 0,
	\end{align}
where $\mathcal{F}_\mathrm{mov}$ represents the spatial discretization in the moving interface approach. In this method, both newborn and vanishing cells must be taken into account to maintain a topologically consistent discretization. As a result, these cells are marked as source cells and agglomerated to their appropriate neighbors.

Compared to each other, the splitting approach presents a more affordable solution, whereas the moving interface method incurs higher costs but provides superior accuracy. Ultimately, the choice of interface evolution involves a trade-off between accuracy and cost. For a comparison of these two methods, the reader is referred to the original work~\cite{Kummer_2018_time}. In this study, we restrict our numerical tests to the splitting approach in order to reduce the computational requirements.

In addition, it is important to refrain newborn and vanishing cells from being agglomeration targets, regardless of their size, as they cause topological disparities. Furthermore, interface movement should be limited to a finite speed to prevent excessive topological changes that can reduce accuracy. Therefore, a limit of one cell per time step for interface movement is implemented.




\subsection{Coincidence of cell boundaries with interfaces}
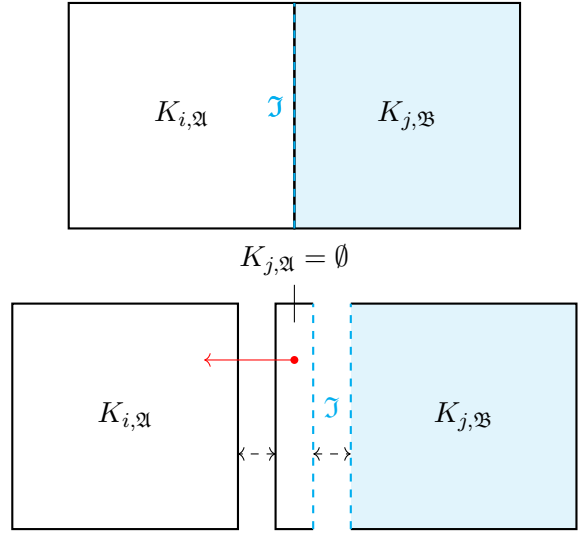
\begin{figure}[t!]
	\centering
	\begin{tikzpicture}

    \pgfmathsetmacro{\xSize}{3}
    \pgfmathsetmacro{\ySize}{3}
	\pgfmathsetmacro{\offSet}{1}

	
	  \fill[cyan!10] (2*\xSize, \ySize) -- (\xSize, \ySize) -- (\xSize,0) -- (2*\xSize,0) -- cycle;
	  \draw[thick] (0,0) rectangle (\xSize,\ySize);
	  \draw[thick] (\xSize,0) rectangle (2*\xSize,\ySize);
	  \draw[domain=0:3, samples=100, cyan, thick, dashed] plot ({\xSize}, {(\x)});

	  \node[text=cyan] at (\xSize*11/12,\ySize*6.5/12) {$\mathfrak{I}$};
	  
	  \node at (\xSize/2,\ySize/2) {$K_{i,\mathfrak{A}}$};
	  \node at (3/2*\xSize,\ySize/2) {$K_{j,\mathfrak{B}}$};
  
	\fill[cyan!10] (2*\xSize+0.75, -\offSet) -- (\xSize+0.75, -\offSet) -- (\xSize+0.75,-\ySize-\offSet) -- (2*\xSize+0.75,-\ySize-\offSet) -- cycle;
	\draw[thick] (0-0.75,-\ySize-\offSet) rectangle (\xSize-0.75,-\offSet);
	\draw[thick] (\xSize+0.25,-\offSet) -- (\xSize-0.25,-\offSet) -- (\xSize-0.25,-\ySize-\offSet) -- (\xSize+0.25,-\offSet-\ySize);
	\draw[thick] (\xSize+0.75,-\offSet) -- (2*\xSize+0.75,-\offSet) -- (2*\xSize+0.75,-\ySize-\offSet) -- (\xSize+0.75,-\ySize-\offSet) ;
	\draw[domain=0:3, samples=100, cyan, thick, dashed] plot ({\xSize+0.25}, {(\x-\ySize-\offSet)});
	\draw[domain=0:3, samples=100, cyan, thick, dashed] plot ({\xSize+0.75}, {(\x-\ySize-\offSet)});
	
	\node[text=cyan] at (\xSize*11/12,\ySize*6.5/12) {$\mathfrak{I}$};
	
	\node at (0.75,\ySize/2-\ySize-\offSet) {$K_{i,\mathfrak{A}}$};
	\draw (\xSize,-5/4*\offSet) -- (\xSize,-3*\offSet/4) node[left,above=0pt] {$K_{j,\mathfrak{A}} = \emptyset$};
	\node at (5.25,\ySize/2-\ySize-\offSet) {$K_{j,\mathfrak{B}}$};
 	\node[text=cyan] at (\xSize*13/12+0.25,\ySize*6.5/12-\ySize-\offSet) {$\mathfrak{I}$};
 	
 	\draw [dashed, <->] (\xSize*11/12,-\ySize) --  (\xSize*9/12,-\ySize);
 	\draw [dashed, <->] (\xSize*13/12,-\ySize) --  (\xSize*15/12,-\ySize) ;

	 \fill [red] (\xSize,-7/4*\offSet) circle (1.5pt);
	 \draw [red,->]  (\xSize,-7*\offSet/4) -- ++(-\xSize*0.4,0);



  



\end{tikzpicture}
	\caption{Illustration of the species coupling in case of a coinciding interface. In the figure, the dashed interface lies on top of the edge between two cells. It is assumed that $K_{j}$ in Equation~\ref{eq:coinFrac} is the cell with the lower index and owns the coinciding interface. The affected edge $\partial K_i \cap \partial K_j$ belongs to the species $\mathfrak{A}$ and takes care of the coupling between cells $K_{i,\mathfrak{A}}$ and $K_{j,\mathfrak{A}}$. The species are then coupled inside $K_{j}$ via the interface, from the empty cell $K_{j,\mathfrak{A}} = \emptyset$ to the full cell $K_{j,\mathfrak{B}} = K_{j}$. Finally, by performing the agglomeration, the discrete system is algebraically modified. This modification eliminates the (edge) contributions on $\partial K_{j,\mathfrak{A}}$ of the empty cell and combines the cell and species coupling, establishing the connection between $K_{i,\mathfrak{A}}$ and $K_{j,\mathfrak{B}}$. The lower part of the figure shows an exploded view of the situation to clarify the connectivity. \label{fig:CoinLevSet}}
\end{figure}
In rare cases, the interface may coincide with an edge of the background mesh, i.e. 
\begin{align}
	\left(\partial K_i \cap \partial K_j\right) \cap \mathfrak{I} = \left(\partial K_i \cap \partial K_j\right).
\end{align}
Then, one species completely fills one cell, while it vanishes in the other, w.l.o.g.:
\begin{align}\label{eq:coinFrac}
	\sfrac(K_{i,\mathfrak{A}}) &= 1 \land \sfrac(K_{j,\mathfrak{A}}) = 0 \\
	\sfrac(K_{i,\mathfrak{B}}) &= 0 \land \sfrac(K_{j,\mathfrak{B}}) = 1.
\end{align}
Such a situation is shown in Fig.~\ref{fig:CoinLevSet}.
Owing to the structure of the discretization couplings between species are occurring within cells, i.e. each part of the interface is assigned to exactly one cell. Moreover, the coupling between cells occurs via edges within one species, thereby each part of a cell boundary is uniquely assigned to a species. 
Since the underlying quadrature rules are determined on a cell-local level, small round-of-errors can then lead to a detection of the interface in both or none of the cells. Therefore, special care has to be taken in the treatment of coinciding interfaces. In practice, this is achieved by assigning coinciding interface parts to the cell with the lower index. This ensures that the interface is only treated once and in a consistent manner across processor boundaries and the edge is also assigned to exactly one species. Then, coupling terms on the interface might link to completely empty cells. However, since these empty cut cells are agglomerated, a well-defined situation is recovered.

%
\subsection{Implementation} \label{sec:agg_implementation}
\begin{figure}[t!]
    \captionsetup{font=large}
	\centering
	\resizebox{0.40\textwidth}{!}{ 
	\begin{subfigure}[b]{0.40\textwidth}
		\centering
		\begin{tikzpicture}
    \pgfmathsetmacro{\xSize}{0.75}
    \pgfmathsetmacro{\ySize}{\xSize}
    \pgfmathsetmacro{\lineR}{0.9}

    \pgfmathsetmacro{\R}{2.7}
    \pgfmathsetmacro{\sR}{0.6}

\fill [red] (-3.8*\xSize,-0.5*\ySize) circle (1.5pt);
\fill [red] (-3.8*\xSize,0.5*\ySize) circle (1.5pt);
\fill [red] (3.8*\xSize,-0.5*\ySize) circle (1.5pt);
\fill [red] (3.8*\xSize,0.5*\ySize) circle (1.5pt);

\fill [red] (3.5*\xSize,-1.5*\ySize) circle (1.5pt);

\fill [red] (3.5*\xSize,1.5*\ySize) circle (1.5pt);

\fill [red] (-3.5*\xSize,-1.5*\ySize) circle (1.5pt);

\fill [red] (-3.5*\xSize,1.5*\ySize) circle (1.5pt);

\fill [red] (-0.5*\ySize, -3.8*\xSize) circle (1.5pt);
\fill [red] (0.5*\ySize, -3.8*\xSize) circle (1.5pt);
\fill [red] (-0.5*\ySize, 3.8*\xSize) circle (1.5pt);
\fill [red] (0.5*\ySize, 3.8*\xSize) circle (1.5pt);

\fill [red] (-1.5*\ySize, -3.5*\xSize) circle (1.5pt);

\fill [red] (1.5*\ySize, -3.5*\xSize) circle (1.5pt);

\fill [red] (-1.5*\ySize, 3.5*\xSize) circle (1.5pt);

\fill [red] (1.5*\ySize, 3.5*\xSize) circle (1.5pt);

\pgfmathsetmacro{\lineR}{0.5}
\fill [red] (2.75*\xSize,-2.75*\ySize) circle (1.5pt);

\fill [red] (-2.75*\xSize,-2.75*\ySize) circle (1.5pt);

\fill [red] (2.75*\xSize,2.75*\ySize) circle (1.5pt);

\fill [red] (-2.75*\xSize,2.75*\ySize) circle (1.5pt);

\draw[cyan,thick] plot[domain=0:360, smooth, samples=100] ({\R*cos(\x)}, {\R*sin(\x)});

\fill[cyan!10] plot[domain=0:360, smooth, samples=100] ({\R*cos(\x)}, {\R*sin(\x)}) -- cycle;

\draw[cyan,thick] plot[domain=0:360, smooth, samples=100] ({\sR*cos(\x)}, {\sR*sin(\x)});

\fill[white] plot[domain=0:360, smooth, samples=100] ({\sR*cos(\x)}, {\sR*sin(\x)}) -- cycle;

\fill [red] (-0.35*\ySize, -0.35*\xSize) circle (1.5pt);
\fill [red] (-0.35*\ySize, 0.35*\xSize) circle (1.5pt);
\fill [red] (0.35*\ySize, -0.35*\xSize) circle (1.5pt);
\fill [red] (0.35*\ySize, 0.35*\xSize) circle (1.5pt);

\draw[step=\xSize, line width=0.03cm] (-3,-3) grid (4*\xSize,4*\xSize);

    \node at (-3.6*\xSize,-3.7*\ySize) {$\mathfrak{A}$};
    \node[text=cyan!90] at (1/2*\xSize,5/2*\ySize) {$\mathfrak{B}$};
   \clip (-3,-3) grid (4*\xSize,4*\xSize);     
  \end{tikzpicture}
		\caption{Source identification}
	\end{subfigure}
}
	\resizebox{0.40\textwidth}{!}{ 
		\begin{subfigure}[b]{0.40\textwidth}
			\centering
			\begin{tikzpicture}
    \pgfmathsetmacro{\xSize}{0.75}
    \pgfmathsetmacro{\ySize}{\xSize}
    \pgfmathsetmacro{\lineR}{0.9}

    \pgfmathsetmacro{\R}{2.7}
    \pgfmathsetmacro{\sR}{0.6}

\fill (-3.8*\xSize,-0.5*\ySize) circle (1.5pt);
\fill (-3.8*\xSize,0.5*\ySize) circle (1.5pt);
\fill (3.8*\xSize,-0.5*\ySize) circle (1.5pt);
\fill (3.8*\xSize,0.5*\ySize) circle (1.5pt);

\fill (3.5*\xSize,-1.5*\ySize) circle (1.5pt);
\draw [red,->]  (3.5*\xSize,-1.5*\ySize) -- ++(0,-\xSize*\lineR);

\fill (3.5*\xSize,1.5*\ySize) circle (1.5pt);
\draw [red,->]  (3.5*\xSize,1.5*\ySize) -- ++(0,\xSize*\lineR);

\fill (-3.5*\xSize,-1.5*\ySize) circle (1.5pt);
\draw [red,->]  (-3.5*\xSize,-1.5*\ySize) -- ++(0,-\xSize*\lineR);

\fill (-3.5*\xSize,1.5*\ySize) circle (1.5pt);
\draw [red,->]  (-3.5*\xSize,1.5*\ySize) -- ++(0,\xSize*\lineR);

\fill (-0.5*\ySize, -3.8*\xSize) circle (1.5pt);
\fill (0.5*\ySize, -3.8*\xSize) circle (1.5pt);
\fill (-0.5*\ySize, 3.8*\xSize) circle (1.5pt);
\fill (0.5*\ySize, 3.8*\xSize) circle (1.5pt);

\fill (-1.5*\ySize, -3.5*\xSize) circle (1.5pt);
\draw [red,->]  (-1.5*\ySize, -3.5*\xSize) -- ++(-\xSize*\lineR, 0);

\fill (1.5*\ySize, -3.5*\xSize) circle (1.5pt);
\draw [red,->]  (1.5*\ySize, -3.5*\xSize) -- ++(\xSize*\lineR, 0);

\fill (-1.5*\ySize, 3.5*\xSize) circle (1.5pt);
\draw [red,->]  (-1.5*\ySize, 3.5*\xSize) -- ++(-\xSize*\lineR, 0);

\fill (1.5*\ySize, 3.5*\xSize) circle (1.5pt);
\draw [red,->]  (1.5*\ySize, 3.5*\xSize) -- ++(\xSize*\lineR, 0);

\pgfmathsetmacro{\lineR}{0.5}
\fill (2.75*\xSize,-2.75*\ySize) circle (1.5pt);
\draw [red,->]  (2.75*\xSize,-2.75*\ySize) -- ++(0,-\xSize*\lineR);

\fill (-2.75*\xSize,-2.75*\ySize) circle (1.5pt);
\draw [red,->]  (-2.75*\xSize,-2.75*\ySize) -- ++(-\xSize*\lineR,0);

\fill (2.75*\xSize,2.75*\ySize) circle (1.5pt);
\draw [red,->]  (2.75*\xSize,2.75*\ySize) -- ++(0,\xSize*\lineR);

\fill (-2.75*\xSize,2.75*\ySize) circle (1.5pt);
\draw [red,->]  (-2.75*\xSize,2.75*\ySize) -- ++(-\xSize*\lineR,0);

\draw[cyan,thick] plot[domain=0:360, smooth, samples=100] ({\R*cos(\x)}, {\R*sin(\x)});

\fill[cyan!10] plot[domain=0:360, smooth, samples=100] ({\R*cos(\x)}, {\R*sin(\x)}) -- cycle;

\draw[cyan,thick] plot[domain=0:360, smooth, samples=100] ({\sR*cos(\x)}, {\sR*sin(\x)});

\fill[white] plot[domain=0:360, smooth, samples=100] ({\sR*cos(\x)}, {\sR*sin(\x)}) -- cycle;

\fill (-0.35*\ySize, -0.35*\xSize) circle (1.5pt);
\fill (-0.35*\ySize, 0.35*\xSize) circle (1.5pt);
\fill (0.35*\ySize, -0.35*\xSize) circle (1.5pt);
\fill (0.35*\ySize, 0.35*\xSize) circle (1.5pt);

\draw[step=\xSize, line width=0.03cm] (-3,-3) grid (4*\xSize,4*\xSize);

    \node at (-3.6*\xSize,-3.7*\ySize) {$\mathfrak{A}$};
    \node[text=cyan!90] at (1/2*\xSize,5/2*\ySize) {$\mathfrak{B}$};
   \clip (-3,-3) grid (4*\xSize,4*\xSize);     
  \end{tikzpicture}
			\caption{Direct agglomeration}
		\end{subfigure}
	}
		\vfill
		\vfill
	\resizebox{0.40\textwidth}{!}{ 
		\begin{subfigure}[b]{0.40\textwidth}
			\centering
			\begin{tikzpicture}
  \pgfmathsetmacro{\xSize}{0.75}
  \pgfmathsetmacro{\ySize}{\xSize}
  \pgfmathsetmacro{\lineR}{0.9}

  \pgfmathsetmacro{\R}{2.7}
  \pgfmathsetmacro{\sR}{0.6}

\fill (3.5*\xSize,-1.5*\ySize) circle (1.5pt);
\draw [->]  (3.5*\xSize,-1.5*\ySize) -- ++(0,-\xSize*\lineR);

\fill (3.5*\xSize,1.5*\ySize) circle (1.5pt);
\draw [->]  (3.5*\xSize,1.5*\ySize) -- ++(0,\xSize*\lineR);

\fill (-3.5*\xSize,-1.5*\ySize) circle (1.5pt);
\draw [->]  (-3.5*\xSize,-1.5*\ySize) -- ++(0,-\xSize*\lineR);

\fill (-3.5*\xSize,1.5*\ySize) circle (1.5pt);
\draw [->]  (-3.5*\xSize,1.5*\ySize) -- ++(0,\xSize*\lineR);


\fill (-1.5*\ySize, -3.5*\xSize) circle (1.5pt);
\draw [->]  (-1.5*\ySize, -3.5*\xSize) -- ++(-\xSize*\lineR, 0);

\fill (1.5*\ySize, -3.5*\xSize) circle (1.5pt);
\draw [->]  (1.5*\ySize, -3.5*\xSize) -- ++(\xSize*\lineR, 0);

\fill (-1.5*\ySize, 3.5*\xSize) circle (1.5pt);
\draw [->]  (-1.5*\ySize, 3.5*\xSize) -- ++(-\xSize*\lineR, 0);

\fill (1.5*\ySize, 3.5*\xSize) circle (1.5pt);
\draw [->]  (1.5*\ySize, 3.5*\xSize) -- ++(\xSize*\lineR, 0);

\pgfmathsetmacro{\lineR}{0.5}
\fill (2.75*\xSize,-2.75*\ySize) circle (1.5pt);
\draw [->]  (2.75*\xSize,-2.75*\ySize) -- ++(0,-\xSize*\lineR);

\fill (-2.75*\xSize,-2.75*\ySize) circle (1.5pt);
\draw [->]  (-2.75*\xSize,-2.75*\ySize) -- ++(-\xSize*\lineR,0);

\fill (2.75*\xSize,2.75*\ySize) circle (1.5pt);
\draw [->]  (2.75*\xSize,2.75*\ySize) -- ++(0,\xSize*\lineR);

\fill (-2.75*\xSize,2.75*\ySize) circle (1.5pt);
\draw [->]  (-2.75*\xSize,2.75*\ySize) -- ++(-\xSize*\lineR,0);

  \pgfmathsetmacro{\lineR}{1.9}
\fill (-3.8*\xSize,-0.5*\ySize) circle (1.5pt);
\draw [red,->]  (-3.8*\xSize,-0.5*\ySize) -- ++(0,-\xSize*\lineR);

\fill (-3.8*\xSize,0.5*\ySize) circle (1.5pt);
\draw [red,->]  (-3.8*\xSize,0.5*\ySize) -- ++(0,\xSize*\lineR);

\fill (3.8*\xSize,-0.5*\ySize) circle (1.5pt);
\draw [red,->]  (3.8*\xSize,-0.5*\ySize) -- ++(0,-\xSize*\lineR);

\fill (3.8*\xSize,0.5*\ySize) circle (1.5pt);
\draw [red,->]  (3.8*\xSize,0.5*\ySize) -- ++(0,\xSize*\lineR);

\fill (-0.5*\ySize, -3.8*\xSize) circle (1.5pt);
\draw [red,->]  (-0.5*\ySize, -3.8*\xSize) -- ++(-\xSize*\lineR, 0);

\fill (0.5*\ySize, -3.8*\xSize) circle (1.5pt);
\draw [red,->]  (0.5*\ySize, -3.8*\xSize) -- ++(\xSize*\lineR, 0);

\fill (-0.5*\ySize, 3.8*\xSize) circle (1.5pt);
\draw [red,->]  (-0.5*\ySize, 3.8*\xSize) -- ++(-\xSize*\lineR, 0);

\fill (0.5*\ySize, 3.8*\xSize) circle (1.5pt);
\draw [red,->]  (0.5*\ySize, 3.8*\xSize) -- ++(\xSize*\lineR, 0);

\draw[cyan,thick] plot[domain=0:360, smooth, samples=100] ({\R*cos(\x)}, {\R*sin(\x)});

\fill[cyan!10] plot[domain=0:360, smooth, samples=100] ({\R*cos(\x)}, {\R*sin(\x)}) -- cycle;

\draw[cyan,thick] plot[domain=0:360, smooth, samples=100] ({\sR*cos(\x)}, {\sR*sin(\x)});

\fill[white] plot[domain=0:360, smooth, samples=100] ({\sR*cos(\x)}, {\sR*sin(\x)}) -- cycle;

\fill (-0.35*\ySize, -0.35*\xSize) circle (1.5pt);
\fill (-0.35*\ySize, 0.35*\xSize) circle (1.5pt);
\fill (0.35*\ySize, -0.35*\xSize) circle (1.5pt);
\fill (0.35*\ySize, 0.35*\xSize) circle (1.5pt);

  \draw[step=\xSize, line width=0.03cm] (-3,-3) grid (4*\xSize,4*\xSize);

  \node at (-3.6*\xSize,-3.7*\ySize) {$\mathfrak{A}$};
  \node[text=cyan!90] at (1/2*\xSize,5/2*\ySize) {$\mathfrak{B}$};

\end{tikzpicture}
			\caption{Chain formation}
		\end{subfigure}
	}
	\resizebox{0.40\textwidth}{!}{ 
	\begin{subfigure}[b]{0.40\textwidth}
		\centering
		\begin{tikzpicture}
  \pgfmathsetmacro{\xSize}{0.75}
  \pgfmathsetmacro{\ySize}{\xSize}
  \pgfmathsetmacro{\lineR}{0.9}

  \pgfmathsetmacro{\R}{2.7}
  \pgfmathsetmacro{\sR}{0.6}




\fill (3.5*\xSize,-1.5*\ySize) circle (1.5pt);
\draw [->]  (3.5*\xSize,-1.5*\ySize) -- ++(0,-\xSize*\lineR);

\fill (3.5*\xSize,1.5*\ySize) circle (1.5pt);
\draw [->]  (3.5*\xSize,1.5*\ySize) -- ++(0,\xSize*\lineR);

\fill (-3.5*\xSize,-1.5*\ySize) circle (1.5pt);
\draw [->]  (-3.5*\xSize,-1.5*\ySize) -- ++(0,-\xSize*\lineR);

\fill (-3.5*\xSize,1.5*\ySize) circle (1.5pt);
\draw [->]  (-3.5*\xSize,1.5*\ySize) -- ++(0,\xSize*\lineR);


\fill (-1.5*\ySize, -3.5*\xSize) circle (1.5pt);
\draw [->]  (-1.5*\ySize, -3.5*\xSize) -- ++(-\xSize*\lineR, 0);

\fill (1.5*\ySize, -3.5*\xSize) circle (1.5pt);
\draw [->]  (1.5*\ySize, -3.5*\xSize) -- ++(\xSize*\lineR, 0);

\fill (-1.5*\ySize, 3.5*\xSize) circle (1.5pt);
\draw [->]  (-1.5*\ySize, 3.5*\xSize) -- ++(-\xSize*\lineR, 0);

\fill (1.5*\ySize, 3.5*\xSize) circle (1.5pt);
\draw [->]  (1.5*\ySize, 3.5*\xSize) -- ++(\xSize*\lineR, 0);

\pgfmathsetmacro{\lineR}{0.5}
\fill (2.75*\xSize,-2.75*\ySize) circle (1.5pt);
\draw [->]  (2.75*\xSize,-2.75*\ySize) -- ++(0,-\xSize*\lineR);

\fill (-2.75*\xSize,-2.75*\ySize) circle (1.5pt);
\draw [->]  (-2.75*\xSize,-2.75*\ySize) -- ++(-\xSize*\lineR,0);

\fill (2.75*\xSize,2.75*\ySize) circle (1.5pt);
\draw [->]  (2.75*\xSize,2.75*\ySize) -- ++(0,\xSize*\lineR);

\fill (-2.75*\xSize,2.75*\ySize) circle (1.5pt);
\draw [->]  (-2.75*\xSize,2.75*\ySize) -- ++(-\xSize*\lineR,0);

  \pgfmathsetmacro{\lineR}{1.9}
\fill (-3.8*\xSize,-0.5*\ySize) circle (1.5pt);
\draw [->]  (-3.8*\xSize,-0.5*\ySize) -- ++(0,-\xSize*\lineR);

\fill (-3.8*\xSize,0.5*\ySize) circle (1.5pt);
\draw [->]  (-3.8*\xSize,0.5*\ySize) -- ++(0,\xSize*\lineR);

\fill (3.8*\xSize,-0.5*\ySize) circle (1.5pt);
\draw [->]  (3.8*\xSize,-0.5*\ySize) -- ++(0,-\xSize*\lineR);

\fill (3.8*\xSize,0.5*\ySize) circle (1.5pt);
\draw [->]  (3.8*\xSize,0.5*\ySize) -- ++(0,\xSize*\lineR);

\fill (-0.5*\ySize, -3.8*\xSize) circle (1.5pt);
\draw [->]  (-0.5*\ySize, -3.8*\xSize) -- ++(-\xSize*\lineR, 0);

\fill (0.5*\ySize, -3.8*\xSize) circle (1.5pt);
\draw [->]  (0.5*\ySize, -3.8*\xSize) -- ++(\xSize*\lineR, 0);

\fill (-0.5*\ySize, 3.8*\xSize) circle (1.5pt);
\draw [->]  (-0.5*\ySize, 3.8*\xSize) -- ++(-\xSize*\lineR, 0);

\fill (0.5*\ySize, 3.8*\xSize) circle (1.5pt);
\draw [->]  (0.5*\ySize, 3.8*\xSize) -- ++(\xSize*\lineR, 0);

\draw[cyan,thick] plot[domain=0:360, smooth, samples=100] ({\R*cos(\x)}, {\R*sin(\x)});
\fill[cyan!10] plot[domain=0:360, smooth, samples=100] ({\R*cos(\x)}, {\R*sin(\x)}) -- cycle;
\draw[cyan,thick] plot[domain=0:360, smooth, samples=100] ({\sR*cos(\x)}, {\sR*sin(\x)});
\fill[white] plot[domain=0:360, smooth, samples=100] ({\sR*cos(\x)}, {\sR*sin(\x)}) -- cycle;

\pgfmathsetmacro{\lineR}{0.55}

\fill (-0.35*\ySize, -0.35*\xSize) circle (1.5pt);

\fill (-0.35*\ySize, 0.35*\xSize) circle (1.5pt);
\draw [red,->]    (-0.35*\ySize, 0.35*\xSize) -- ++(0, -\xSize*\lineR);

\fill (0.35*\ySize, -0.35*\xSize) circle (1.5pt);
\draw [red,->]   (0.35*\ySize, -0.35*\xSize) -- ++(-\xSize*\lineR, 0);

\fill (0.35*\ySize, 0.35*\xSize) circle (1.5pt);
\draw [red,->]   (0.35*\ySize, 0.35*\xSize) -- ++(-\xSize*\lineR, -\xSize*\lineR);

  \draw[step=\xSize, line width=0.03cm] (-3,-3) grid (4*\xSize,4*\xSize);

  \node at (-3.6*\xSize,-3.7*\ySize) {$\mathfrak{A}$};
  \node[text=cyan!90] at (1/2*\xSize,5/2*\ySize) {$\mathfrak{B}$};

\end{tikzpicture}
		\caption{Group formation} 
	\end{subfigure}
}
    \captionsetup{font=normalsize}
	\caption{A visual representation of the agglomeration algorithm for $\mathfrak{A}$ (white) on a 2D stationary cut-cell mesh with $\mathfrak{B}$ (cyan).
	The algorithm starts with identifying source cells (shown with dots) and searches target cells (shown with arrows) among non-source immediate neighbors. Then, it proceeds with the chain formation while automatically performing the level reduction. Lastly, a group is formed with the remaining source cells with the biggest cell chosen as the target. \label{fig:AggAlgorithm}}
\end{figure}
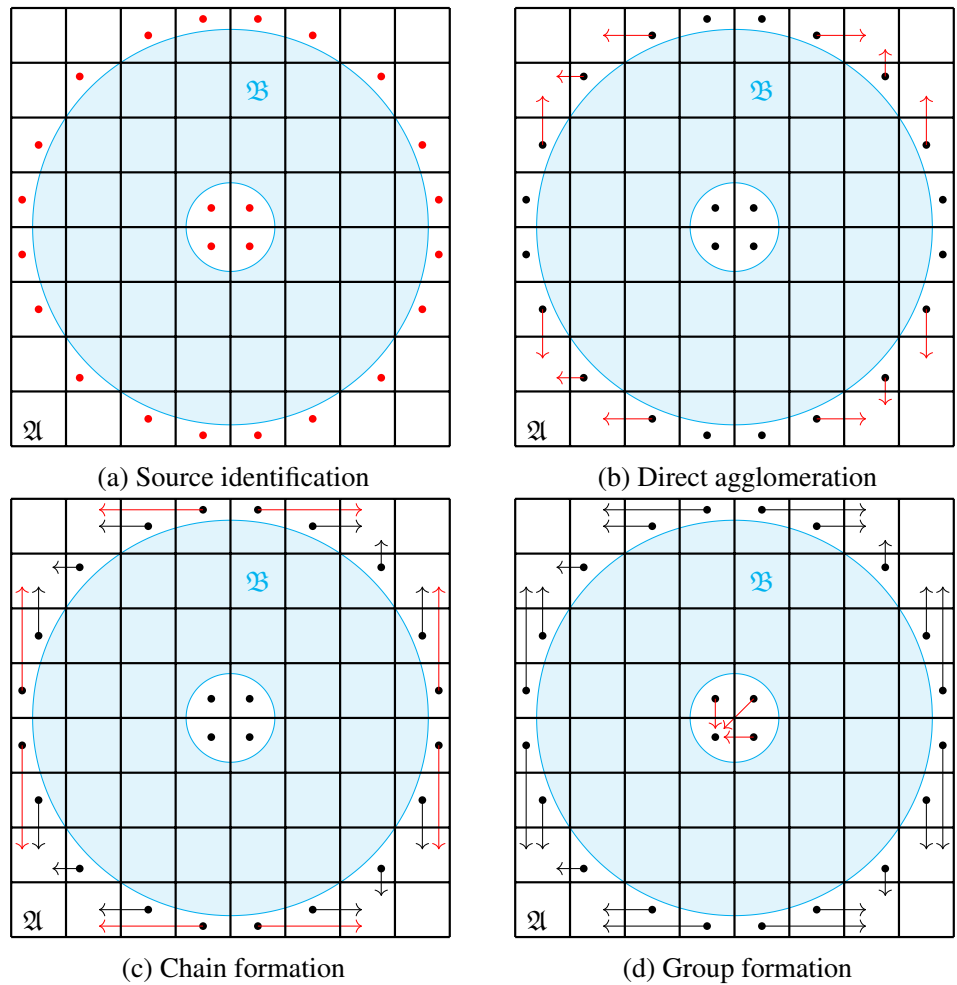

The agglomeration algorithm is implemented into the BoSSS (short for Bounded Support Spectral Solver) software package~\footnote{https://doi.org/10.5281/zenodo.8409677 \label{note1}} (see \cite{Kummer_2021_BoSSS}), which aims to solve partial differential equations by means of the discontinuous Galerkin (DG) methods. 
 For multiphase problems, BoSSS employs the eXtended Discontinuous Galerkin method (XDG)~\cite{kummer_extended_2017} and utilizes the cut cell approach to handle embedded geometries or interfaces. In the following, we present the algorithm of the proposed cell agglomeration strategy in pseudo codes. The source code of the actual implementation can be reached at the BoSSS GitHub repository~\footnote{DOI(TBD)}, which is written in C\# language using Message Passing Interface (MPI) for parallelization in distributed memory.


 The structure of the agglomeration algorithm is divided into several consecutive routines under \hyperref[subsec:SourceIdentification]{source identification} (Alg.~\ref{alg:SourceFilter}), \hyperref[subsec:TargetIdentification]{target identification} (Alg.~\ref{alg:DirectAgg},~\ref{alg:ChainFormation},~\ref{alg:GroupFormation}), and \hyperref[subsec:LevelDetermination]{level determination} (Alg.~\ref{alg:LevelDetermination}). It starts with the determination of local agglomeration sources by identifying the small-cut, vanishing, and newborn cells (Alg.~\ref{alg:SourceFilter}). 
 
 Subsequently, the algorithm attempts to create agglomeration mappings $\mathfrak{A}_{\mathrm{map},\mathfrak{s}}$ from these cells to suitable target cells via target identification routines. Initially, it searches target cells among the non-source immediate neighbors through the local direct target identification routine (Alg.~\ref{alg:DirectAgg}). If suitable targets cannot be found among the immediate neighbors, the algorithm then checks possible routes to form agglomeration chains with the already paired cells in the mapping (Alg.~\ref{alg:ChainFormation}). Lastly, it creates source-to-source agglomeration groups for cells without a route to a non-source cell via the group formation routine, with the largest cells designated as the target (Alg.~\ref{alg:GroupFormation}). 
 
 Once the mappings are created, the levels for inter-processor pairs are determined (Alg.~\ref{alg:LevelDetermination}), after which the agglomeration algebra is performed by following the sequential order (see §\ref{sec:aggLevels}). Figure~\ref{fig:AggAlgorithm} illustrates the different stages of the agglomeration algorithm. To ensure the uniqueness of the agglomeration mappings, an additional selection criterion based on the unique cell number is applied to all subroutines when all other selection criteria yield identical weights. For the sake of simplicity, this criterion is omitted in the presented algorithms. Furthermore, the algorithms are illustrated for local partitions, assuming the domain is already decomposed into processors.
 
 The underlying reason to divide the algorithm into sections lies in the communication requirements and the levels of agglomeration.
 The presented source and direct target identification routines correspond to more local and simple approaches, where no additional effort is needed in parallelization. This is because these routines only require information about neighboring cells of a partition, which is typically exchanged for flux calculations by default (i.e., ghost cells) in parallel numerical algorithms. 

 However, the information exchange between distant cells is often found redundant and resource-intensive in large simulations. This is particularly relevant in the context of DG methods, where high-order accuracy can be achieved by interacting primarily with zero-th-level neighbor cells. As a result, the formation of chains and determination of levels entails more expensive and less desirable processes, especially for inter-processor agglomerations, since they cannot be substituted by lower levels and call for parallel information exchange.
 Therefore, the algorithm first attempts to create a mapping with the direct agglomeration routine. Nevertheless, in complicated simulations, such as those with dynamic boundaries or multiple bodies, it is not always possible to create an appropriate mapping between immediate neighbors, making chain agglomerations unavoidable.




 \subsubsection*{Source identification} \label{subsec:SourceIdentification}
 \begin{algorithm}[ht!]
	\caption{Source cell filtration}
	 \label{alg:SourceFilter}
      \ttfamily \fontsize{9}{12}\selectfont{
    	 \begin{algorithmic}
    	 \Require A multi-step cut cell mesh $\mathfrak{K}_h^{X} = \{ \mathfrak{K}_h^{X,n-1}, \, \mathfrak{K}_h^{X,n} \}$
    	 \Ensure The lists of cells to be agglomerated $\mathfrak{K}_\mathrm{van,\mathfrak{s}},\mathfrak{K}_\mathrm{new,\mathfrak{s}}, \mathfrak{K}_{\mathrm{sml},\mathfrak{s}}, \mathfrak{K}_{\mathrm{src},\mathfrak{s}}$ 
    	 \ForEach {species $\mathfrak{s} \in \{\mathfrak{A}, \mathfrak{B}\}$}
    	 \State Initiate $\mathfrak{K}_\mathrm{van,\mathfrak{s}} := \{ K_{i,\mathfrak{s}} \;  | \; \sfrac(K_{i,\mathfrak{s}}^{n-1}) \geq 0$ and $\sfrac(K_{i,\mathfrak{s}}^{n}) = 0$\}
    	 \State Initiate $\mathfrak{K}_\mathrm{new,\mathfrak{s}} := \{ K_{i,\mathfrak{s}} \;  | \; \sfrac(K_{i,\mathfrak{s}}^{n-1}) = 0$ and $\sfrac(K_{i,\mathfrak{s}}^{n}) \geq 0$\}
    	 \State Initiate $\mathfrak{K}_\mathrm{sml,\mathfrak{s}} := \{ K_{i,\mathfrak{s}} \;  | \; \sfrac(K_{i,\mathfrak{s}}^{n-1}) < \alpha$ or $\sfrac(K_{i,\mathfrak{s}}^{n}) < \alpha$\}
    
    	 \State Initiate $\mathfrak{K}_{\mathrm{src},\mathfrak{s}} := \mathfrak{K}_\mathrm{van,\mathfrak{s}} \cup \mathfrak{K}_\mathrm{new,\mathfrak{s}} \cup \mathfrak{K}_{\mathrm{sml},\mathfrak{s}}$
    	 \State Communicate $\mathfrak{K}_{\mathrm{src},\mathfrak{s}}$ \Comment{exchange for cells at processor boundaries}
    	 \EndFor
    	 \end{algorithmic}
      }
 \end{algorithm}
 Finding the agglomeration sources follows a straightforward approach that identifies small-cut, vanishing, and newborn cells. Initially, the topological changes are determined and added to the lists $\mathfrak{K}_\mathrm{van,\mathfrak{s}}$ and $\mathfrak{K}_\mathrm{new,\mathfrak{s}}$. Subsequently, small-cut cells are detected and added to the source cells list $\mathfrak{K}_{\mathrm{sml},\mathfrak{s}}$. In the final stage, the lists are merged to $\mathfrak{K}_{\mathrm{src},\mathfrak{s}}$ and communicated across processors. Algorithm~\ref{alg:SourceFilter} provides the procedure for examining cells for vanishing, newborn, and small-cut cells in the case of the moving interface method. When using the splitting approach, the algorithm omits the labeling of vanishing cells into the source cell lists. In the case of a static mesh, the algorithm omits the labeling of both types of topological changes.

\subsubsection*{Target identification} \label{subsec:TargetIdentification}

\begin{algorithm}[ht!]
    \caption{Direct target identification}
    \label{alg:DirectAgg}
    \ttfamily \fontsize{9}{12}\selectfont{
        \begin{algorithmic}
            \Require {A cut cell mesh $\mathfrak{K}_h^X$\\
            \text{\phantom{Input}}The lists of cells to be agglomerated $\mathfrak{K}_{\mathrm{src},\mathfrak{s}}$} \hspace{0.8cm}
            \Ensure Agglomeration mappings $\mathfrak{A}_{\mathrm{map},\mathfrak{s}}$\\
            \text{\phantom{Output}}The lists of cells to be chain agglomerated $\mathfrak{K}_{\mathrm{chn},\mathfrak{s}}$
            \ForEach {species $\mathfrak{s} \in \{\mathfrak{A}, \mathfrak{B}\}$}
    		   \State Initiate $\mathfrak{A}_{\mathrm{map},\mathfrak{s}}$, $\mathfrak{K}_{\mathrm{chn},\mathfrak{s}}:= \{ \}$
                \State Flag non-source cells as candidate cells
    			\Comment (local and ghost cells)
                \ForEach {cell in $\mathfrak{K}_{\mathrm{src},\mathfrak{s}}$} 
                    \If {the cell has at least one shared edge with candidate cells}
                        \State {Choose the cell with the greatest fraction among them and add the pair into $\mathfrak{A}_{\mathrm{map},\mathfrak{s}}$}
                    \ElsIf {the cell has a shared edge with another source cell}
                        \State {Add the cell to the chain agglomeration cell list $\mathfrak{K}_{\mathrm{chn},\mathfrak{s}}$}
                    \EndIf
                \EndFor
            \EndFor
        \end{algorithmic}
    }
\end{algorithm}
To find agglomeration targets, the algorithm starts with the direct agglomeration procedure, wherein a suitable target cell $K_\mathrm{tar}$ is sought among the immediate neighbors of source cells, as shown in Alg.~\ref{alg:DirectAgg} with pseudo-codes. 
 Since all the candidates are adjacent cells in this stage, the distance creation is disregarded and the edges are weighted by the cell fraction of neighbor cells, calculated by $\sfrac()$ operator.
 Moreover, to ensure that the threshold $\alpha$ is passed and to prevent mapping to newborn or vanishing cells, the candidates are limited to non-source cells.
 If such a target is not found, a possible connection to another target is sought through the neighboring cells to form an agglomeration chain, and if so, the cell is forwarded to the chain-forming process shown in Alg.~\ref{alg:ChainFormation}.
 
 As a result, the direct agglomeration routine (i.e., Algorithm~\ref{alg:DirectAgg}) always produces pairs between adjacent cells, so it requires no additional information exchange to ensure connectivity.
 Moreover, the process is cycle-free because the targets are chosen among non-source cells, which are not mapped to another cell. These target cells serve as root cells in the agglomeration groups.

\begin{algorithm}[ht!]
    \caption{Chain formation}
    \label{alg:ChainFormation}
    \ttfamily \fontsize{9}{12}\selectfont{
        \begin{algorithmic}
    	\\
        \Require {A cut cell mesh $\mathfrak{K}_h^X$ \\
        \text{\phantom{Input}}The lists of cells to be chain agglomerated $\mathfrak{K}_{\mathrm{chn}, \mathfrak{s}}$\\
        \text{\phantom{Input}}Agglomeration mappings $\mathfrak{A}_{\mathrm{map},\mathfrak{s}}$}
        \Ensure {Updated agglomeration mappings $\mathfrak{A}_{\mathrm{map},\mathfrak{s}}$ \\
    	\text{\phantom{Output}}The list of remaining chain source cells $\mathfrak{K}_{\mathrm{chn}, \mathfrak{s}}$}
            \ForEach {species $\mathfrak{s} \in \{\mathfrak{A}, \mathfrak{B}\}$}
                 \Do
    			 	\State Communicate the local mapping $\mathfrak{A}_{\mathrm{map},\mathfrak{s}}$
    				\Comment{exchange boundary pairs with neighbor processors}
    				\State {Find all the edges between $\mathfrak{K}_{\mathrm{chn}, \mathfrak{s}}$ and the pairs in $\mathfrak{A}_{\mathrm{map},\mathfrak{s}}$}
                    \State {Weight the edges by the distance to the final targets and then the target sizes}
                    \State {Choose the best edge, identifying the chain source cell and the corresponding pair}
    				\If {the target of the pair is in this processor} 
    				\Comment (local and ghost cells)
    				\State Add the chain source and the target of the pair into $\mathfrak{A}_{\mathrm{map},\mathfrak{s}}$
    				\Else  				
    				\State Add the chain source and the source of the pair into $\mathfrak{A}_{\mathrm{map},\mathfrak{s}}$
    				\EndIf
                    \State {Remove the source cell from $\mathfrak{K}_{\mathrm{chn}, \mathfrak{s}}$}                
                \doWhile{ (any change in $\mathfrak{A}_{\mathrm{map},\mathfrak{s}}$)} 
    	\EndFor  
        \end{algorithmic}
    }
\end{algorithm}
On the other hand, chain agglomeration pairs are more intricate and result in higher levels of agglomeration through indirect neighborship. 
 Because of this, chain formation can produce cycles, so mappings should be selected with care to prevent them. 
 To circumvent such an issue, we restrict the linking of source cells in the chain formation to those that are already matched, either through direct agglomeration or previous iterations of chain agglomeration. Hence, chain agglomeration creates leaves of an existing agglomeration group and directs toward the targets of direct agglomeration pairs. 
 This ensures that a source cell is always mapped to an appropriate final target cell at the end. 
 However, in parallel executions, the instant state of the pairs on other processors is not explicitly recognizable, necessitating additional information exchange regarding neighbor cells outside of the partition. Algorithm~\ref{alg:ChainFormation} provides the proposed routine for chain agglomeration.
 
 Moreover, the possible agglomeration routes are weighted based on the properties of the final target cells, rather than of the neighboring source cells since all the source cells in an agglomeration group are going to be merged with the final target at the end. In this stage of agglomeration, edges are weighted with respect to distance and cell fraction as discussed in §\ref{sec:FormingChains}. To automatically reduce the agglomeration levels in the mappings, the pairs are established directly from source cells to final target cells (i.e., ($K_\mathrm{src},\mathfrak{A}_\mathrm{map}(K_\mathrm{ngb})$) ) unless they are in the different processors. Hence, the non-zero agglomeration levels narrow to the inter-processor pairs.
 

The chain agglomeration section presented in Alg.~\ref{alg:ChainFormation} is inspired by Kruskal's algorithm~\cite{kruskal1956shortest} to provide a minimum spanning forest.
 Usually, parallel variants of such algorithms try to minimize the sharing of information during the sorting stages. 
 However, the agglomeration algorithm differs from general approaches in that it works in a unilateral direction (i.e., from source cells to target cells), and cells are allowed to attach a chain only if the neighbor has already a designated target to prevent the formation of cycles.
 Hence, this procedure requires parallel communication between the partitions of the domain to supply the necessary information about the newly formed neighbor agglomeration pairs at the end of each iteration.
 
On the contrary, if cells were always agglomerated to a target cell within its partition, then no complicated information exchange would be required.
 Nonetheless, there is no guarantee that there will always be a suitable target, and even if there is, it will be the most preferable one. This is particularly relevant when using dynamic load balancing of the cut cells. Since cut cells require more computational effort than ordinary phase cells due to their irregular shapes, it is aimed to distribute them as evenly as possible across the processors and this eventually leads to neighboring source cells being located on different processors. 
Therefore, the proposed approach involves exchanging information about the mapped sources that lie in processor boundaries for agglomeration routines. This also ensures that the resulting agglomeration mappings remain independent of the number of processors. 
 Note that this information exchange only involves newly formed pairs at each iteration and is limited to the boundary processors that have adjacent source cells to those pairs.


\begin{algorithm}[htb!]
    \caption{Group formation}
    \label{alg:GroupFormation}
    \ttfamily \fontsize{9}{12}\selectfont{
        \begin{algorithmic}
    
    	\\
    	\Require {A cut cell mesh $\mathfrak{K}_h^X$ \\
    	\text{\phantom{Input}}Lists of source cells $\mathfrak{K}_{\mathrm{chn}, \mathfrak{s}},  \mathfrak{K}_\mathrm{van,\mathfrak{s}}, \mathfrak{K}_\mathrm{new,\mathfrak{s}} $\\
    	\text{\phantom{Input}}Agglomeration mappings $\mathfrak{A}_{\mathrm{map},\mathfrak{s}}$}
    	\Ensure {Updated agglomeration mappings $\mathfrak{A}_{\mathrm{map},\mathfrak{s}}$}
            \ForEach {species $\mathfrak{s} \in \{\mathfrak{A}, \mathfrak{B}\}$}
                 \State {Flag non-topological cells in $\mathfrak{K}_{\mathrm{chn}, \mathfrak{s}}$ as candidate cells}
                 \While {there is a candidate cell at least in one processor}
                    \State Select the biggest cell and communicate it across processors
    				\If {the cell is the biggest in all processors}
    					\State Assign the cell as the target and create an agglomeration group 
    					\Do
    					\Comment{create local group} 
    					\State Find all the neighbors of the agglomeration group
    					\State Add the neighbors pairing with the target cell into $\mathfrak{A}_{\mathrm{map},\mathfrak{s}}$
    					\State Remove the agglomerated cells from $\mathfrak{K}_{\mathrm{chn}, \mathfrak{s}}$
    					\doWhile{ (any change in the agglomeration group)} 
    				\EndIf			
    				\State {Call chain formation routine with $\mathfrak{K}_{\mathrm{chn}, \mathfrak{s}}$,$\mathfrak{A}_{\mathrm{map},\mathfrak{s}}$,$\mathfrak{K}_h^X$} 
                    \State {Re-flag non-topological cells in $\mathfrak{K}_{\mathrm{chn}, \mathfrak{s}}$ as candidate cells}
                \EndWhile
    	\EndFor  
        \end{algorithmic}
     }
\end{algorithm}

Until this point, the source cells that can be linked to a non-source cell should have already been mapped with the aforementioned routines. Nevertheless, there may still be source cells disconnected from any non-source cell. 
For these cells, local agglomeration groups consisting solely of source cells are established, and the cell with the highest fraction without requiring a topology change is designated as the final target in each group. Afterward, the formed agglomeration groups are parallelized by calling the chain formation routine in Alg.~\ref{alg:ChainFormation}, which can also be omitted depending on necessity. The respective procedure is shown in Algorithm~\ref{alg:GroupFormation}.

\subsubsection*{Level determination} \label{subsec:LevelDetermination}
\begin{algorithm}[ht!]
    \caption{Level determination}
    \label{alg:LevelDetermination}
    \ttfamily \fontsize{9}{12}\selectfont{
        \begin{algorithmic}
    	\\
    	\Require {A cut cell mesh $\mathfrak{K}_h^X$ and agglomeration mappings $\mathfrak{A}_{\mathrm{map},\mathfrak{s}}$}
    	\Ensure {Agglomerated cut cell mesh $\mathfrak{K}_h^{\mathrm{X},agg}$}
            \ForEach {species $\mathfrak{s} \in \{\mathfrak{A}, \mathfrak{B}\}$}
                 \State {Set levels to 0 for every source cell in $\mathfrak{A}_{\mathrm{map},\mathfrak{s}}$}
                \Do 
                    \ForEach {pair in $\mathfrak{A}_{\mathrm{map},\mathfrak{s}}$}
                        \While {there is a subsequent pair}
    						\If {the level is not bigger in the subsequent pair}
                        		\State {Assign one level higher to the source cell in the subsequent pair}
                    		\EndIf
                            \State Assign the subsequent pair as the current pair
                        \EndWhile
                    \EndFor
                \State Communicate levels for source cells
    			\Comment exchange for cells at processor boundaries
                \doWhile{ (any change in the levels)} 
    			\State Create $\mathfrak{K}_h^{\mathrm{X},agg}$ via $\mathfrak{A}_{\mathrm{map},\mathfrak{s}}$ following ascending order of the levels
    
    	\EndFor  
        \end{algorithmic}
    }
\end{algorithm}

Following the creation of agglomeration mappings, another routine shown in Algorithm~\ref{alg:LevelDetermination} tries to determine the levels of pairs in the inter-processor agglomeration chains. This routine traverses the agglomeration chains for each pair and then increases the levels for the subsequent pairs. Once the levels of agglomeration pairs are determined, the source cells are agglomerated to their target cell with an increasing order of levels. Hence, the sequential order of the agglomeration algebra is ensured and the operations are performed accordingly.



\section{Physical model equations and test cases} \label{sec:MethodologyPhysics}
In this study, we present several complex and dynamic test cases that are embedded in a fluid domain with a Cartesian background grid to showcase the presented agglomeration strategy. The subdomains $\mathfrak{A}$ and $\mathfrak{B}$ are reserved for fluid and embedded geometries, respectively. The fluid behavior is governed by the incompressible Navier-Stokes equation supplemented with the continuity equation given as:
\begin{align}
	\rho  \frac{\partial \vec{u}}{\partial t}+ \rho  \nabla \cdot (  \vec{u} \otimes \vec{u}) = \nabla \secondranktensor{\sigma}  + \vec{b} \ \ \text{in } \ \ \mathfrak{A}(t) \times (0,T), \\
	\nabla \cdot \vec{u} = 0 \ \ \text{in } \ \ \mathfrak{A}(t) \times (0,T) ,
\end{align}
where $\rho$ denotes the density of the fluid, $\vec{u}=(u,v,w)$ denotes the velocity vector, while $\secondranktensor{\sigma}$ represents surface forces acting on the fluid. The fluid is considered to be Newtonian and hence the stress tensor $\secondranktensor{\sigma}$ is modeled by:
\begin{equation}
	\secondranktensor{\sigma} = -P \secondranktensor{I} + \mu [\nabla \vec{u} + (\nabla \vec{u})^T],
\end{equation}
where $P$ is the pressure, $\secondranktensor{I}$ is the second rank unit tensor and $\mu$ is the dynamic viscosity. As a non-dimensional parameter, the Reynolds number is introduced and defined with respect to the maximum velocity $U$ and the characteristic length scale $L$ as $Re={\rho U L} / {\mu}$. In order to fulfill the Lady\u{z}enskaja-Babu\u{s}ka-Brezzi condition~\cite{Babuska_1973, Brezzi_1974}, the pressure $P$ is discretized in a space of lower degree space $p'=p-1$ for a velocity of degree $p$.
The fluid domain is initialized with uniform fields of zero value as initial conditions.:
\begin{equation}
	\begin{split}
\vec{u}(\vec{x},0) = \vec{0} \quad \text{on } \ \mathfrak{A}(t=0), \\ 
P(\vec{x},0) = 0 \quad \text{on } \ \mathfrak{A}(t=0).	
	\end{split}
	\label{eq:torusIC}
\end{equation} 		
The boundary conditions in test cases correspond to the pressure outlet on outer boundaries and a no-slip boundary condition on the interface:
\begin{align}
	\label{eq:rotTorusBC}
	\secondranktensor{\sigma} \cdot \vec{n}_{\Gamma_\mathrm{N}} &= 0 \ \ \ \text{on } \Gamma_\mathrm{N}=\partial \Omega, \\
	\vec{u} &= \vec{u_g} \  \text{on } \mathfrak{I},
\end{align}
where $\vec{u_g}$ denotes the velocity of the embedded geometry at the fluid boundary. 

The embedded geometries are described as rigid bodies by respective level-set functions (fluid: $\psi < 0$, geometry: $\psi > 0$) and coupled by means of the immersed boundary method (IBM) (for details see \cite{Muller_2017}). In this study, they are treated as a pseudo-solid phase implemented to create a complex, dynamic geometry for agglomeration, whereby internal physical phenomena in this region are not modeled as this study does not cover them. However, the presented agglomeration strategy is designed to be applicable for realistic multiphase interactions, provided with a cut-cell structure. In the following, several test cases inspired by similar applications are introduced.

 
\subsection{Vanishing sphere}
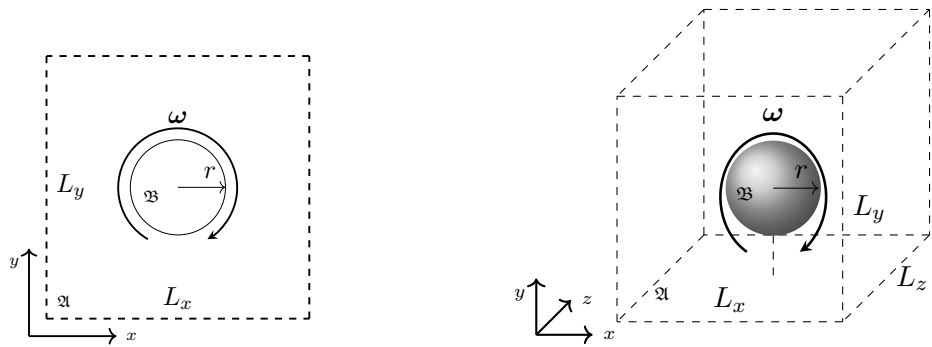
\begin{figure}[h!]
	\centering
	\begin{subfigure}[b]{0.45\textwidth}
		\centering
		\begin{tikzpicture}[scale=0.9]
\begin{axis}[
axis equal,
hide axis,
xmin=-0.2, xmax=2,
ymin=-0.2, ymax=2,
xlabel={$y$},
ylabel={$x$},
ticks=none,
]

\draw[black] (axis cs:0.75,0.75) circle (20pt);
\draw[->] (axis cs:0.75,0.75) -- ++(20pt,0pt) node[above left] {\(r\)};

\coordinate (center) at (axis cs:0.75,0.75);
\cercle{center}{25pt}{240}{-60};
\node[align=center] at (axis cs:0.75,1.15) {$\vec{\omega}$};

\node[align=center, font=\scriptsize] at (axis cs:0.1,0.1) {$\mathfrak{A}$};
\node[align=center, font=\scriptsize] at (axis cs:0.6,0.7) {$\mathfrak{B}$};

\draw[thick, dashed] (axis cs:0,0) -- (axis cs:1.5,0) node[midway, above=0pt] {$L_x$};
\draw[thick, dashed] (axis cs:1.5,0) -- (axis cs:1.5,1.5);
\draw[thick, dashed] (axis cs:1.5,1.5) -- (axis cs:0,1.5);
\draw[thick, dashed] (axis cs:0,1.5) -- (axis cs:0,0) node[midway, right=0pt] {$L_y$};

\pgfmathsetmacro{\xPos}{-0.1}
\pgfmathsetmacro{\yPos}{-0.1}

\draw[->, thick] (axis cs:\yPos,\xPos) -- (axis cs:\yPos+0.5,\xPos) node[right,font=\scriptsize] {$x$};
\draw[->, thick] (axis cs:\yPos,\xPos) -- (axis cs:\yPos,\xPos+0.5) node[below left,font=\scriptsize]{$y$};


\end{axis}
\end{tikzpicture}

	\end{subfigure}
	\begin{subfigure}[b]{0.45\textwidth}
		\centering
		\tdplotsetmaincoords{110}{70} 
\begin{tikzpicture}[scale=1.5]
\pgfplotsset{compat=1.18}
\renewcommand{\AxisRotator}[1][rotate=0]{%
    \tikz [x=24pt,y=20pt,line width=.2ex,-stealth,#1] \draw (0,0) arc (150:-150:1 and 1);%
}
\draw[dashed] (0,0,0) -- (0,2,0);
\draw[dashed] (0,2,0) -- (2,2,0);
\draw[dashed] (2,2,0) -- (2,0,0);
\draw[dashed] (2,0,0) -- (0,0,0);
\draw[dashed] (0,0,0) -- (0,0,2);
\draw[dashed] (2,0,0) -- (2,0,2) node[midway, right=0pt] {$L_z$};
\draw[dashed] (2,2,0) -- (2,2,2);
\draw[dashed] (0,2,0) -- (0,2,2);
\draw[dashed] (0,0,2) -- (0,2,2);
\draw[dashed] (0,2,2) -- (2,2,2);
\draw[dashed] (2,2,2) -- (2,0,2) node[midway, right=0pt] {$L_y$};
\draw[dashed] (2,0,2) -- (0,0,2) node[midway, above=0pt] {$L_x$};

\draw[dashed] (1,0.8,1)  -- ++(0,-0.8,0);
\node[align=center] at (1,1.45,1) {$\vec{\omega}$};

\shade[ball color=gray!50] (1,0.8,1) circle (12pt);


\node[font=\scriptsize] at (0.4,0.25,2) {$\mathfrak{A}$};
\node[font=\scriptsize] at (0.75,0.75,1) {$\mathfrak{B}$};

\draw[->] (1,0.8,1) -- ++(0.4pt,0pt,0pt) node[above left] {\(r\)};




    \pgfmathsetmacro{\xPos}{-0.6}
    \pgfmathsetmacro{\yPos}{0}
    \pgfmathsetmacro{\zPos}{2.3}
    \draw[->, thick] (\xPos,\yPos,\zPos) -- (\xPos+0.5,\yPos,\zPos) node[right,font=\scriptsize] {$x$};
    \draw[->, thick] (\xPos,\yPos, \zPos) -- (\xPos,\yPos+0.5,\zPos) node[below left,font=\scriptsize]{$y$};
    \draw[->, thick] (\xPos,\yPos, \zPos) -- (\xPos,\yPos,\zPos-0.8) node[right,font=\scriptsize]{$z$};

\node at (1,0.8,1)  {\AxisRotator[rotate=-270]};
\coordinate (center) at (1.1,0.8,1);
\end{tikzpicture}

	\end{subfigure}
	\caption{Schematic of the vanishing sphere test case. A sphere is placed at the origin of the domain with a pre-described rotation around its center with $\vec{u_g}(\vec{x},t) = \vec{\omega}(t) \times \vec{p}$. The radius of the sphere is defined dynamically as $r_\mathrm{s}(t)=(1-\frac{dr}{dt} t) r_{\mathrm{s},0}$. Pressure outlet boundary conditions are imposed on the outer boundary.  \label{fig:VanishingSphereSetup}}
\end{figure}
This test case examines a rotating sphere centered at the origin and characterized by a dynamically changing radius over time as illustrated in Figure~\ref{fig:VanishingSphereSetup}. It is described by the level set function given as:
\begin{equation}
	\psi(\vec{x},t) = r_\mathrm{s}(t)^2 - (\Vec{x}) \cdot (\Vec{x}),
\end{equation}
where $r_\mathrm{s}(t)$ denotes the time-dependent radius of the sphere. The radius is controlled by a constant shrinking rate $\frac{dr}{dt}$ and is defined as $r_\mathrm{s}(t)=(1-\frac{dr}{dt} t) r_{\mathrm{s},0}$ with respect to the initial radius $r_{\mathrm{s},0}$. Moreover, the sphere is subjected to a pre-described rotation around the origin by the angular velocity $\vec{\omega}$ leading to the velocity vector:
\begin{equation} \label{eq:rotation}
	\vec{u_g}(\vec{x},t) = \vec{\omega}(t) \times \vec{p} \quad \text{on } \ \mathfrak{B}(t), 
\end{equation}
where $\vec{p}$ is the position vector relative to the origin.
\subsection{Colliding spheres}
\begin{figure}[h!]
	\centering
	\begin{subfigure}[b]{0.45\textwidth}
		\centering
		\begin{tikzpicture}[scale=0.9]
  \begin{axis}[
    axis equal,
    hide axis,
    xmin=-0.4, xmax=4.1,
    ymin=-0.2, ymax=2,
    xlabel={$y$},
    ylabel={$x$},
    ticks=none,
  ]
  
  
  
  
  \draw[->] (axis cs:2.5,1) -- ++(15pt,0pt) node[above left] {\(r\)};
  \draw[->] (axis cs:1.5,1) -- ++(15pt,0pt) node[above left] {\(r\)};
  
  \draw[black] (axis cs:2.5,1) circle (15pt);
  \draw[black] (axis cs:1.5,1) circle (15pt);
  
  \node[align=center, font=\scriptsize] at (axis cs:0.1,0.1) {$\mathfrak{A}$};
  \node[align=center, font=\scriptsize] at (axis cs:1.3,1) {$\mathfrak{B}$};
  \node[align=center, font=\scriptsize] at (axis cs:2.3,1) {$\mathfrak{B}$};
  \draw[->, thick] (axis cs:1.4,1.45) -- ++(10pt,0pt) node[above left,font=\scriptsize] {$u$}; 
  \draw[->, thick] (axis cs:2.6,1.45) -- ++(-10pt,0pt) node[above right,font=\scriptsize] {$u$}; 

     

  \draw[thick, dashed] (axis cs:0,0) -- (axis cs:4,0) node[midway, above=0pt] {$L_x$};
  \draw[thick, dashed] (axis cs:4,0) -- (axis cs:4,2);
  \draw[thick, dashed] (axis cs:4,2) -- (axis cs:0,2);
  \draw[thick, dashed] (axis cs:0,2) -- (axis cs:0,0) node[midway, right=0pt] {$L_y$};
  
  \pgfmathsetmacro{\xPos}{-0.2}
  \pgfmathsetmacro{\yPos}{-0.2}
  
  \draw[->, thick] (axis cs:\yPos,\xPos) -- (axis cs:\yPos+0.5,\xPos) node[right,font=\scriptsize] {$x$};
  \draw[->, thick] (axis cs:\yPos,\xPos) -- (axis cs:\yPos,\xPos+0.5) node[below left,font=\scriptsize]{$y$};

  \end{axis}
  \end{tikzpicture}

	\end{subfigure}
	\begin{subfigure}[b]{0.45\textwidth}
		\centering
		\tdplotsetmaincoords{110}{70} 
\begin{tikzpicture}[scale=1.2]
\pgfplotsset{compat=1.18}

\draw[dashed] (0,0,0) -- (0,2,0);
\draw[dashed] (0,2,0) -- (4,2,0);
\draw[dashed] (4,2,0) -- (4,0,0);
\draw[dashed] (4,0,0) -- (0,0,0);
\draw[dashed] (0,0,0) -- (0,0,2);
\draw[dashed] (4,0,0) -- (4,0,2) node[midway, right=0pt] {$L_z$};
\draw[dashed] (4,2,0) -- (4,2,2);
\draw[dashed] (0,2,0) -- (0,2,2);
\draw[dashed] (0,0,2) -- (0,2,2);
\draw[dashed] (0,2,2) -- (4,2,2);
\draw[dashed] (4,2,2) -- (4,0,2) node[midway, right=-1pt] {$L_y$};
\draw[dashed] (4,0,2) -- (0,0,2) node[midway, above=0pt] {$L_x$};

\shade[ball color=gray!50] (2.5,0.8,1) circle (12pt);
\shade[ball color=gray!50] (1.4,0.8,1) circle (12pt);


\node[font=\scriptsize] at (0.4,0.25,2) {$\mathfrak{A}$};
\node[font=\scriptsize] at (1.15,0.75,1) {$\mathfrak{B}$};
\node[font=\scriptsize] at (2.25,0.75,1) {$\mathfrak{B}$};

\draw[->] (2.5,0.8,1) -- ++(0.4pt,0pt,0pt) node[above left] {\(r\)};
\draw[->] (1.4,0.8,1) -- ++(0.4pt,0pt,0pt) node[above left] {\(r\)};
\draw[->, thick] (1.25,1.3,1) -- ++(0.25pt,0pt,0pt) node[above left,font=\scriptsize] {$u$}; 
\draw[->, thick] (2.65,1.3,1) -- ++(-0.25pt,0pt,0pt) node[above right,font=\scriptsize] {$u$}; 




    \pgfmathsetmacro{\xPos}{-0.6}
    \pgfmathsetmacro{\yPos}{-0.1}
    \pgfmathsetmacro{\zPos}{2.1}
    \draw[->, thick] (\xPos,\yPos,\zPos) -- (\xPos+0.5,\yPos,\zPos) node[right,font=\scriptsize] {$x$};
    \draw[->, thick] (\xPos,\yPos, \zPos) -- (\xPos,\yPos+0.5,\zPos) node[below left,font=\scriptsize]{$y$};
    \draw[->, thick] (\xPos,\yPos, \zPos) -- (\xPos,\yPos,\zPos-0.8) node[right,font=\scriptsize]{$z$};


\end{tikzpicture}

	\end{subfigure}
	\caption{Schematic of colliding spheres test case. The spheres are placed into the domain with uniform mesh and pre-described with a velocity $u$ towards each other. Pressure outlet boundary conditions are imposed on the outer boundary. The radius for the spheres is given by $r_\mathrm{s}$.  \label{fig:CollidingSpheresSetup}}
\end{figure}
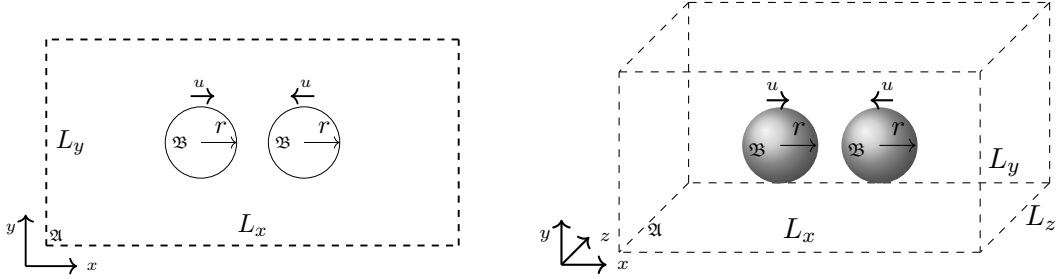
Another test case is introduced to examine the agglomeration strategy in the vicinity of merging, which often occurs during collisions of particles or droplets. To mimic the behavior of these isolated pockets in the fluid domain, two spheres are placed at a distance and subjected to velocities toward each other. Moreover, the spheres are allowed to overlap upon contact, thereby illustrating a prime instance of topology changes within the time-dependent domain. Accordingly, two spheres moving towards each other are defined in 2D and 3D spaces as depicted in Figure~\ref{fig:CollidingSpheresSetup}, which are described by the zero-th level sets of:
\begin{align}
	\psi(\vec{x},t) &= \max(\psi_\mathrm{L},\, \psi_\mathrm{R}),\\
	\psi_\mathrm{L}(\vec{x},t) &= r_\mathrm{s}^2 - (\Vec{x}-\Vec{x}_\mathrm{L}(t)) \cdot (\Vec{x}-\Vec{x}_\mathrm{L}(t)) ,\\
	\psi_\mathrm{R}(\vec{x},t) &= r_\mathrm{s}^2 - (\Vec{x}-\Vec{x}_\mathrm{R}(t)) \cdot (\Vec{x}-\Vec{x}_\mathrm{R}(t)),
\end{align}
where $r_\mathrm{s}$ denotes the radius of spheres, while $\Vec{x}_\mathrm{L}(t)$ and $\Vec{x}_\mathrm{R}(t)$ are the time-dependent position vectors describing the centers of spheres that are exposed to constant movement in the $x$-direction with a scalar velocity $u_\mathrm{s}$ as:
\begin{align}
	\Vec{x}_\mathrm{L}(t) &= (-1.5 r_\mathrm{s} + u_\mathrm{s} t, 0, 0), \\
	\Vec{x}_\mathrm{R}(t) &= (+1.5 r_\mathrm{s} - u_\mathrm{s} t, 0, 0).\
\end{align}
Hence, the velocity in the pseudo-solid phase (i.e., $\mathfrak{B}$) is determined as:
\begin{equation}
	\vec{u_g}(\vec{x},t) =
	 \begin{cases}
		\vec{u}_{g,\mathrm{L}} =  (u_\mathrm{s},0,0)  \quad \quad  \text{on } \ \vec{x} \, | \, \psi_\mathrm{L}(\vec{x},t)>0 \text{ and } \psi_\mathrm{R}(\vec{x},t)<0, \\
		\vec{u}_{g,\mathrm{R}} =  (-u_\mathrm{s},0,0) \quad \, \text{on } \ \vec{x} \, | \, \psi_\mathrm{L}(\vec{x},t)<0 \text{ and } \psi_\mathrm{R}(\vec{x},t)>0, \\
		\vec{u}_{g,\mathrm{LR}} =  \vec{0} \quad \quad \quad \quad \; \; \text{on } \ \vec{x} \, | \, \psi_\mathrm{L}(\vec{x},t) \geq 0 \text{ and } \psi_\mathrm{R}(\vec{x},t) \geq 0.

	\end{cases}
\end{equation}
\subsection{Rotating popcorn}
This test case introduces a popcorn-like shape inspired by the studies \cite{Burman_2015,Gurkan_2020} into the fluid domain described by the level set:
\begin{equation}
	\psi_{3D}(\vec{x}) = -(\sqrt{\Vec{x} \cdot \Vec{x}}-r_\mathrm{p}-\sum_{k=0}^{11} A \exp(- (\Vec{x} - \Vec{x}_{k,3D}) \cdot (\Vec{x} - \Vec{x}_{k,3D}) / \lambda^2 )),
\end{equation}
where 
\begin{equation}
	\Vec{x}_{k,3D}:=
	\begin{cases}
		\frac{r}{\sqrt{5}}(2 \cos (\frac{2k\pi}{5}),2 \sin (\frac{2k\pi}{5}),1) &\text{ if } \; \, \; \, \; \, 0 \leq k \leq 4,\\
		\frac{r}{\sqrt{5}}(2 \cos (\frac{(2(k-5)-1)\pi}{5}),2 \sin (\frac{(2(k-5)-1)\pi}{5}),-1) &\text{ if } \; \, \; \, \; \, 5 \leq k \leq 9, \\
		(0,0,r) &\text{ if } \; \, \; \, \; \, k = 10,\\
		(0,0,-r) &\text{ if } \; \, \; \, \; \, k = 11.
	\end{cases}
\end{equation} 
The parameters controlling the shape are given as $A=2$ and $\lambda=0.2$. Moreover, a 2D instance of the given level-set function is derived as:
\begin{equation}
	\psi_{2D}(\vec{x}) = -(\sqrt{\Vec{x} \cdot \Vec{x}}-r_\mathrm{p}-\sum_{k=0}^{4} A \exp(- (\Vec{x} - \Vec{x}_{k,2D}) \cdot (\Vec{x} - \Vec{x}_{k,2D}) / \lambda^2 )),
\end{equation}
with
\begin{equation}
	\Vec{x}_{k,2D}:=
		\frac{r}{\sqrt{5}}(2 \cos (\frac{2k\pi}{5}),2 \sin (\frac{2k\pi}{5})).
\end{equation} 
In line with the vanishing sphere, the computational domain remains consistent with the representation in Figure~\ref{fig:VanishingSphereSetup}. Similarly, a continuous rotation around the origin, as defined in Eq.\eqref{eq:rotation}, is applied, along with a corresponding motion of the geometry.
For the motion of the shape, the coordinate vector $\vec{x}$ is manipulated according to the rotation caused by $\vec{\omega}$ at time $t$ to account for movement. The resulting output $\vec{x'}(\vec{x},\vec{\omega},t)=(x',y',z')$ is then supplied to the given level-set equation.

\subsection{Rotating torus}
In this test case, a rotating torus is placed inside a fluid domain to simulate a moving boundary problem. The rotus is described by the zero-th level of the function:  
\begin{equation}
	\psi(\vec{x}) = -\sqrt{{\bigl(\sqrt{x^2+y^2}-r_{\mathrm{ma}}\bigr)}^2 + z^2} + r_{\mathrm{mi}},
\end{equation}
where $r_{\mathrm{ma}}$ and $r_{\mathrm{mi}}$ represent the major and minor radiuses, respectively, while $\vec{x}=(x,y,z)$ is the coordinate vector.
To create a more challenging and dynamic cut-cell structure, the torus is tilted in the $x$-axis with an angle of $\pi/4$ and then rotated with an angular velocity $\vec{\omega}$. 
The motion of the solid body is pre-described by the angular velocity $\vec{\omega}$ as given by Eq.~\ref{eq:rotation}. 



\section{Numerical results} \label{sec:Results}
In this section, the numerical results of the test cases described in the preceding section are presented. The test cases are simulated with varying agglomeration thresholds ($\alpha=0$, $0.1$, $0.2$, $0.3$, $0.4$, $0.5$), polynomial orders ($p=1, 2, 3$), and spatial dimensions ($D=2, 3$), except for the rotating torus case, which exclusively involves a 3D setup. To quantify the effect of agglomeration, we calculate the condition number of the mass and spatial operator matrices corresponding to the terms presented in Eq.~\ref{eq:weakFormulation}. Two condition numbers are computed according to the simulation scale: the global condition number, denoted by $\kappa_\mathrm{g}$, and the stencil condition number, denoted by $\kappa_\mathrm{s}$. Global condition numbers are obtained by providing the global mass and operator matrices of the overall system using MATLAB's built-in \textit{cond} command. The stencil condition numbers are computed for each cut cell by providing its stencil structure with direct neighbors using LAPACK's \textit{dgecon} routine. For cases with agglomeration threshold $\alpha=0$, it is observed that infinite values (`Inf') arise during the computations due to entries smaller than the machine epsilon in the discretizations corresponding to extremely small cuts. These entries are excluded from the results and their simulations are indicated with an asterisk in the provided tables.

Furthermore, the proposed agglomeration algorithm demonstrates consistency in condition numbers across different processor counts since it effectively exchanges information about cut cells in processor boundaries. Consequently, the number of processors becomes indifferent for the condition number study, so the simulations are performed with varying numbers of processors: 1 for low-degree 2D setups and 4, 8, 16, and 32 for more demanding test cases.



\subsection{Vanishing sphere}
\begin{figure}[ht!]
	\centering
	\includegraphics[width=0.8\textwidth]{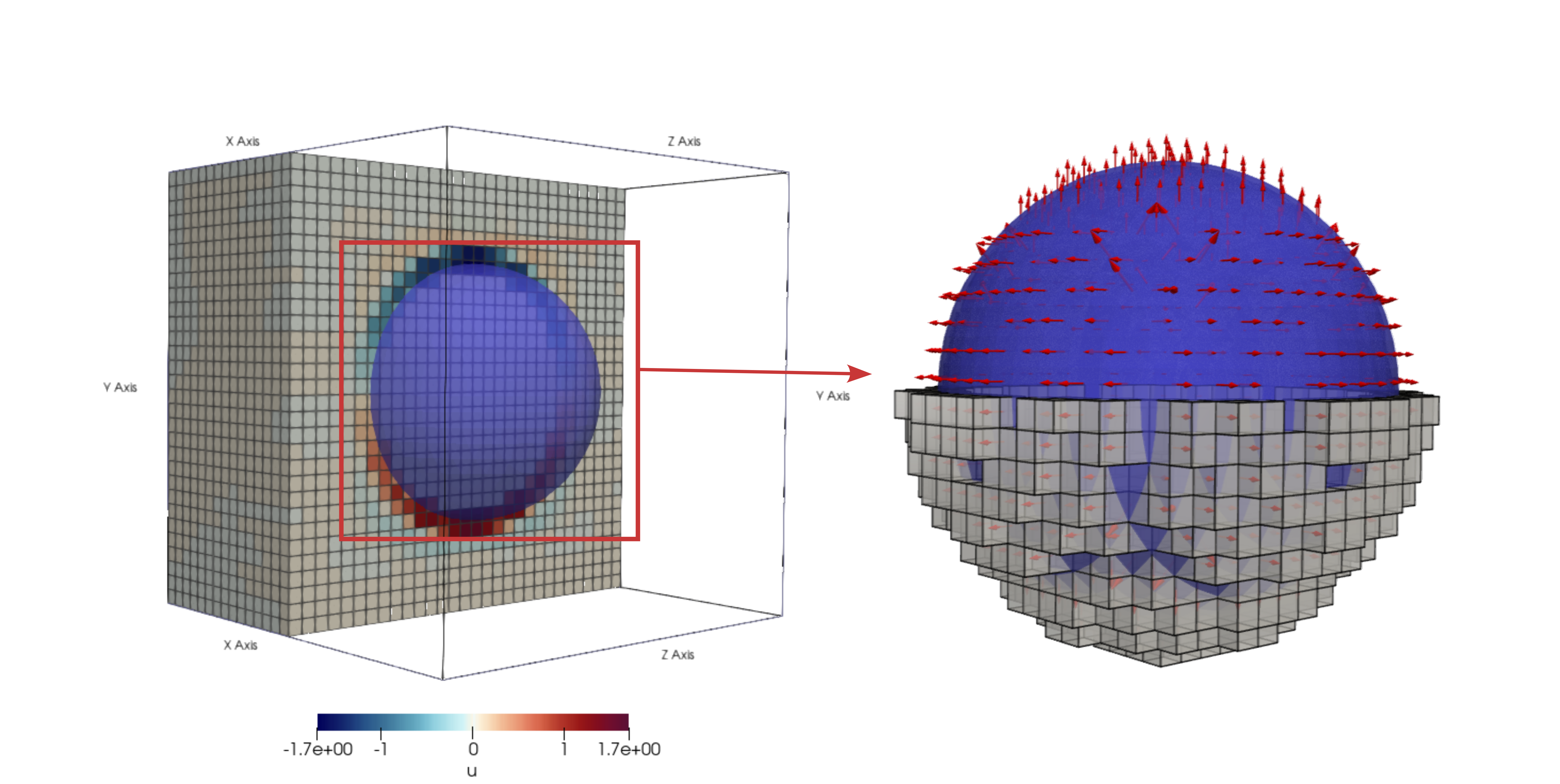}

	\caption{Instantaneous visualization of the 3D vanishing sphere test case with $p=1$ and $\alpha=0.5$ at the first time step. The embedded sphere is shown on the left with background mesh and the $x$-component of the velocity resulting from the rotation, while the agglomeration groups with cut cells are shown on the right. The interface between subdomains $\mathfrak{A}$ and $\mathfrak{B}$ is highlighted with blue color. The red glyphs illustrate the agglomeration mapping from source to target cells in $\mathfrak{A}$.  \label{fig:VanishingSphereInstant}}
 \end{figure}
 The vanishing sphere test case is simulated within a computational domain characterized by equal dimensions along each axis ($L_x=L_y=L_z$), as depicted in Figure~\ref{fig:VanishingSphereSetup}. The domain is discretized using a uniform Cartesian grid with $30 \times 30$ cells in 2D and $30 \times 30 \times 30$ cells in 3D. The sphere's initial radius is set as $r_{\mathrm{s},0} = 0.6 L_x$, and a constant shrinkage rate is defined such that the sphere fully vanishes at the end of the simulation, i.e. $\frac{dr}{dt}=1/T$. The sphere is subjected to an angular velocity around $z$-axis, resulting in $Re=2 \rho U r_{\mathrm{s}} / \mu$=500. The simulation period is determined as one revolution of the sphere around its axis, i.e. $T=2 \pi / |\Vec{\Omega}|$, while temporal discretization is achieved using an implicit Euler time scheme with 100 time steps, i.e. $dt=T/100$. 
\begin{figure}[htbp]
	\centering
	\begin{subfigure}[b]{0.35\textwidth}
		\centering
		\resizebox{\linewidth}{!}{\begin{tikzpicture}
    \begin{axis}[
        ylabel={$\% $},
        xmin=0, xmax=2*pi,
        ymin=-1, ymax=8,
        grid=both,
        legend pos=north east,
        cycle list name=exotic,
        legend columns=3,
        xtick={0, pi/2, pi, 3*pi/2, 2*pi}, 
        xticklabels={0, $\frac{T}{4}$, $\frac{T}{2}$, $\frac{3T}{4}$, $T$}, 
    ]
    
    \pgfmathsetmacro{\TotalCellNumber}{9} 

    \addplot table[x expr={\coordindex /80 * 2 * pi}, y expr={\thisrowno{0}/ \TotalCellNumber}] {data/2D_VanishingSphere_k1_Agg0/AggNoForCombinedA.csv};
    \addlegendentry{$\alpha = 0.0$};

   \addplot table[x expr={\coordindex /80 * 2 * pi}, y expr={\thisrowno{0}/ \TotalCellNumber}] {data/2D_VanishingSphere_k1_Agg0.1/AggNoForCombinedA.csv};
   \addlegendentry{$\alpha = 0.1$};

   \addplot table[x expr={\coordindex /80 * 2 * pi}, y expr={\thisrowno{0}/ \TotalCellNumber}] {data/2D_VanishingSphere_k1_Agg0.2/AggNoForCombinedA.csv};
   \addlegendentry{$\alpha = 0.2$};

   \addplot table[x expr={\coordindex /80 * 2 * pi}, y expr={\thisrowno{0}/ \TotalCellNumber}] {data/2D_VanishingSphere_k1_Agg0.3/AggNoForCombinedA.csv};
   \addlegendentry{$\alpha = 0.3$};

   \addplot table[x expr={\coordindex /80 * 2 * pi}, y expr={\thisrowno{0}/ \TotalCellNumber}] {data/2D_VanishingSphere_k1_Agg0.4/AggNoForCombinedA.csv};
   \addlegendentry{$\alpha = 0.4$};

   \addplot table[x expr={\coordindex /80 * 2 * pi}, y expr={\thisrowno{0}/ \TotalCellNumber}] {data/2D_VanishingSphere_k1_Agg0.5/AggNoForCombinedA.csv};
   \addlegendentry{$\alpha = 0.5$};

    \end{axis}
\end{tikzpicture}}
		\label{fig:VanishingSphereAggNo2D}
	\end{subfigure}
	\begin{subfigure}[b]{0.35\textwidth}
		\centering
		\resizebox{\linewidth}{!}{\begin{tikzpicture}
    \begin{axis}[
        ylabel={$\% $},
        xmin=0, xmax=2*pi,
        ymin=-0.5, ymax=4,
        grid=both,
        legend pos=north east,
        cycle list name=exotic,
        legend columns=3,
        xtick={0, pi/2, pi, 3*pi/2, 2*pi}, 
        xticklabels={0, $\frac{T}{4}$, $\frac{T}{2}$, $\frac{3T}{4}$, $T$}, 
    ]
    
    \pgfmathsetmacro{\TotalCellNumber}{270} 
    \addplot table[x expr={\coordindex /80* 2 * pi}, y expr={\thisrowno{0}/ \TotalCellNumber}]
 {data/VanishingSphere_k1_Agg0/AggNoForCombinedA.csv};
    \addlegendentry{$\alpha = 0.0$};

   \addplot table[x expr={\coordindex /80* 2 * pi}, y expr={\thisrowno{0}/ \TotalCellNumber}]
 {data/VanishingSphere_k1_Agg0.1/AggNoForCombinedA.csv};
   \addlegendentry{$\alpha = 0.1$};

   \addplot table[x expr={\coordindex /80* 2 * pi}, y expr={\thisrowno{0}/ \TotalCellNumber}]
 {data/VanishingSphere_k1_Agg0.2/AggNoForCombinedA.csv};
   \addlegendentry{$\alpha = 0.2$};

   \addplot table[x expr={\coordindex /80* 2 * pi}, y expr={\thisrowno{0}/ \TotalCellNumber}]
 {data/VanishingSphere_k1_Agg0.3/AggNoForCombinedA.csv};
   \addlegendentry{$\alpha = 0.3$};

   \addplot table[x expr={\coordindex /80* 2 * pi}, y expr={\thisrowno{0}/ \TotalCellNumber}]
 {data/VanishingSphere_k1_Agg0.4/AggNoForCombinedA.csv};
   \addlegendentry{$\alpha = 0.4$};

   \addplot table[x expr={\coordindex /80* 2 * pi}, y expr={\thisrowno{0}/ \TotalCellNumber}]
 {data/VanishingSphere_k1_Agg0.5/AggNoForCombinedA.csv};
   \addlegendentry{$\alpha = 0.5$};

    \end{axis}
\end{tikzpicture}}
		\label{fig:VanishingSphereAggNo3D}
	\end{subfigure}
	\caption{Percentage of the agglomerated cut cells for the fluid phase ($\mathfrak{A}$) during the vanishing sphere test case in 2D (left) and 3D (right). \label{fig:VanishingSphereAggNo}}
 \end{figure}
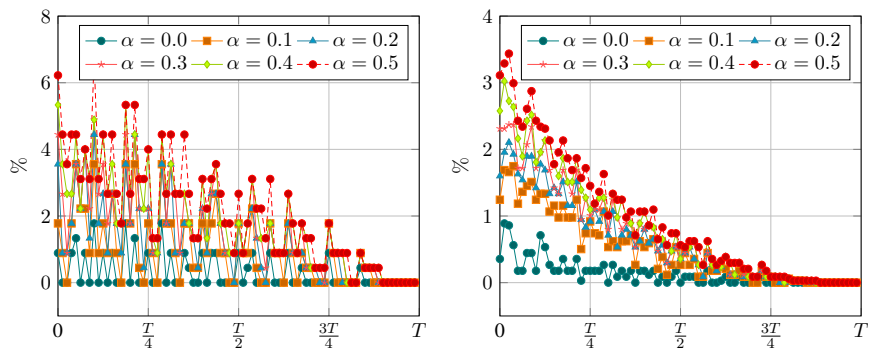

 Figure~\ref{fig:VanishingSphereInstant} illustrates an example of the instantaneous simulation results for the chosen 3D simulation configuration with $\alpha=0.5$. The agglomeration groups are demonstrated with the respective mapping via glyphs directing from the source cells to their paired target cells. Furthermore, Figure~\ref{fig:VanishingSphereAggNo} depicts the corresponding agglomerated cell numbers for the subdomain $\mathfrak{A}$ over the course of simulations, which see a constant decrease due to the shrinkage of the sphere. The variations of the agglomerated cell numbers in the cases with $\alpha=0$ correspond to the topological changes occurring occasionally and are the same for every simulation, regardless of the small-cut threshold.

\begin{figure}[ht!]
	\centering
	\begin{subfigure}[b]{0.35\textwidth}
		\centering
		\resizebox{\linewidth}{!}{\begin{tikzpicture}
    \begin{axis}[
        ylabel={$\kappa_\mathrm{g}$},
        xmin=0, xmax=2*pi,
        ymin=1e3, ymax=1e11,
        ymode=log,
        cycle list name=exotic,
        grid=none,
        legend columns=3,
        xtick={0, pi/2, pi, 3*pi/2, 2*pi}, 
        xticklabels={0, $\frac{T}{4}$, $\frac{T}{2}$, $\frac{3T}{4}$, $T$}, 
    ]
    \addplot table[col sep=comma, /pgf/number format/use comma, x expr=\thisrow{File Index}/80*2*pi] {data/2D_VanishingSphere_k1_Agg0/MassMtxCondNo_data.csv};
    \addlegendentry{$\alpha = 0.0$};

   \addplot table[col sep=comma, /pgf/number format/use comma, x expr=\thisrow{File Index}/80*2*pi] {data/2D_VanishingSphere_k1_Agg0.1/MassMtxCondNo_data.csv}; \addlegendentry{$\alpha = 0.1$};

   \addplot table[col sep=comma, /pgf/number format/use comma, x expr=\thisrow{File Index}/80*2*pi] {data/2D_VanishingSphere_k1_Agg0.2/MassMtxCondNo_data.csv};
   \addlegendentry{$\alpha = 0.2$};

   \addplot table[col sep=comma, /pgf/number format/use comma, x expr=\thisrow{File Index}/80*2*pi] {data/2D_VanishingSphere_k1_Agg0.3/MassMtxCondNo_data.csv};
   \addlegendentry{$\alpha = 0.3$};

   \addplot table[col sep=comma, /pgf/number format/use comma, x expr=\thisrow{File Index}/80*2*pi] {data/2D_VanishingSphere_k1_Agg0.4/MassMtxCondNo_data.csv};
   \addlegendentry{$\alpha = 0.4$};

   \addplot table[col sep=comma, /pgf/number format/use comma, x expr=\thisrow{File Index}/80*2*pi] {data/2D_VanishingSphere_k1_Agg0.5/MassMtxCondNo_data.csv};
   \addlegendentry{$\alpha = 0.5$};

    \end{axis}
\end{tikzpicture}}
		\label{fig:VanishingSphere2DMassCondk1}
	\end{subfigure}
	\begin{subfigure}[b]{0.35\textwidth}
		\centering
		\resizebox{\linewidth}{!}{\begin{tikzpicture}
    \begin{axis}[
       ylabel={$\kappa_\mathrm{g}$},
        xmin=0, xmax=2*pi,
        ymode=log,
        cycle list name=exotic,
        grid=none,
        legend columns=3,
        legend pos=south east,
        xtick={0, pi/2, pi, 3*pi/2, 2*pi}, 
        xticklabels={0, $\frac{T}{4}$, $\frac{T}{2}$, $\frac{3T}{4}$, $T$}, 
    ]

     \addplot table[col sep=comma, /pgf/number format/use comma, x expr=\thisrow{File Index}/80*2*pi] {data/2D_VanishingSphere_k1_Agg0/TotCondNo-Vars0.1.2_data.csv};
     \addlegendentry{$\alpha = 0.0$};

   \addplot table[col sep=comma, /pgf/number format/use comma, x expr=\thisrow{File Index}/80*2*pi] {data/2D_VanishingSphere_k1_Agg0.1/TotCondNo-Vars0.1.2_data.csv};
   \addlegendentry{$\alpha = 0.1$};

   \addplot table[col sep=comma, /pgf/number format/use comma, x expr=\thisrow{File Index}/80*2*pi] {data/2D_VanishingSphere_k1_Agg0.2/TotCondNo-Vars0.1.2_data.csv};
   \addlegendentry{$\alpha = 0.2$};

   \addplot table[col sep=comma, /pgf/number format/use comma, x expr=\thisrow{File Index}/80*2*pi] {data/2D_VanishingSphere_k1_Agg0.3/TotCondNo-Vars0.1.2_data.csv};
   \addlegendentry{$\alpha = 0.3$};

   \addplot table[col sep=comma, /pgf/number format/use comma, x expr=\thisrow{File Index}/80*2*pi] {data/2D_VanishingSphere_k1_Agg0.4/TotCondNo-Vars0.1.2_data.csv};
   \addlegendentry{$\alpha = 0.4$};

   \addplot table[col sep=comma, /pgf/number format/use comma, x expr=\thisrow{File Index}/80*2*pi] {data/2D_VanishingSphere_k1_Agg0.5/TotCondNo-Vars0.1.2_data.csv};
   \addlegendentry{$\alpha = 0.5$};

    \end{axis}
\end{tikzpicture}}
		\label{fig:VanishingSphere2DOpCondk1}
	\end{subfigure}
	\begin{subfigure}[b]{0.35\textwidth}
		\centering
		\resizebox{\linewidth}{!}{\begin{tikzpicture}
    \begin{axis}[
        ylabel={$\kappa_\mathrm{g}$},
        xmin=0, xmax=2*pi,
        ymin=1e3, ymax=1e11,
        ymode=log,
        cycle list name=exotic,
        grid=none,
        legend columns=3,
        xtick={0, pi/2, pi, 3*pi/2, 2*pi}, 
        xticklabels={0, $\frac{T}{4}$, $\frac{T}{2}$, $\frac{3T}{4}$, $T$}, 
    ]
    \addplot table[col sep=comma, /pgf/number format/use comma, x expr=\thisrow{File Index}/80*2*pi] {data/VanishingSphere_k1_Agg0/MassMtxCondNo_data.csv};
    \addlegendentry{$\alpha = 0.0$};

   \addplot table[col sep=comma, /pgf/number format/use comma, x expr=\thisrow{File Index}/80*2*pi] {data/VanishingSphere_k1_Agg0.1/MassMtxCondNo_data.csv}; \addlegendentry{$\alpha = 0.1$};

   \addplot table[col sep=comma, /pgf/number format/use comma, x expr=\thisrow{File Index}/80*2*pi] {data/VanishingSphere_k1_Agg0.2/MassMtxCondNo_data.csv};
   \addlegendentry{$\alpha = 0.2$};

   \addplot table[col sep=comma, /pgf/number format/use comma, x expr=\thisrow{File Index}/80*2*pi] {data/VanishingSphere_k1_Agg0.3/MassMtxCondNo_data.csv};
   \addlegendentry{$\alpha = 0.3$};

   \addplot table[col sep=comma, /pgf/number format/use comma, x expr=\thisrow{File Index}/80*2*pi] {data/VanishingSphere_k1_Agg0.4/MassMtxCondNo_data.csv};
   \addlegendentry{$\alpha = 0.4$};

   \addplot table[col sep=comma, /pgf/number format/use comma, x expr=\thisrow{File Index}/80*2*pi] {data/VanishingSphere_k1_Agg0.5/MassMtxCondNo_data.csv};
   \addlegendentry{$\alpha = 0.5$};

    \end{axis}
\end{tikzpicture}}
		\label{fig:VanishingSphereMassCondk1}
	\end{subfigure}
	\begin{subfigure}[b]{0.35\textwidth}
		\centering
		\resizebox{\linewidth}{!}{\begin{tikzpicture}
    \begin{axis}[
        ylabel={$\kappa_\mathrm{g}$},
        ymax=10^3.7,
        xmin=0, xmax=2*pi,
        ymode=log,
        cycle list name=exotic,
        grid=none,
        legend columns=3,
        xtick={0, pi/2, pi, 3*pi/2, 2*pi}, 
        xticklabels={0, $\frac{T}{4}$, $\frac{T}{2}$, $\frac{3T}{4}$, $T$}, 
    ]

     \addplot table[col sep=comma, /pgf/number format/use comma, x expr=\thisrow{File Index}/80*2*pi] {data/VanishingSphere_k1_Agg0/TotCondNo-Vars0.1.2.3_data.csv};
     \addlegendentry{$\alpha = 0.0$};

   \addplot table[col sep=comma, /pgf/number format/use comma, x expr=\thisrow{File Index}/80*2*pi] {data/VanishingSphere_k1_Agg0.1/TotCondNo-Vars0.1.2.3_data.csv};
   \addlegendentry{$\alpha = 0.1$};

   \addplot table[col sep=comma, /pgf/number format/use comma, x expr=\thisrow{File Index}/80*2*pi] {data/VanishingSphere_k1_Agg0.2/TotCondNo-Vars0.1.2.3_data.csv};
   \addlegendentry{$\alpha = 0.2$};

   \addplot table[col sep=comma, /pgf/number format/use comma, x expr=\thisrow{File Index}/80*2*pi] {data/VanishingSphere_k1_Agg0.3/TotCondNo-Vars0.1.2.3_data.csv};
   \addlegendentry{$\alpha = 0.3$};

   \addplot table[col sep=comma, /pgf/number format/use comma, x expr=\thisrow{File Index}/80*2*pi] {data/VanishingSphere_k1_Agg0.4/TotCondNo-Vars0.1.2.3_data.csv};
   \addlegendentry{$\alpha = 0.4$};

   \addplot table[col sep=comma, /pgf/number format/use comma, x expr=\thisrow{File Index}/80*2*pi] {data/VanishingSphere_k1_Agg0.5/TotCondNo-Vars0.1.2.3_data.csv};
   \addlegendentry{$\alpha = 0.5$};

    \end{axis}
\end{tikzpicture}}
		\label{fig:VanishingSphereOpCondk1}
	\end{subfigure}
		
	\caption{Global condition numbers for the mass (left) and operator (right) matrices in the transient 2D (top) and 3D (bottom) vanishing sphere test cases with polynomial degree $p=1$ and varying agglomeration threshold $\alpha$. Note that `Inf' entries computed for the mass matrix, which correspond to the jumps at operator matrices, are excluded in the figure. \label{fig:VanishingSphereCondk1}}
 \end{figure}
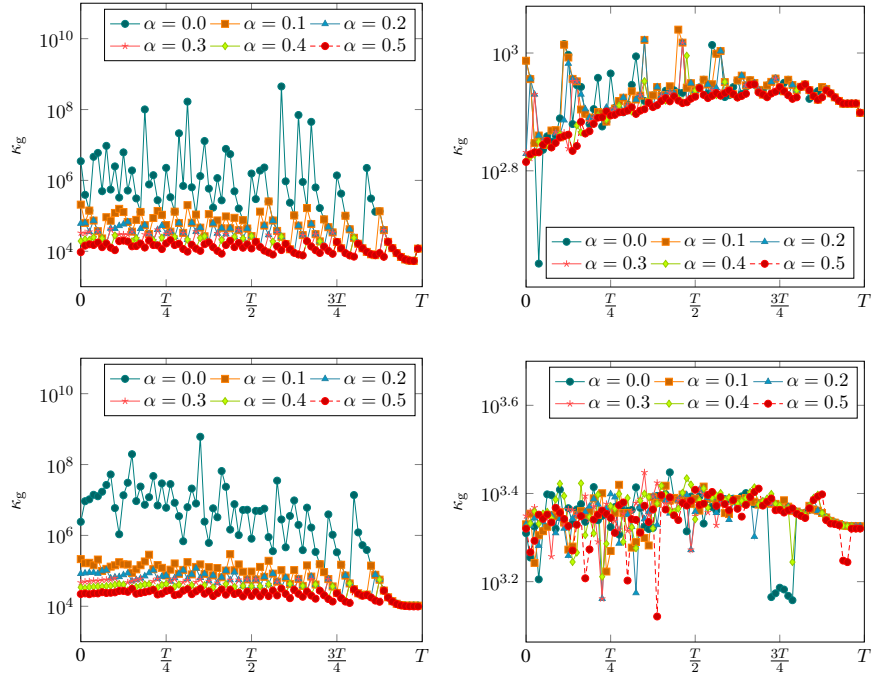

The global condition numbers for the simulations with $p=1$ in 2D and 3D are presented in Figure~\ref{fig:VanishingSphereCondk1} in time series since the simulation considers a dynamically changing shape over time. As can be seen in the figure, the condition numbers exhibit slight variations corresponding to changes in the shape and the associated cut cell structure. Notably, the most significant fluctuations occur in the instances where no agglomeration is present for small-cut cells, denoted by $\alpha=0$, showcasing remarkable jumps over time. This disparity becomes more pronounced, particularly in mass matrices, owing to their orthonormalization. As the number of agglomerated cells gradually decreases due to the shrinkage of the sphere's domain ($\mathfrak{B}$), the effect of agglomeration disappears and the condition numbers match with each other toward the end of the simulations, resulting in the same minimum values. These qualitative trends are observed for all configurations of different polynomial degrees, with only discrepancies in the magnitudes of the condition numbers. Therefore, the maximum condition numbers for the simulations are summarized in Table~\ref{tab:VanishingResults}. 

\begin{table}[htb] 
    \centering
    \caption{Maximum condition numbers for the vanishing sphere test case}
    \label{tab:VanishingResults}
	\small
	\begin{tabular}{l|c:c|c:c|c:c}  
		\toprule
			\multicolumn{1}{c}{} & \multicolumn{6}{c}{2D global numbers($\kappa_\mathrm{g}$)}  \\ 
	  \cmidrule(rl){2-7}  
	
	\multicolumn{1}{c}{} & \multicolumn{2}{c}{$p=1$}  & \multicolumn{2}{c}{$p=2$} & \multicolumn{2}{c}{$p=3$} \\ 
	
	$\alpha$ & Mass & Op  & Mass & Op & Mass & Op  \\ \hline
	
	0\phantom{.1}  \;  & 4.48e8 & 1.04e3 & 1.06e9* & 4.36e20 & 5.85e8* & 7.12e19   \\ \hline
	0.1 \;  & 2.55e5 & 1.10e3 & 2.91e6 & 5.80e3 & 1.41e7 & 1.39e4 \\ \hline	
	0.2 \;  & 7.37e4 & 1.05e3 & 8.82e5 & 4.91e3 & 4.05e6 & 1.14e4 \\ \hline
	0.3 \;  & 4.25e4 & 1.04e3 & 4.87e5 & 4.36e3 & 2.48e6 & 1.01e4 \\ \hline
	0.4 \;  & 2.80e4 & 9.89e2 & 3.17e5 & 4.09e3 & 1.37e6 & 8.99e3 \\ \hline
	0.5 \;  & 2.06e4 & 8.85e2 & 2.31e5 & 3.90e3 & 9.87e5 & 8.99e3 \\ \hline
	
		\end{tabular}
		
			\begin{tabular}{l|c:c|c:c|c:c}  
		\multicolumn{7}{c}{} \\
			\multicolumn{1}{c}{} & \multicolumn{6}{c}{3D stencil numbers($\kappa_\mathrm{s}$)}  \\ 
	  \cmidrule(rl){2-7}  
	
	0\phantom{.1}  \;  & 6.00e8 & 7.88e1 & 2.29e10 & 3.95e2
	 & 4.68e13* & 2.05e22   \\ \hline
	0.1 \;  & 2.71e5 & 8.31e1 & 4.03e6 & 3.90e2 & 2.41e7 & 1.40e3 \\ \hline
	0.2 \;  & 1.07e5 & 7.13e1 & 1.33e6 & 3.52e2 & 6.14e6 & 1.33e3 \\ \hline
	0.3 \;  & 5.71e4 & 6.40e1 & 7.14e5 & 3.05e2 & 3.61e6 & 1.17e3 \\ \hline
	0.4 \;  & 4.11e4 & 6.36e1 & 8.29e5 & 2.77e2 & 3.61e6 & 1.08e3 \\ \hline
	0.5 \;  & 3.16e4 & 6.32e1 & 4.11e5 & 2.65e2 & 1.45e6 & 9.59e2 \\ \hline

	
		\end{tabular}
    \begin{tabular}{l}
		\footnotesize *Maximum value excluding `Inf' entries     \phantom{afasdfadadfasdsfasdfasdfassdfdfasdfasdf}
	\end{tabular}
\end{table}
It is evident from the numerical values, presented in Fig.~\ref{fig:VanishingSphereCondk1} and Tab.~\ref{tab:VanishingResults}, that the condition numbers are strongly dependent on the presence of agglomeration. Specifically, cases with $\alpha=0.0$ exhibit considerably higher condition numbers for the mass matrices, often surpassing additional three orders of magnitude. Remarkably, these disparities are more accentuated in the maximum condition numbers of operator matrices, occasionally reaching magnitudes such as $7.12 \times 10^{19}$ or $7.12 \times 10^{20}$ during simulations. This situation is attributed to the presence of infinite entries in the mass matrices corresponding to these extreme values, which are excluded in the table. These specific simulations are denoted with an asterisk.
Regarding the impact of the agglomeration threshold, it is notable that an increased agglomeration threshold corresponds to diminished condition numbers. A clear trend emerges for the mass matrices, indicating a consistent decrease in condition numbers as more cut cells are agglomerated as shown in Fig.~\ref{fig:VanishingSphereAggNo}. These mass matrices exhibit an approximate halving of their condition numbers from $\alpha=0.1$ to $\alpha=0.2$. This trend persists until the agglomeration threshold $\alpha=0.3$, although the numerical values then exhibit a gradual shift with higher thresholds. For operator matrices, the condition numbers do not yield a particular trend with $\alpha$ and display similar values with each other as a result of the underlying orthonormalization of the bases. If the bases were not orthonormalized, the influence of agglomeration would be reflected in the operator matrices, while the values for the mass matrices would be lower, but the overall system would have higher condition numbers.

\subsection{Colliding spheres}
\begin{figure}[ht!]
	\centering
		\centering
		\includegraphics[width=0.75\textwidth]{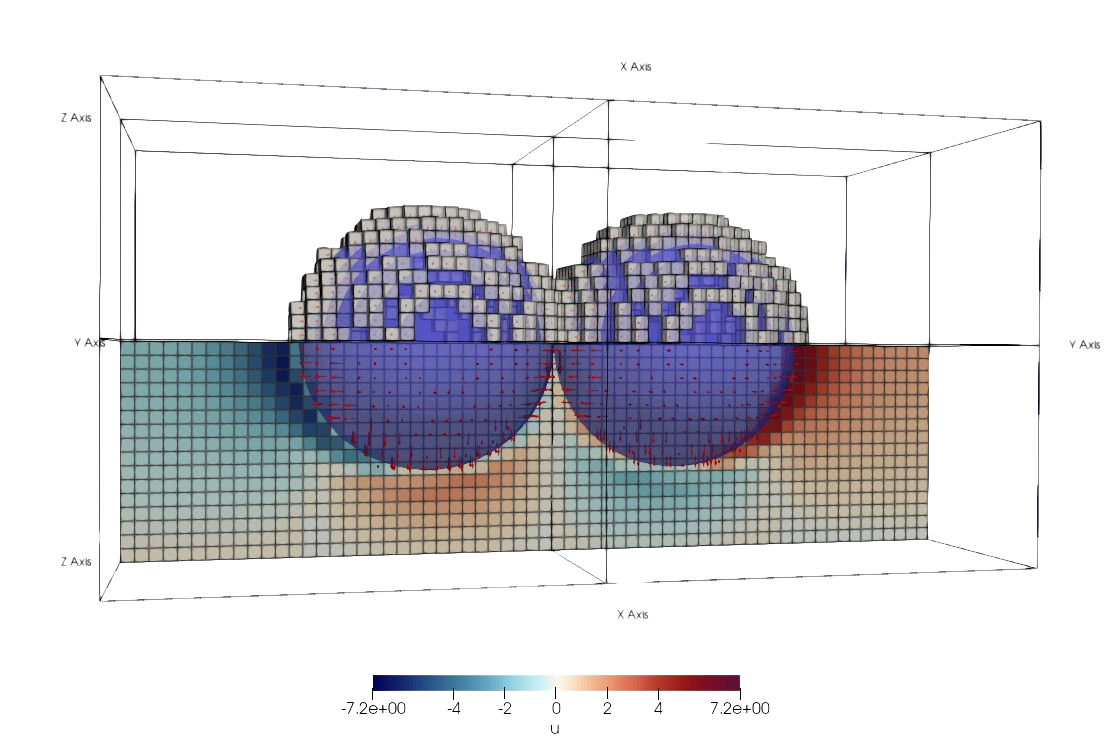}

		\caption{Instantaneous visualization of the 3D colliding spheres test case with $p=1$ and $\alpha=0.3$ just prior to the collision at $t \approx T/6$. The simulation's black outline depicts the individual processor boundaries, symmetrically dividing each axis into two. The background mesh is shown partially with the $x$-component of the velocity resulting from the rotation. The agglomeration groups with cut cells are shown with red glyphs indicating the agglomeration mapping from source to target cells in $\mathfrak{A}$, while the interface between subdomains $\mathfrak{A}$ and $\mathfrak{B}$ is highlighted with blue color. \label{fig:CollidingSpheresInstant}}

\end{figure}

To simulate the colliding sphere test case, a rectangular computational domain with $L_x=2L_y=2L_z$ is discretized into a uniform grid of $64 \times 32$ cells in 2D and $64 \times 32 \times 32$ cells in 3D. The spheres, each with a radius $r_{\mathrm{s}} = 0.15 L_x$, are positioned at $\pm 1.5 r_{\mathrm{s}}$ relative to the origin in $x$-axis. They are assigned with velocities directed toward each other, resulting in $Re=   \rho  u_\mathrm{s} 2 r_{\mathrm{s}} / \mu=1000$. The simulation period is determined such that the spheres interchange their initial positions at the end of the simulations, i.e. $T=3 r_{\mathrm{s}} / u_{\mathrm{s}}$ and discretized into 100 time steps with an implicit Euler time scheme.
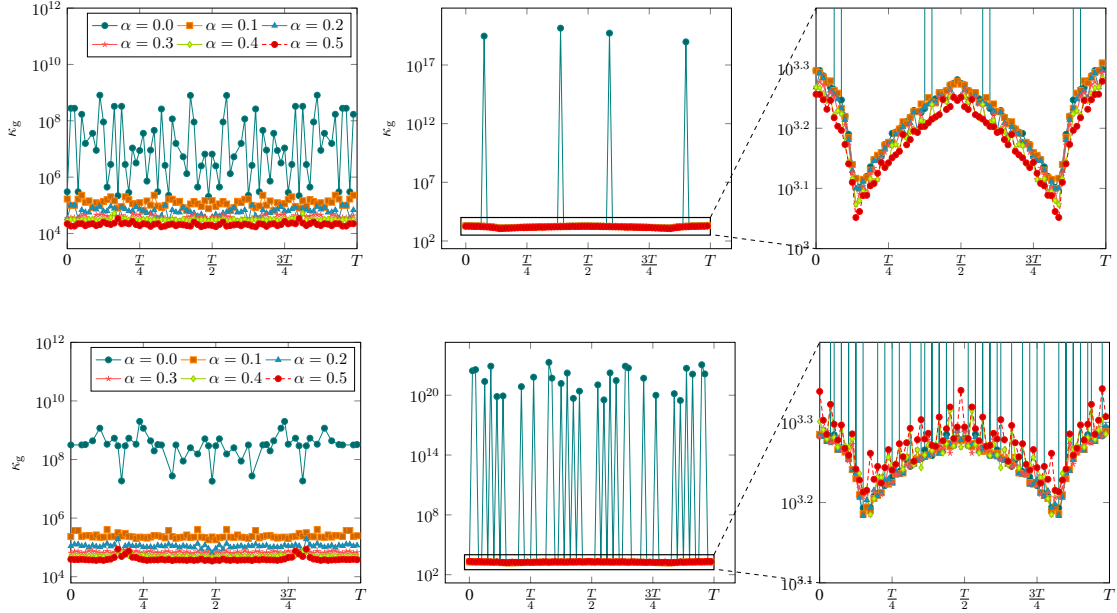
\begin{figure}[tpb]
	\centering

	\begin{subfigure}[b]{0.9\textwidth}
		\centering
		\resizebox{\linewidth}{!}{\begin{tikzpicture}

    \begin{axis}[
        name=ax0,
        ylabel={$\kappa_\mathrm{g}$},
        xmin=0, xmax=2*pi,
        ymax=10^12,
        ymode=log,
        cycle list name=exotic,
        grid=none,
        legend columns=3,
        xtick={0, pi/2, pi, 3*pi/2, 2*pi}, 
        xticklabels={0, $\frac{T}{4}$, $\frac{T}{2}$, $\frac{3T}{4}$, $T$}, 
    ]
    \addplot table[col sep=comma, /pgf/number format/use comma, x expr=\thisrow{File Index}/80*2*pi] {data/CollidingSpheres2D_k1_Agg0/MassMtxCondNo_data.csv};
    \addlegendentry{$\alpha = 0.0$};

   \addplot table[col sep=comma, /pgf/number format/use comma, x expr=\thisrow{File Index}/80*2*pi] {data/CollidingSpheres2D_k1_Agg0.1/MassMtxCondNo_data.csv}; \addlegendentry{$\alpha = 0.1$};

   \addplot table[col sep=comma, /pgf/number format/use comma, x expr=\thisrow{File Index}/80*2*pi] {data/CollidingSpheres2D_k1_Agg0.2/MassMtxCondNo_data.csv};
   \addlegendentry{$\alpha = 0.2$};

   \addplot table[col sep=comma, /pgf/number format/use comma, x expr=\thisrow{File Index}/80*2*pi] {data/CollidingSpheres2D_k1_Agg0.3/MassMtxCondNo_data.csv};
   \addlegendentry{$\alpha = 0.3$};

   \addplot table[col sep=comma, /pgf/number format/use comma, x expr=\thisrow{File Index}/80*2*pi] {data/CollidingSpheres2D_k1_Agg0.4/MassMtxCondNo_data.csv};
   \addlegendentry{$\alpha = 0.4$};

   \addplot table[col sep=comma, /pgf/number format/use comma, x expr=\thisrow{File Index}/80*2*pi] {data/CollidingSpheres2D_k1_Agg0.5/MassMtxCondNo_data.csv};
   \addlegendentry{$\alpha = 0.5$};

    \end{axis}

    \begin{axis}[
        name=ax1,
        ylabel={$\kappa_\mathrm{g}$},
        ymode=log,
        cycle list name=exotic,
        grid=none,
        legend columns=3,
        legend pos=north east,
        at={($(ax0.south east)+(2cm,0)$)},
        xtick={0, pi/2, pi, 3*pi/2, 2*pi}, 
        xticklabels={0, $\frac{T}{4}$, $\frac{T}{2}$, $\frac{3T}{4}$, $T$}, 
    ]

   \addplot table[col sep=comma, /pgf/number format/use comma, x expr=\thisrow{File Index}/80*2*pi] {data/CollidingSpheres2D_k1_Agg0/TotCondNo-Vars0.1.2_data.csv};

   \addplot table[col sep=comma, /pgf/number format/use comma, x expr=\thisrow{File Index}/80*2*pi] {data/CollidingSpheres2D_k1_Agg0.1/TotCondNo-Vars0.1.2_data.csv};

   \addplot table[col sep=comma, /pgf/number format/use comma, x expr=\thisrow{File Index}/80*2*pi] {data/CollidingSpheres2D_k1_Agg0.2/TotCondNo-Vars0.1.2_data.csv};

   \addplot table[col sep=comma, /pgf/number format/use comma, x expr=\thisrow{File Index}/80*2*pi] {data/CollidingSpheres2D_k1_Agg0.3/TotCondNo-Vars0.1.2_data.csv};

   \addplot table[col sep=comma, /pgf/number format/use comma, x expr=\thisrow{File Index}/80*2*pi] {data/CollidingSpheres2D_k1_Agg0.4/TotCondNo-Vars0.1.2_data.csv};

   \addplot table[col sep=comma, /pgf/number format/use comma, x expr=\thisrow{File Index}/80*2*pi] {data/CollidingSpheres2D_k1_Agg0.5/TotCondNo-Vars0.1.2_data.csv};

  \coordinate (c1) at (axis cs:-pi/25,10^2.5);
  \coordinate (c2) at (axis cs:2.04*pi,10^2.5);
  \coordinate (c3) at (axis cs:2*pi,10^4);    

  \draw[thick] (c1) rectangle (c3);

    \end{axis}

    \begin{axis}[
        name=ax2,
        ymin=10^3, ymax=10^3.4,
        xmin=0, xmax=pi*2,
        ymode=log,
        cycle list name=exotic,
        grid=none,
        at={($(ax1.south east)+(2cm,0)$)},
        label style={color=black,font=\bfseries},
        xtick={0, pi/2, pi, 3*pi/2, 2*pi}, 
        xticklabels={0,$\frac{T}{4}$, $\frac{T}{2}$,  $\frac{3T}{4}$, $T$}, 
    ]
    \addplot table[col sep=comma, /pgf/number format/use comma, x expr=\thisrow{File Index}/80*2*pi] {data/CollidingSpheres2D_k1_Agg0/TotCondNo-Vars0.1.2_data.csv};

   \addplot table[col sep=comma, /pgf/number format/use comma, x expr=\thisrow{File Index}/80*2*pi] {data/CollidingSpheres2D_k1_Agg0.1/TotCondNo-Vars0.1.2_data.csv};

   \addplot table[col sep=comma, /pgf/number format/use comma, x expr=\thisrow{File Index}/80*2*pi] {data/CollidingSpheres2D_k1_Agg0.2/TotCondNo-Vars0.1.2_data.csv};

   \addplot table[col sep=comma, /pgf/number format/use comma, x expr=\thisrow{File Index}/80*2*pi] {data/CollidingSpheres2D_k1_Agg0.3/TotCondNo-Vars0.1.2_data.csv};

   \addplot table[col sep=comma, /pgf/number format/use comma, x expr=\thisrow{File Index}/80*2*pi] {data/CollidingSpheres2D_k1_Agg0.4/TotCondNo-Vars0.1.2_data.csv};

   \addplot table[col sep=comma, /pgf/number format/use comma, x expr=\thisrow{File Index}/80*2*pi] {data/CollidingSpheres2D_k1_Agg0.5/TotCondNo-Vars0.1.2_data.csv};
\end{axis}
\draw [dashed] (c2) -- (ax2.south west);
\draw [dashed] (c3) -- (ax2.north west);
\end{tikzpicture}}
		\label{fig:CollidingSpheres2DOpCondk1}
	\end{subfigure}
	\medskip	

	\centering
	\begin{subfigure}[b]{0.9\textwidth}
		\centering
		\resizebox{\linewidth}{!}{\begin{tikzpicture}
    \begin{axis}[
        name=ax0,
        ylabel={$\kappa_\mathrm{g}$},
        xmin=0, xmax=2*pi,
        ymax=10^12,
        ymode=log,
        cycle list name=exotic,
        grid=none,
        legend columns=3, 
        xtick={0, pi/2, pi, 3*pi/2, 2*pi}, 
        xticklabels={0, $\frac{T}{4}$, $\frac{T}{2}$, $\frac{3T}{4}$, $T$}, 
    ]
    \addplot table[col sep=comma, /pgf/number format/use comma, x expr=\thisrow{File Index}/80*2*pi] {data/CollidingSpheres3D_k1_Agg0/MassMtxCondNo_data.csv};
    \addlegendentry{$\alpha = 0.0$};

   \addplot table[col sep=comma, /pgf/number format/use comma, x expr=\thisrow{File Index}/80*2*pi] {data/CollidingSpheres3D_k1_Agg0.1/MassMtxCondNo_data.csv}; \addlegendentry{$\alpha = 0.1$};

   \addplot table[col sep=comma, /pgf/number format/use comma, x expr=\thisrow{File Index}/80*2*pi] {data/CollidingSpheres3D_k1_Agg0.2/MassMtxCondNo_data.csv};
   \addlegendentry{$\alpha = 0.2$};

   \addplot table[col sep=comma, /pgf/number format/use comma, x expr=\thisrow{File Index}/80*2*pi] {data/CollidingSpheres3D_k1_Agg0.3/MassMtxCondNo_data.csv};
   \addlegendentry{$\alpha = 0.3$};

   \addplot table[col sep=comma, /pgf/number format/use comma, x expr=\thisrow{File Index}/80*2*pi] {data/CollidingSpheres3D_k1_Agg0.4/MassMtxCondNo_data.csv};
   \addlegendentry{$\alpha = 0.4$};

   \addplot table[col sep=comma, /pgf/number format/use comma, x expr=\thisrow{File Index}/80*2*pi] {data/CollidingSpheres3D_k1_Agg0.5/MassMtxCondNo_data.csv};
   \addlegendentry{$\alpha = 0.5$};

    \end{axis}

    \begin{axis}[
        name=ax1,
        ymode=log,
        cycle list name=exotic,
        grid=none,
        legend columns=3,
        legend pos=north east,
        at={($(ax0.south east)+(2cm,0)$)},
        xtick={0, pi/2, pi, 3*pi/2, 2*pi}, 
        xticklabels={0, $\frac{T}{4}$, $\frac{T}{2}$, $\frac{3T}{4}$, $T$}, 
    ]

   \addplot table[col sep=comma, /pgf/number format/use comma, x expr=\thisrow{File Index}/80*2*pi] {data/CollidingSpheres3D_k1_Agg0/TotCondNo-Vars0.1.2.3_data.csv};

   \addplot table[col sep=comma, /pgf/number format/use comma, x expr=\thisrow{File Index}/80*2*pi] {data/CollidingSpheres3D_k1_Agg0.1/TotCondNo-Vars0.1.2.3_data.csv};

   \addplot table[col sep=comma, /pgf/number format/use comma, x expr=\thisrow{File Index}/80*2*pi] {data/CollidingSpheres3D_k1_Agg0.2/TotCondNo-Vars0.1.2.3_data.csv};

   \addplot table[col sep=comma, /pgf/number format/use comma, x expr=\thisrow{File Index}/80*2*pi] {data/CollidingSpheres3D_k1_Agg0.3/TotCondNo-Vars0.1.2.3_data.csv};

   \addplot table[col sep=comma, /pgf/number format/use comma, x expr=\thisrow{File Index}/80*2*pi] {data/CollidingSpheres3D_k1_Agg0.4/TotCondNo-Vars0.1.2.3_data.csv};

   \addplot table[col sep=comma, /pgf/number format/use comma, x expr=\thisrow{File Index}/80*2*pi] {data/CollidingSpheres3D_k1_Agg0.5/TotCondNo-Vars0.1.2.3_data.csv};

  \coordinate (c1) at (axis cs:-pi/25,10^2.5);
  \coordinate (c2) at (axis cs:2.04*pi,10^2.5);
  \coordinate (c3) at (axis cs:2*pi,10^4);    

  \draw[thick] (c1) rectangle (c3);

    \end{axis}

    \begin{axis}[
        name=ax2,
        ymin=10^3.1, ymax=10^3.4,
        xmin=0, xmax=pi*2,
        ymode=log,
        cycle list name=exotic,
        grid=none,
        at={($(ax1.south east)+(2cm,0)$)},
        label style={color=black,font=\bfseries},
        xtick={0, pi/2, pi, 3*pi/2, 2*pi}, 
        xticklabels={0,$\frac{T}{4}$, $\frac{T}{2}$,  $\frac{3T}{4}$, $T$}, 
    ]
    \addplot table[col sep=comma, /pgf/number format/use comma, x expr=\thisrow{File Index}/80*2*pi] {data/CollidingSpheres3D_k1_Agg0/TotCondNo-Vars0.1.2.3_data.csv};

   \addplot table[col sep=comma, /pgf/number format/use comma, x expr=\thisrow{File Index}/80*2*pi] {data/CollidingSpheres3D_k1_Agg0.1/TotCondNo-Vars0.1.2.3_data.csv};

   \addplot table[col sep=comma, /pgf/number format/use comma, x expr=\thisrow{File Index}/80*2*pi] {data/CollidingSpheres3D_k1_Agg0.2/TotCondNo-Vars0.1.2.3_data.csv};

   \addplot table[col sep=comma, /pgf/number format/use comma, x expr=\thisrow{File Index}/80*2*pi] {data/CollidingSpheres3D_k1_Agg0.3/TotCondNo-Vars0.1.2.3_data.csv};

   \addplot table[col sep=comma, /pgf/number format/use comma, x expr=\thisrow{File Index}/80*2*pi] {data/CollidingSpheres3D_k1_Agg0.4/TotCondNo-Vars0.1.2.3_data.csv};

   \addplot table[col sep=comma, /pgf/number format/use comma, x expr=\thisrow{File Index}/80*2*pi] {data/CollidingSpheres3D_k1_Agg0.5/TotCondNo-Vars0.1.2.3_data.csv};
\end{axis}
\draw [dashed] (c2) -- (ax2.south west);
\draw [dashed] (c3) -- (ax2.north west);
\end{tikzpicture}}
		\label{fig:CollidingSpheres3DOpCondk1}
	\end{subfigure}
	\caption{Global condition numbers for the mass (left) and operator matrices (middle, right) in the 2D (top) and 3D (bottom) colliding spheres test cases with polynomial degree $p=1$ and varying agglomeration threshold $\alpha$. Note that `Inf' entries computed for the mass matrix, which correspond to the jumps at operator matrices, are excluded in the figure. \label{fig:CollidingSpheresCondk1}}
 \end{figure}

 Figure~\ref{fig:CollidingSpheresInstant} illustrates an example of 3D simulations executed on eight processors, at the instance just before the collision. The global condition numbers calculated for the colliding sphere test cases with $p=1$ are presented in Figure~\ref{fig:CollidingSpheresCondk1}. For both the 2D and 3D setups, it is evident that cases with $\alpha=0.0$ exhibit recurring jumps during the simulations, reaching values 15 orders of magnitude larger than the rest due to the small-cut problem. Moreover, `Inf' entries emerge in the MATLAB calculations for the mass matrices corresponding to these jumps, mirroring patterns observed in other test cases. This behavior becomes more noticeable in the 3D case due to its larger mass structure. However, the application of agglomeration significantly reduces the condition numbers, demonstrating its efficiency in addressing the small-cut problem. Overall, the increasing agglomeration thresholds lead to decreasing condition numbers for mass matrices, which is consistent with observations in the other scenarios. For the operator matrices, a particular trend is not observed with agglomeration threshold $\alpha$. A summary of the calculated condition numbers is provided in Table~\ref{tab:CollidingResults}.
  \begin{table}[htb!] 
    \centering
    \caption{Maximum condition numbers for the colliding spheres test case}
    \label{tab:CollidingResults}
	\small
    \begin{tabular}{l|c:c|c:c|c:c}  
    \toprule
		\multicolumn{1}{c}{} & \multicolumn{6}{c}{2D global numbers($\kappa_\mathrm{g}$)}  \\ 
  \cmidrule(rl){2-7}  
\multicolumn{1}{c}{} & \multicolumn{2}{c}{$p=1$}  & \multicolumn{2}{c}{$p=2$} & \multicolumn{2}{c}{$p=3$} \\ 
        $\alpha$ & Mass & Op  & Mass & Op & Mass & Op \\ \hline
        0\phantom{.1}  \; & 8.10e8* & 1.38e20 & 1.94e9* & 1.70e20* & 6.09e8* & 1.68e22\\ \hline
		0.1 \; & 2.28e5 & 2.04e3 & 2.47e6 & 3.93e3 & 1.22e7 & 7.81e3 \\ \hline
        0.2 \; & 6.96e4 & 1.98e3 & 1.05e6 & 2.52e3 & 4.87e6 & 6.36e3 \\ \hline
        0.3 \; & 5.15e4 & 1.96e3 & 5.67e5 & 1.96e3 & 3.06e6 & 4.56e3 \\ \hline
        0.4 \; & 5.15e4 & 1.92e3 & 5.67e5 & 1.95e3 & 3.06e6 & 4.06e3 \\ \hline
        0.5 \; & 3.38e4 & 1.90e3 & 3.35e5 & 1.95e3 & 1.63e6 & 3.53e3 \\ \hline
    \end{tabular}
    \begin{tabular}{l|c:c|c:c|c:c}  
		\multicolumn{7}{c}{} \\
		\multicolumn{1}{c}{} & \multicolumn{6}{c}{3D stencil numbers ($\kappa_\mathrm{s}$)}  \\ 
  \cmidrule(rl){2-7} 
        0\phantom{.1}  \; & 2.74e10 & 1.40e2 & 3.47e10* & 1.43e23 & 1.07e10* & 1.11e23*\\ \hline
		0.1 \; & 3.58e5 & 8.28e1 & 4.31e6 & 4.09e2 & 2.45e7 & 1.47e3 \\ \hline
        0.2 \; & 1.91e5 & 7.47e1 & 2.43e6 & 3.28e2 & 1.39e7 & 1.22e3 \\ \hline
        0.3 \; & 8.09e4 & 6.15e1 & 9.66e5 & 2.87e2 & 5.60e6 & 1.10e3 \\ \hline
        0.4 \; & 7.55e4 & 6.15e1 & 9.48e5 & 2.65e2 & 5.66e6 & 9.94e2 \\ \hline
        0.5 \; & 8.69e4 & 7.36e1 & 9.55e5 & 2.55e2 & 5.95e6 & 8.41e2 \\ \hline
    \end{tabular}
    \begin{tabular}{l}
		\footnotesize *Maximum number excluding `Inf' entries      \phantom{afasdfadfadafasfsfaadfafasdfasdfasdf}
	\end{tabular}
\end{table}

\subsection{Rotating popcorn and torus}
\begin{figure}[ht!]
	\centering
	\includegraphics[width=0.8\textwidth]{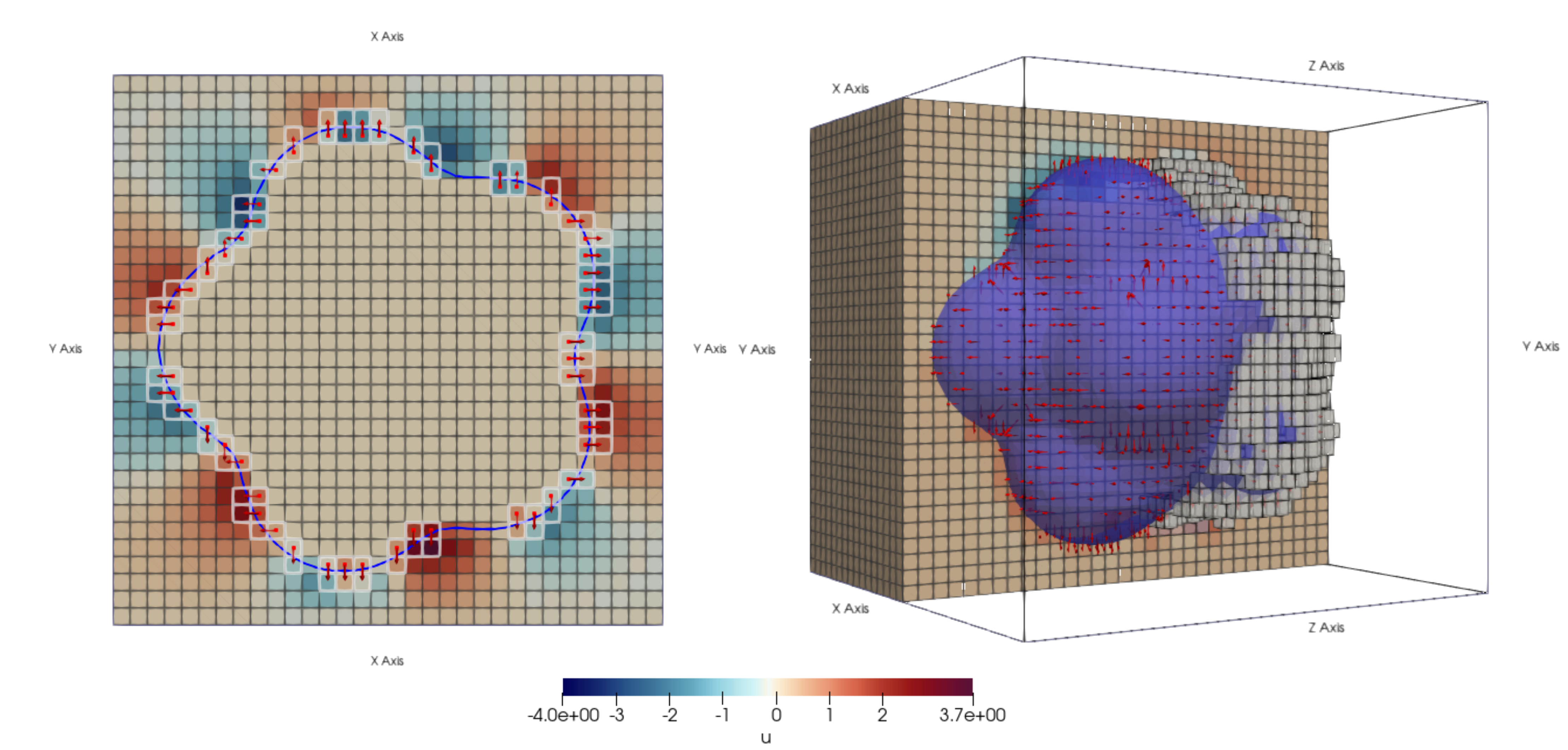}

	\caption{Instantaneous visualization of the 2D (left) and 3D (right) rotating popcorn test case with $p=3$ and $\alpha=0.5$ at the third time step resolved on domains with a resolution of 32 cells in each axes. The embedded shapes are displayed on the background mesh with the $x$-component of the velocity resulting from the rotation. The interface between subdomains $\mathfrak{A}$ and $\mathfrak{B}$ is highlighted with blue color. The red glyphs illustrate the agglomeration mappings from source to target cells in $\mathfrak{A}$.  \label{fig:PopcornInstant}}
 \end{figure}
 \begin{figure}[ht!]
	\centering
	\includegraphics[width=0.8\textwidth]{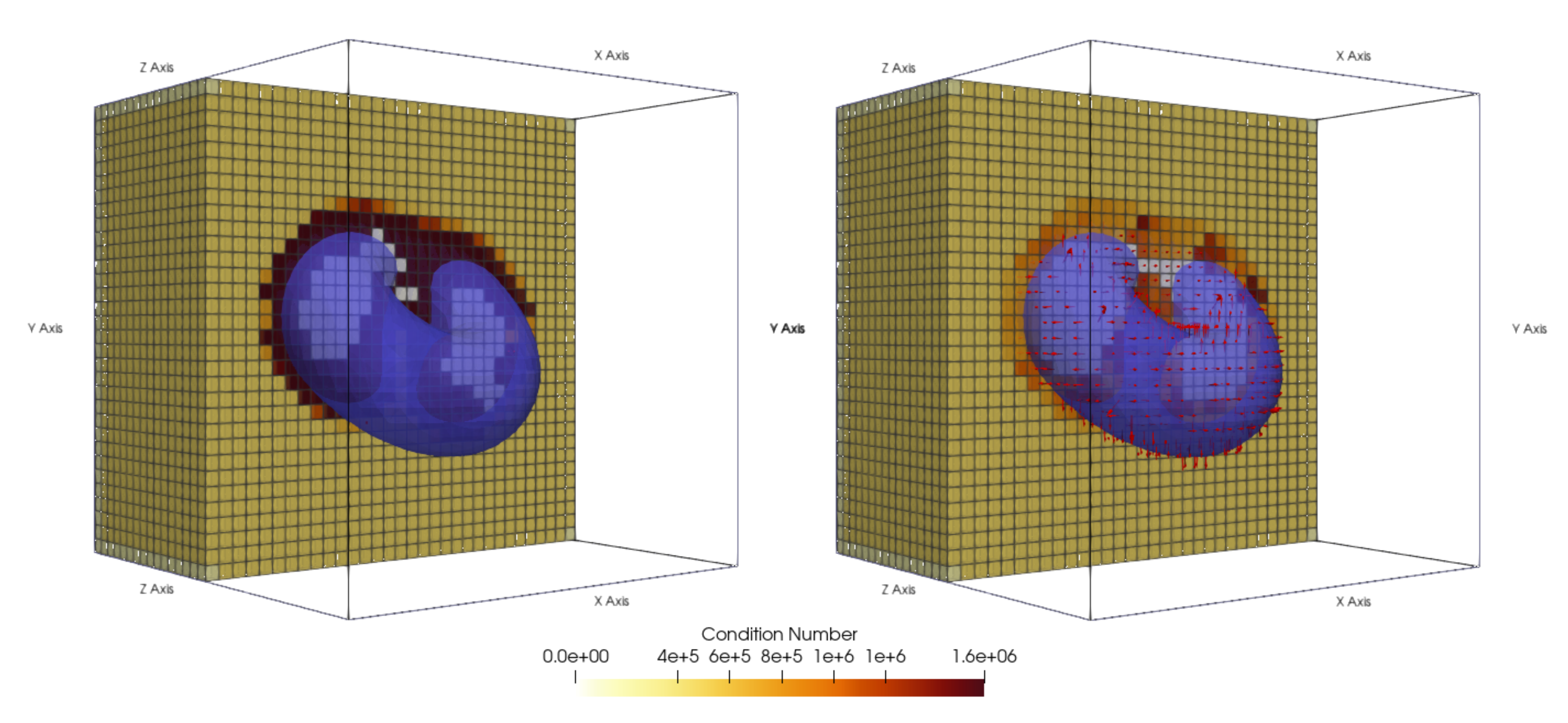}

	\caption{Instantaneous visualization of the 3D rotating torus test case with $32\times32\times32$ cells, $p=3$, $\alpha=0$ (left) and $\alpha=0.5$ (right) at the first time step. The embedded shapes are displayed on the background mesh with the stencil condition number $\kappa_\mathrm{s}$. The interface between subdomains $\mathfrak{A}$ and $\mathfrak{B}$ is highlighted with blue color. The red glyphs illustrate the agglomeration mappings from source to target cells in $\mathfrak{A}$.  \label{fig:RotTorusInstant}}
 \end{figure}
2D and 3D popcorn-like shapes are described in square and cubic computational domains with a dimension of $L_x$ as illustrated in Figure~\ref{fig:PopcornInstant} and rotated in $z$-axis. The core radius of the shape is determined as $r_\mathrm{p}=0.6 L_x$. An angular velocity corresponding to Reynold number $(2 \rho U  r_{\mathrm{p}}/ \mu)$ of 1000 is applied to the pseudo-solid phase. The mesh resolution is varied as $16\times16$ and $32\times32$, $64\times64$, and $128\times128$ in 2D and as $16\times16\times16$, $32\times32\times32$ in 3D. The simulations are conducted until the rotated shapes reach their symmetry axis (i.e., $T=2 \pi / |5\Vec{\Omega}|$) and discretized into 80 time steps.

Analogously, the 3D rotating torus is placed in a cubic domain with $L_x$ as illustrated in Figure~\ref{fig:RotTorusInstant} and rotated in $y$-axis with an angular velocity corresponding to $Re= 2 \rho r_{\mathrm{ma}} U / \mu = 1000$. The major and minor radiuses are defined as $r_{\mathrm{ma}}=0.39L_x$ and $r_{\mathrm{mi}}=0.26L_x$, respectively. The mesh resolution is varied as $16\times16\times16$ and $32\times32\times32$. The simulation period is determined as $T= \pi  / |4\Vec{\Omega}|$ and discretized into 25 time steps. 


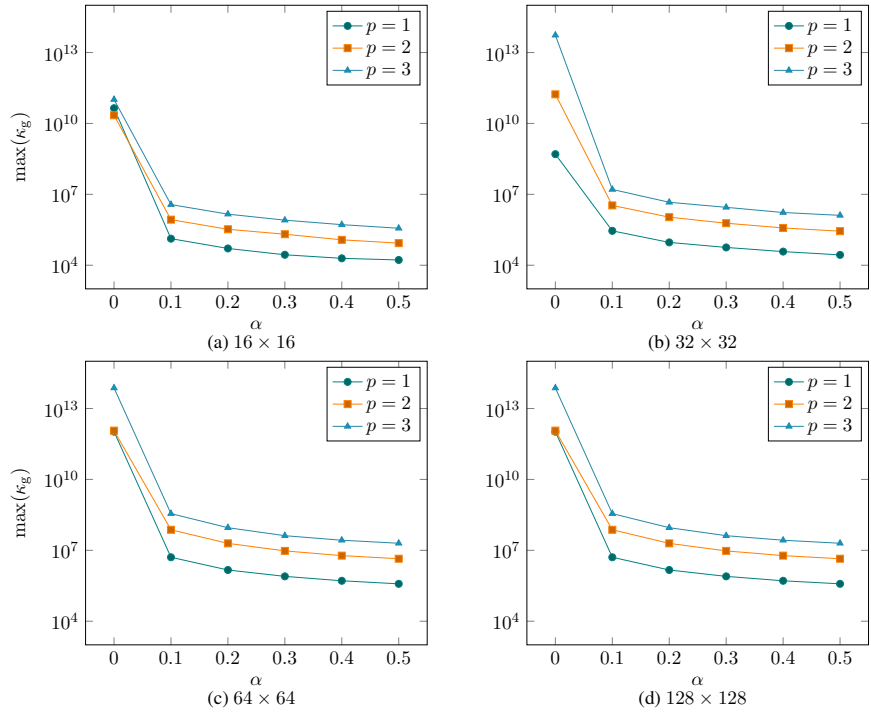
\begin{figure}[htb!]
	\centering
	\begin{subfigure}[b]{0.7\textwidth}
		\centering
		\resizebox{\linewidth}{!}{\begin{tikzpicture}
    \begin{axis}[
        name=ax0,
        xlabel={$\alpha$},
        ylabel={$\max(\kappa_\mathrm{g})$},
        ymin=1e3, ymax=1e15,
        ymode=log,
        cycle list name=exotic,
        grid=none,
        legend columns=1,
    ]
    \addplot coordinates {
        (0,4.41e+10)   (0.1,1.33e+05)  (0.2,5.13e+04)  (0.3,2.77e+04)
        (0.4,1.96e+04) (0.5,1.67e+04)
    };      \addlegendentry{$p = 1$};

    \addplot coordinates {
        (0,2.21e+10)   (0.1,8.49e+05)  (0.2,3.32e+05)  (0.3,2.05e+05)
        (0.4,1.19e+05) (0.5,8.61e+04)
    };      \addlegendentry{$p = 2$};

    \addplot coordinates {
        (0,1.03e+11)   (0.1,3.68e+06)  (0.2,1.45e+06)  (0.3,8.07e+05)
        (0.4,5.21e+05) (0.5,3.67e+05)
    };      \addlegendentry{$p = 3$};   

    \end{axis}

    \node[align=center] at ($(ax0.south)+(-0.1,-1.1cm)$) {(a) $16 \times 16$};

    \begin{axis}[
        name=ax1,
        xlabel={$\alpha$},
        at={($(ax0.south east)+(2cm,0)$)},
        ymin=1e3, ymax=1e15,
        ymode=log,
        cycle list name=exotic,
        grid=none,
        legend columns=1
    ]
    \addplot coordinates {
            (0,5.01e+08)   (0.1,2.85e+05)  (0.2,9.25e+04)  (0.3,5.66e+04)
            (0.4,3.78e+04) (0.5,2.74e+04)
    };      \addlegendentry{$p = 1$};
    
    \addplot coordinates {
            (0,1.70e+11)  
             (0.1,3.40e+06)  (0.2,1.08e+06)  (0.3,6.00e+05)
            (0.4,3.79e+05) (0.5,2.76e+05)
    };      \addlegendentry{$p = 2$};
    
    \addplot coordinates {
            (0,5.42e+13) 
              (0.1,1.61e+07)  (0.2,4.61e+06)  (0.3,2.82e+06)
            (0.4,1.70e+06) (0.5,1.29e+06)
    };      \addlegendentry{$p = 3$};   

    \end{axis}
    \node[align=center] at ($(ax1.south)+(-0.1,-1.1cm)$) {(b) $32 \times 32$};
\end{tikzpicture}}
	\end{subfigure}
	\medskip
	\centering
	\begin{subfigure}[b]{0.7\textwidth}
		\centering
		\resizebox{\linewidth}{!}{\begin{tikzpicture}
    \begin{axis}[
        name=ax0,
        xlabel={$\alpha$},
        ylabel={$\max(\kappa_\mathrm{g})$},
        ymin=1e3, ymax=1e15,
        ymode=log,
        cycle list name=exotic,
        grid=none,
        legend columns=1,
    ]
    \addplot coordinates {
        (0,1.05e+12)   (0.1,5.08e+06)  (0.2,1.47e+06)  (0.3,7.90e+05)
        (0.4,5.14e+05) (0.5,3.82e+05)
    };      \addlegendentry{$p = 1$};

    \addplot coordinates {
        (0,1.16e+12)   (0.1,7.41e+07)  (0.2,1.97e+07)  (0.3,9.45e+06)
        (0.4,5.96e+06) (0.5,4.41e+06)
    };      \addlegendentry{$p = 2$};

    \addplot coordinates {
        (0,7.40e+13)   (0.1,3.57e+08)  (0.2,9.07e+07)  (0.3,4.16e+07)
        (0.4,2.67e+07) (0.5,1.99e+07)
    };      \addlegendentry{$p = 3$};       \addlegendentry{$p = 3$};   

    \end{axis}

    \node[align=center] at ($(ax0.south)+(-0.1,-1.1cm)$) {(c) $64 \times 64$};    
    
    \begin{axis}[
        name=ax1,
        xlabel={$\alpha$},
        at={($(ax0.south east)+(2cm,0)$)},
        ymin=1e3, ymax=1e15,
        ymode=log,
        cycle list name=exotic,
        grid=none,
        legend columns=1
    ]
        \addplot coordinates {
            (0,1.05e+12)   (0.1,5.08e+06)  (0.2,1.47e+06)  (0.3,7.90e+05)
            (0.4,5.14e+05) (0.5,3.82e+05)
        };      \addlegendentry{$p = 1$};
    
        \addplot coordinates {
            (0,1.16e+12)   (0.1,7.41e+07)  (0.2,1.97e+07)  (0.3,9.45e+06)
            (0.4,5.96e+06) (0.5,4.41e+06)
        };      \addlegendentry{$p = 2$};
    
        \addplot coordinates {
            (0,7.40e+13)   (0.1,3.57e+08)  (0.2,9.07e+07)  (0.3,4.16e+07)
            (0.4,2.67e+07) (0.5,1.99e+07)
        };      \addlegendentry{$p = 3$};   
    \end{axis}

    \node[align=center] at ($(ax1.south)+(-0.1,-1.1cm)$) {(d) $128 \times 128$};

\end{tikzpicture}}
	\end{subfigure}
	\caption{Global condition numbers ($\kappa_\mathrm{s}$) for the mass matrices in the 2D rotating popcorn test case with varying grid resolution, polynomial degree $p$, and agglomeration threshold $\alpha$. The infinite condition numbers originating from extremely small-cuts in simulations with $\alpha=0$ are excluded. \label{fig:Popcorn2DCondNumbers}}
\end{figure}

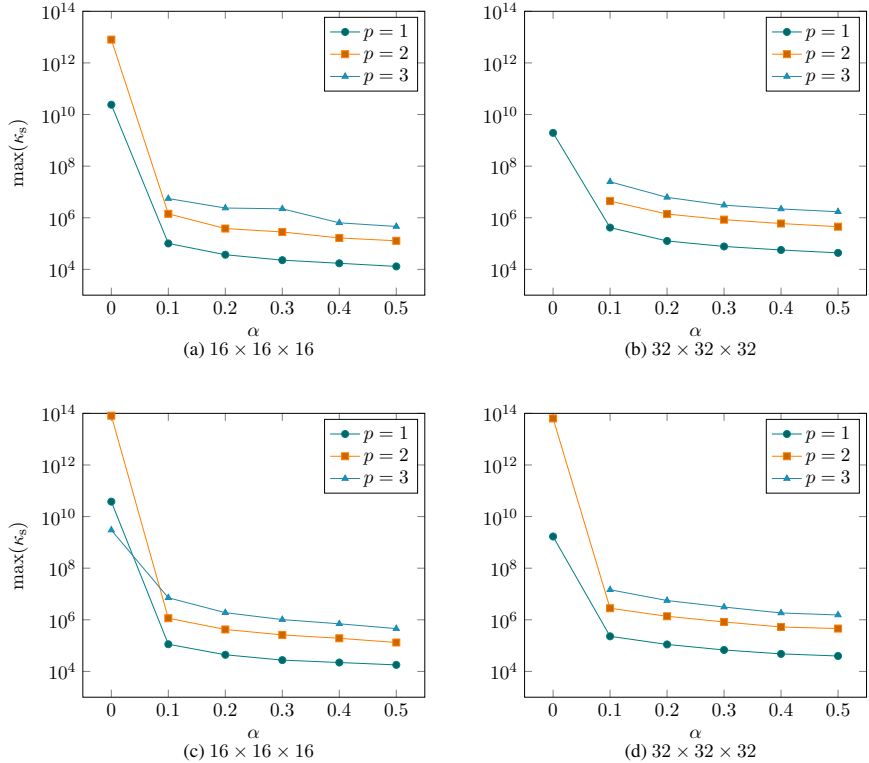
\begin{figure}[htb!]
	\centering
	\begin{subfigure}[b]{0.7\textwidth}
		\centering
		\resizebox{\linewidth}{!}{\begin{tikzpicture}
    \begin{axis}[
        name=ax0,
        xlabel={$\alpha$},
        ylabel={$\max(\kappa_\mathrm{s})$},
        ymin=1e3, ymax=1e14,
        ymode=log,
        cycle list name=exotic,
        grid=none,
        legend columns=1,
    ]
    \addplot coordinates {
        (0,2.40e+10)   (0.1,1.01e+05)  (0.2,3.69e+04)  (0.3,2.27e+04)
        (0.4,1.72e+04) (0.5,1.30e+04)
    };      \addlegendentry{$p = 1$};

    \addplot coordinates {
        (0,8.00e+12)   (0.1,1.42e+06)  (0.2,3.82e+05)  (0.3,2.81e+05)
        (0.4,1.65e+05) (0.5,1.28e+05)
    };      \addlegendentry{$p = 2$};

    \addplot coordinates {
        (0.1,5.55e+06)  (0.2,2.38e+06)  (0.3,2.22e+06)
        (0.4,6.36e+05) (0.5,4.59e+05)
    };       \addlegendentry{$p = 3$};   

    \end{axis}

    \node[align=center] at ($(ax0.south)+(-0.1,-1.1cm)$) {(a) $16 \times 16 \times 16$};

    \begin{axis}[
        name=ax1,
        xlabel={$\alpha$},
        at={($(ax0.south east)+(2cm,0)$)},
        ymin=1e3, ymax=1e14,
        ymode=log,
        cycle list name=exotic,
        grid=none,
        legend columns=1
    ]
    \addplot coordinates {
        (0,1.94e+09)   (0.1,4.16e+05)  (0.2,1.26e+05)  (0.3,7.71e+04)
        (0.4,5.64e+04) (0.5,4.35e+04)
    };      \addlegendentry{$p = 1$};

    \addplot coordinates {
        (0.1,4.42e+06)  (0.2,1.39e+06)  (0.3,8.44e+05)
        (0.4,5.97e+05) (0.5,4.50e+05)
    };      \addlegendentry{$p = 2$};

    \addplot coordinates {
        (0.1,2.47e+07)  (0.2,6.17e+06)  (0.3,3.07e+06)
        (0.4,2.19e+06) (0.5,1.71e+06)
    };      \addlegendentry{$p = 3$};   

    \end{axis}
    \node[align=center] at ($(ax1.south)+(-0.1,-1.1cm)$) {(b) $32 \times 32 \times 32$};
\end{tikzpicture}}
		 \label{fig:Popcorn3DCondNumbers}
	\end{subfigure}

	\centering
	\begin{subfigure}[b]{0.7\textwidth}
		\centering
		\resizebox{\linewidth}{!}{\begin{tikzpicture}
    \begin{axis}[
        name=ax0,
        xlabel={$\alpha$},
        ylabel={$\max(\kappa_\mathrm{s})$},
        ymin=1e3, ymax=1e14,
        ymode=log,
        cycle list name=exotic,
        grid=none,
        legend columns=1,
    ]
    \addplot coordinates {
        (0,3.79e+10)   (0.1,1.13e+05)  (0.2,4.41e+04)  (0.3,2.74e+04)
        (0.4,2.21e+04) (0.5,1.79e+04)
    };      \addlegendentry{$p = 1$};

    \addplot coordinates {
        (0,8.12e+13)   (0.1,1.15e+06)  (0.2,4.22e+05)  (0.3,2.60e+05)
        (0.4,1.92e+05) (0.5,1.32e+05)
    };      \addlegendentry{$p = 2$};

    \addplot coordinates {
        (0,3.00e+09)   (0.1,7.14e+06)  (0.2,1.89e+06)  (0.3,1.02e+06)
        (0.4,7.00e+05) (0.5,4.55e+05)
    };       \addlegendentry{$p = 3$};   

    \end{axis}
    \node[align=center] at ($(ax0.south)+(-0.1,-1.1cm)$) {(c) $16 \times 16 \times 16$};
    
    \begin{axis}[
        name=ax1,
        xlabel={$\alpha$},
        at={($(ax0.south east)+(2cm,0)$)},
        ymin=1e3, ymax=1e14,
        ymode=log,
        cycle list name=exotic,
        grid=none,
        legend columns=1
    ]
    \addplot coordinates {
        (0,1.69e+09)   (0.1,2.29e+05)  (0.2,1.11e+05)  (0.3,6.77e+04)
        (0.4,4.78e+04) (0.5,3.96e+04)
    };      \addlegendentry{$p = 1$};

    \addplot coordinates {
        (0,6.36e+13)  
        (0.1,2.80e+06) (0.2,1.37e+06) (0.3,8.25e+05)  (0.4,5.26e+05)
        (0.5,4.58e+05) 
    };      \addlegendentry{$p = 2$};

    \addplot coordinates {
        (0.1,1.46e+07)  (0.2,5.59e+06)  (0.3,3.12e+06)
        (0.4,1.84e+06) (0.5,1.54e+06)
    };      \addlegendentry{$p = 3$};   

    \end{axis}
    \node[align=center] at ($(ax1.south)+(-0.1,-1.1cm)$) {(d) $32 \times 32 \times 32$};
\end{tikzpicture}}
	\end{subfigure}

	\caption{Stencil condition numbers ($\kappa_\mathrm{s}$) for the mass matrices in the 3D rotating popcorn (a,b) and torus (c,d) test cases with varying grid resolution, polynomial degree $p$, and agglomeration threshold $\alpha$. The infinite condition numbers originating from extremely small-cuts are excluded. The rotating popcorn with $p=2, 3$ and the torus with $p=3$ lead to complete infinite values with a grid resolution of 32. \label{fig:RotPopcornTorusCondNumbers}}

\end{figure}

The global condition numbers for the 2D rotating popcorn are presented in Figure~\ref{fig:Popcorn2DCondNumbers}, whereas the stencil condition numbers for 3D simulations are presented in Figure~\ref{fig:RotPopcornTorusCondNumbers}. Additionally, Figure~\ref{fig:RotTorusInstant} provides a visual comparison of the stencil condition numbers between $\alpha=0.0$ and $\alpha=0.5$ for an example simulation of the 3D rotating torus test case. As can be seen in the figures, the condition numbers demonstrate a consistent behavior with other test cases. The application of small-cut agglomeration significantly reduces the condition numbers, which otherwise reach extreme values and can even become infinite. For instance, the 3D configurations of the popcorn shape with $32\times32\times32$ cells, $\alpha=0$ and $p=2, 3$, and the torus with the same configuration but with $p=3$ lead to completely infinite condition numbers, providing clear examples of ill-conditioned systems. Subsequently, the incremental increase in agglomeration thresholds leads to only gradual reductions in the condition numbers, displaying a gentle slope across various resolutions and polynomial degrees.

Conversely, finer grids and increased polynomial degrees contribute to higher condition numbers. Doubling the mesh resolution leads to an increase in the number of small-cut cells as well as a decrease in their size. Consequently, the presented condition numbers approximately increase to fourfold in the 2D setup and eightfold in the 3D setups for each mesh refinement. Nevertheless, a consistent pattern in the maximum numbers for the simulations with $\alpha=0$ is not discernible since the infinite values are excluded. A similar dynamic is also observed for the increases in polynomial degrees as the number of DOF grows with higher degrees.

\section{Conclusion and outlook} \label{sec:Conclusion}
In this work, we have presented a cell agglomeration approach for eXtended discontinuous Galerkin methods using a parallel algorithm that can seamlessly deal with cut cells originating from complex geometries, moving boundaries, and topology changes in both 2- and 3-dimensional spaces. Our recipe distinguishes direct and chain agglomeration routines and provides complementary sections to create appropriate agglomeration mappings without leading to cycles.
Furthermore, agglomeration levels are introduced in order to maintain the sequential order of the mathematical operations associated with the pairs in chain and inter-processor agglomerations. It is then demonstrated that the local agglomeration pairs can be substituted with their lower-level equivalents, resulting in reduced computational costs. The proposed algorithm can be implemented using sparse matrix-vector and matrix-matrix operations as well as basic graph theory. Therefore, it does not require mesh manipulation, making it easier to implement.

Additionally, we have exhibited selection criteria that result in efficient and unique mappings, independent of the processor numbers, and implemented our algorithm into the BoSSS software package. Subsequently, the algorithm is successfully tested against parallel simulations of demanding test cases inspired by realistic applications, which are characterized by significant irregularities or topological transformations.
During these test cases, the condition numbers corresponding to varying agglomeration thresholds were computed and compared to each other. The application of cell agglomeration significantly reduced the condition numbers, addressing the small-cut problem. Furthermore, a clear trend between increasing threshold and decreasing condition numbers is observed. 

In our future work, we aim to further enhance the capabilities of BoSSS by introducing advanced routines for memory usage and parallelization in iterative solvers, eventually making it a robust tool for addressing complex multiphase problems with high-performance computing.

\section*{Acknowledgments}
The work of M. Toprak is funded by the Deutsche Forschungsgemeinschaft (DFG, German Research Foundation) – Project-ID 492661287 – TRR 361. The work of M. Rieckmann is funded by the Deutsche Forschungsgemeinschaft
(DFG, German Research Foundation) – Project-ID 265191195 – Collaborative Research Center 1194 (CRC 1194). Additionally, the work of the authors is supported by the Graduate
School CE within the Centre for Computational Engineering at
TU Darmstadt.

\section*{Data availability}
The data presented in this study can be reproduced by using Jupyter notebooks presented at this address\footnote{DOI(TBD)} accompanied with BoSSS software. The notebooks can be found in the folder \textit{/public/examples/AgglomerationTestcases} in the respective address.

\section*{Conflict of interest}

The authors declare no potential conflict of interests.


\bibliographystyle{unsrtnat}
\bibliography{template}  






\end{document}